\documentclass[10pt]{amsart}
\usepackage{amsmath}
\usepackage{amsfonts}
\usepackage{amssymb}

\def \cplane{\mathbb{C}}
\def \halfplane{\mathbb{H}}
\def \integer{\mathbb{Z}}
\def \uni{\mathrm{uni}}

\def \mk{\mathfrak}
\def \aut{\mathrm{Aut}}

\def \cplane{\mathbb C}
\def \dbar{\bar{\partial}_{J,j}}
\def \ei{e^{i\theta}}
\def \inv{^{-1}}

\def \gs{\mathrm{gs}}
\def \halfplane{\mathbb{H}}
\def \integer{\mathbb{Z}}

\def \mkj{\mathfrak{j}}
\def \om{\overline{\mathcal M}}
\def \pert{\mathrm{pert}}
\def \pgl{\mathrm{pgl}}
\def \rank{\mathrm{rank}}
\def \rone{\mathbb R}
\def \rs{Riemann surface}
\def \sing{\mathrm{Sing}}
\def \strata{\mathbb{S}}
\def \uni{\mathrm{uni}}

\def \wtm{\widetilde{\mathcal M}}
\def \xto{\xrightarrow}

\def \mr{\mathrm}

\def \mc{\mathcal}

\def \mk{\mathfrak}
\def \mr{\mathrm}

\def \M{\mathcal M}

\def \aut{\mathrm{Aut}}

\def \cplane{\mathbb C}
\def \dbar{\bar{\partial}_{J,i}}
\def \ei{e^{i\theta}}

\def \gs{\mathrm{gs}}
\def \inv{^{-1}}

\def \mkj{\mathfrak{j}}
\def \om{\overline{\mathcal M}}
\def \pert{\mathrm{pert}}
\def \pgl{\mathrm{pgl}}
\def \rone{\mathbb R}
\def \sing{\mathrm{Sing}}
\def \strata{\mathbb{S}}
\def \uni{\mathrm{uni}}

\def \wtm{\widetilde{\mathcal M}}

\numberwithin{equation}{section}

\newtheorem{theorem}{Theorem}[section]

\newtheorem{claim}[theorem]{Claim}

\newtheorem{corollary}[theorem]{Corollary}
\newtheorem{defn}[theorem]{Definition}
\newtheorem{example}[theorem]{Example}

\newtheorem{lemma}[theorem]{Lemma}
\newtheorem{prop}[theorem]{Proposition}
\newtheorem{remark}[theorem]{Remark}
\newtheorem{assumption}[theorem]{Assumption}

\def \inv{^{-1}}
\def \s{stabilization}
\def \mcb{\mathcal{B}}
\def \mcf{\mathcal{F}}
\def \mc{\mathcal}
\def \v{\vskip 0.1in}
\def \n{\noindent}

\def \mco{\mathcal{O}}
\def \mfk{\mathfrak}

\def \rone{\mathbb{R}}

\begin{document}

\title{Symplectic virtual localization of Gromov-Witten invariants}
\author{Bohui Chen}
\address{Department of Mathematics, Sichuan University,
        Chengdu,610064, China}
\email{bohui@cs.wisc.edu}
\author{An-Min Li}
\address{Department of Mathematics, Sichuan University,
        Chengdu,610064, China}
\email{math$\_$li@yahoo.com.cn}
\date{}
\abstract
We show that  moduli spaces of stable maps admits virtual
orbifold structure. The symplectic version of virtual localization formula
is obtained.
\endabstract

\maketitle

Given a compact closed symplectic manifold $(M^{2n},\omega)$ and
an $\omega$-tames almost complex structure  $J$, one can define
the celebrated Gromov-Witten invariants using the moduli spaces of
$J$-holomorphic curves. Such invariants were first discovered by
Ruan-Tian on monotone  manifolds(\cite{Ru-Ti}), then later defined
on general manifolds independently  by several different groups
Fukaya-Ono(\cite{FO}), Li-Tian(\cite{LT2}), Liu-Tian(\cite{Liu}),
Ruan(\cite{R}) and etc. The break-through tool for their works is
now well-known as {\em virtual techniques}. On the other hand, the
algebraic version of the theory was first given by
Li-Tian(\cite{LT}).

Since the theory of Gromov-Witten invariants is set up, the computation
of invariants has been one of the main issue of this area. One of the main
tools of the computation is the localization technique. If the symplectic
manifold admits a torus action,  the action can be
induced on  the moduli spaces of $J$-holomorphic
curves. Since the invariants are obtained via "integration" on the moduli
spaces, Kontsevich observed that one may apply the Atiyah-Bott
localization formula for computation(\cite{K}). To fulfill
such an idea, we need to combine
the Atiyah-Bott localization formula with virtual techniques. We call
such a combination as the virtual localization.
This has been done for algebraic varieties(\cite{GP}). But such a
formula has not been set-up in symplectic category. Our main goal of this
paper is to prove a virtual localization formula on general symplectic
manifolds. We remark that in this paper the group can act on the symplectic
manifold as
symplectomorphisms other than just  Hamiltonian ones.

The ingredients of proving such a virtual localization formula are:
(1), a modified gluing theory which provides smooth structures
on moduli spaces, (2), virtual manifolds/orbifolds and (equivariant)
integration theory on them.
The abstract theory of virtual manifold/orbifolds has been established
in \cite{CT}.  In this paper, we mainly explain how to obtain
smooth structures on moduli spaces via the gluing theory, and then generalize
it to virtual moduli spaces accordingly.

The paper is organized as following: in Part I, we introduce some
preliminary materials that is needed to understand the moduli spaces;
in Part II, we describe the moduli spaces of the stable maps;
in Part III we explain the full package of the gluing theory that
provides {\em a} smooth structure on the moduli spaces;
in the part IV, we develop the virtual theory on the moduli spaces and
localization formula, at the end, as an application, we compute an
example.

Acknowledge. The idea of  this paper and  that of \cite{CT} was emerged
4 years ago. The drafts of papers have been written for quite a while, by some
reason, they have not been completed until recently.
First of all, special
thanks to Y. Ruan and G. Tian for their long time support on this project.
During this
long term preparation of papers, we would like to thank many people's
encouragement and discussion. The list includes G. Liu, K. Liu,  M. Liu,
 W. Zhang, G. Zhao, Q. Zheng and etc. The material of the paper was
explained as lectures in University of Wisconsin-Madison, Peking University.
We would like to thank their hospitality. We would also like  to W. Li,
Y. Long and J. Robbin for their interests in the lectures.
\vskip 0.2in
\begin{center}{\bf Part I. Preliminary}
\end{center}

\section{Complex structures on $\rone^2$}\label{section_1}

\def \rtwo{\mathbb{R}^2}

\subsection{Complex structures on $\rone^2$}\label{section_1.1}

A complex structure $j$ on $\rtwo$ is a linear automorphism of
$\rtwo$ with $j^2=-1$. It induces an orientation $o(j)$ on $\rtwo$
given by $v\wedge jv$ for any $0\not= v\in \rtwo$. Now fix an
orientation $or$ on $\rtwo$. Set
$$
J_{or}(\rtwo)=\{j| j^2=-1, o(j)=or \}.
$$

Fix a complex structure $j_o\in J_{or}(\rtwo)$. With a proper
chosen basis, we may write $j_o$ in terms of matrix as
$$
j_o=\left(
\begin{array}{cc}
0 & -1\\
1 & 0
\end{array}
\right).
$$
$(\rtwo, j_o)$ can be identified with a complex plane
$(\cplane,z)$ via
\begin{eqnarray*}
&&\phi_o: \rtwo\to \cplane;\\
&&z=\phi_o(x,y)=x+y\sqrt{-1}.
\end{eqnarray*}
Let $GL^+(2,\rone)<GL(2,\rone)$ be the subgroup that preseres
$or$. $GL^+(2,\rone)$ acts transitively on $J_{or}(\rtwo)$ via
$$
g\cdot j= gjg\inv.
$$
Via identification $\phi_o$, $GL(1,\cplane)$ is embedded in
$GL^+(2,\rone)$ as a subgroup. Then
\begin{lemma}\label{lemma_1.1.1}
$ J_{or}(\rtwo)\cong GL^+(2,\rone)/GL(1,\cplane)$.
\end{lemma}
{\bf Proof. }Since the isotropic group of the $GL^+(2,\rone)$ at
$j_o$ is $GL(1,\cplane)$, the lemma follows. q.e.d.

\subsection{Beltrami coefficients}\label{section_1.2}

Let $f: (\cplane,z)\to (\cplane, w)$ be a linear isomorphism
between two complex planes. Suppose
$$
w=f(z)=\alpha z+\beta \bar z.
$$
Then
\begin{equation}\label{eqn_1.2.1}
\mu_f= \alpha\inv\beta.
\end{equation}
is called the {\em Beltrami coefficient} of $f$ with respect to
coordinates $z$ and $w$.

Suppose that we change the coordinate of $w$-plane to $\hat
w$-plane by $\hat w=\gamma w,\gamma\in \cplane$. $f$ is
transformed to
$$
\hat f: (\cplane, z)\xrightarrow{f}(\cplane,w)\xrightarrow{\gamma}
(\cplane,\hat w).
$$
We find
$$
\mu_{\hat f}=\mu_f.
$$
This says that $\mu_f$ is independent of the coordinate choice of
$w$-plane.

Now Suppose that we change the coordinate of $z$-plane
 to $\hat z$-plane
by $ \hat z=\gamma\inv  z,\gamma\in \cplane$. Then $f$ is
transformed to
$$
\hat f: (\cplane, \hat
z)\xrightarrow{\gamma}(\cplane,z)\xrightarrow{f} (\cplane,w).
$$
We have
$$
\mu_{\hat f}= \mu_f\frac{\bar\gamma}{\gamma}.
$$
This implies that
\begin{equation}\label{eqn_1.2.2}
\omega_f= \mu_f \frac{d\bar z}{dz}
\end{equation}
is invariant on the first plane.  $\omega_f$ is an $(-1,1)$-form
on $z$-plane. We call it is the {\em Beltrami form} of $f$.

Given  a complex structure $j\in J_{or}(\rtwo)$ and an
identification
$$
\phi: (\rtwo,j)\to (\cplane,w),
$$
The identity map on $\rtwo$ induces a map $A_j$ via the diagram
$$
\begin{array}{ccc}
(\rtwo, j_o) &\xrightarrow{id} & (\rtwo,j)\\
\Big\downarrow\vcenter{\rlap{$\phi_o$}} &&
\Big\downarrow\vcenter{\rlap{$\phi$}} \\
(\cplane,z)&\xrightarrow{A_j} &(\cplane, w).
\end{array}
$$
We define
\begin{eqnarray*}
&&\mu:  J_{or}(\rtwo)\to \cplane;\\
&&\mu(j)= \mu_{A_j}.
\end{eqnarray*}
This map is well defined since $\mu(j)$ depends  on $\phi_o$, but
not on $\phi$.
\begin{prop}\label{prop_1.2.1}
$\mu$ is injective and $Image(\mu)=D$, the unit disk in $\cplane$.
\end{prop}
{\bf Proof. }Since $id$ (or $A_j$) is orientation preserving map,
one can check that $|\mu(j)|<1$. Hence $Image(\mu)\subset D$.

$\mu$ is injective: Suppose $\mu(j_1)=\mu(j_2)$. We have diagram
$$
A_{j_2}: (\cplane,z)\xrightarrow{A_{j_1}}(\cplane,w_1)
\xrightarrow{f_{12}}(\cplane,w_2),
$$
where $f_{12}$ is defined by the equation. Suppose
$$
w_2=f_{12}(w_1)= \alpha w_1+\beta \bar w_1.
$$
Then one can  check directly that
$$
\mu(j_1)=\mu(j_2) \iff \beta=0.
$$
This says that $f_{12}$ is holomorphic, and  so $j_1=j_2$.

$Image(\mu)=D$: let $\gamma$ be any complex number in $D$, we
solve $j$ such that $\mu(j)=\gamma$.
 Suppose $j=g\inv j_og$. Then we have
$A_g\inv$ defined by
$$
\begin{array}{ccc}
(\rtwo, j) &\xrightarrow{g} & (\rtwo,j_o)\\
\Big\downarrow\vcenter{\rlap{$\phi$}} &&
\Big\downarrow\vcenter{\rlap{$\phi_o$}} \\
(\cplane,w)&\xrightarrow{A_g} &(\cplane, z).
\end{array}
$$
By the definition of $g$, $A_g$ is holomorphic. Furthermore
$$
\phi_o g\phi_o\inv= A_g\circ A_j\inv: (\cplane,z)\to (\cplane,z).
$$
Hence,
$$
\mu_{\phi_o g\phi_o\inv}=\mu_{A_j}.
$$
Suppose
$$
g=\left(
\begin{array}{cc}
a & b\\
c & d
\end{array}
\right)
$$
and $\gamma=\alpha+\beta\sqrt{-1}$. Then
$$
\phi_o g\phi_o\inv(z)= (\frac{a+d}{2}+\frac{c-b}{2}\sqrt{-1})z
+(\frac{a-d}{2}+\frac{b+c}{2}\sqrt{-1})\bar z.
$$
Now set
$$
a=1+\alpha,d=1-\alpha,  b=c=\beta.
$$
We see that $\det(g)=1-\alpha^2-\beta^2>0$, which says that $g\in
GL^+(2,\rone)$, and
$$
\mu_{\phi_o g\phi_o\inv}=\gamma.
$$
This solves $\mu(j)=\gamma$. q.e.d.

\subsection{A Universal family of $J_{or}(\rtwo)$}\label{section_1.3}
We combine the result of previous two subsections:
$$
GL^+(2,\rone)/GL(1,\cplane)\cong J_{or}(\rtwo)\cong D.
$$
The second isomorphism is given by $\mu$ and $\mu(j_o)=0$. The
next proposition says that there exists  a canonical section
$\sigma$ (with respect to $j_o$) for the principle bundle
$$
\begin{array}{ccc}
GL(1,\cplane) &\rightarrow & GL^+(2,\rtwo)\\
&& \Big\downarrow \\
&& D
\end{array}.
$$
\begin{prop}\label{prop_1.3.1}
$\sigma(\gamma)= 1-jj_o$ for $\mu(j)=\gamma$.
\end{prop}
{\bf Proof. }Clearly $j\sigma(\gamma)=\sigma(\gamma) j_o$. It
remains to show that $\sigma(\gamma)\in GL^+(2,\rone)$.

Suppose $j=g j_og\inv$, where
$$
g=\left(
\begin{array}{cc}
a & b\\
c & d
\end{array}
\right)
$$
Without loss of generality, we assume that $\det(g)=1$. Then
$$
j=\left(
\begin{array}{cc}
ac+bd & -a^2-b^2\\
c^2+d^2 & -ac-bd
\end{array}
\right),
$$
and
$$
1-jj_o= \left(
\begin{array}{cc}
1 +a^2+b^2 & ac+bd\\
ac+bd & 1+c^2+d^2
\end{array}
\right).
$$
Using the fact $ad-bc=1$, we have $\det(1-jj_o)>0$. q.e.d.

\v Let $\mc R$ be a tautological family of $\rtwo$ with complex
structure parameterized by $D(\cong J_{or}(\rtwo))$: namely, we
have
$$
\mathfrak r:\mc R=(D\times \rtwo,\mc J)\to D,
$$
where $\mc J$ is the fiber-wise complex structure such that
$$
\mc J|_{\mathfrak r\inv \gamma}= \mu\inv(\gamma).
$$
\begin{prop}\label{prop_1.3.2}
There is a canonical trivialization with respect to $j_o$
$$
\Phi_o: \mc R_o:=D\times (\rtwo,j_o)\to \mc R.
$$
\end{prop}
{\bf Proof. } We set $\Phi_o$ fiber-wisely as
$$
\sigma(\gamma): \gamma\times (\rtwo,j_o)\to \mathfrak
r\inv(\gamma).
$$
q.e.d.


\section{Teichmuller spaces}\label{section_2}

\def \rtwo{\mathbb{R}^2}

\subsection{Complex structures on $\Sigma_g$}\label{section_2.1}

Let $\Sigma$ be an oriented  genus-$g$ surface  with orientation
$or(g)$. A complex structure $j$ on $\Sigma$ is a family of
complex structures on $\rtwo\cong T_x\Sigma$ parameterized by
$x\in \Sigma$. A complex structure $j$ induces an orientation
$o(j)$ on $\Sigma$. Set
$$
J(\Sigma)=\{j|o(j)=or(g)\}
$$
to be the set of complex structures on $\Sigma$ that is compatible
with the given orientation $or(g)$.

Now fix a point $j_o\in  J(\Sigma)$. Let
$\Omega^{-1,1}(\Sigma,j_o)$ denote the space of $(-1,1)$-forms on
Riemann surface $(\Sigma,j_o)$. A $(-1,1)$-form is locally
expressed in the form $f(z)d\bar z/dz$ in terms of local complex
coordinate $(\cplane, z)$. For any $j\in J(\Sigma)$, it yields
Beltrami coefficients $\mu(j(x))$ point-wisely, which
 gives a
$(-1,1)$-form, denoted by $\omega_j$, on $\Sigma$. Set
$$
B\Omega^{-1,1}(\Sigma,j_o)=\{\omega\in \Omega^{-1,1}(\Sigma,j_o)|
\|\omega\|<1 \}.
$$
Here $\|\omega\|=\max_{x} |f(x)|$ for $\omega=f(x)d\bar z/dz$.
Then by proposition \ref{prop_1.2.1}, we have
\begin{prop}\label{prop_2.1.1}
$j\to \omega_j$ gives an isomorphism $J(\Sigma)\cong
B\Omega^{-1,1}(\Sigma,j_o)$. Moreover $\omega_{j_o}=0$.
\end{prop}
Conversely, given a form $\omega\in B\Omega^{-1,1}(\Sigma,j_o)$,
we denote the corresponding complex structure by $j_\omega$.

\subsection{Teichmuller spaces $\mc T_g$}\label{section_2.2}

Here we give an informal review of Techmuller spaces $\mc
T_g$.

$\mc T_0$ consists of only one element, i.e, the standard sphere $S^2
=\cplane\cup \{\infty\}$.

$\mc T_1\cong \halfplane$, the upper half plane of $\cplane$.
Given $\lambda\in \halfplane $, we define a lattice
$$
L_\lambda= \{m+n\lambda|m,n\in \integer\};
$$
then the corresponding torus is
$$
T_\lambda= \frac{\cplane}{L_\lambda}.
$$

For $g\geq 2$, define
$$
\mc T_g = \frac{ J(\Sigma)}{Diff^+_0(\Sigma)}
=\frac{B\Omega^{-1,1}(\Sigma,j_o)}{Diff^+_0(\Sigma)}.
$$
Here $Diff^+_0(\Sigma)$ is the component of
$or(g)$-preserving-diffeomorphism group $Diff^+(\Sigma)$ that
contains 1. A classical theory on Teichmuller spaces says that the
quotient
$$
\frac{B\Omega^{-1,1}(\Sigma,j_o)}{Diff^+_0(\Sigma)}
$$
has a global slice. Let $H^{-1,1}(\Sigma,j_o)\subset
\Omega^{-1,1}(\Sigma,j_o)$ denote the space of holomorphic forms
and set
$$
BH(j_o)= B\Omega^{-1,1}(\Sigma,j_o)\cap H^{-1,1}(\Sigma,j_o).
$$
Then
\begin{theorem}\label{theorem_2.2.1}
$$
\frac{B\Omega^{-1,1}(\Sigma,j_o)}{Diff^+_0(\Sigma)}= BH(j_o).
$$
\end{theorem}
This theorem says that $BH(j_o)$ is a global slice of the
quotient. By the Riemann-Roch theorem, we know that $BH(j_o)$ is a
$6g-6$ dimensional ball.
\begin{corollary}\label{corollary_2.2.1}
$\dim \mc T_g= 6g-6, g\geq 2$.
\end{corollary}

We can also consider $\mc T_{g,m}$, the Teichmuller space of genus
$g$ Riemann surfaces with $m$ marked points. We give a complete
list.
\begin{enumerate}
\item[] $\mc T_{0,1}=\{(S^2,\infty)\}$; \item[] $\mc
T_{0,2}=\{(S^2,0,\infty)\}$; \item[] $\mc
T_{0,3}=\{(S^2,0,1,\infty)\}$; \item[] $\mc T_{0,m}= \mc
T_{0,3}\times ((S^2-\{0,1,\infty\})^{m-3}-\Delta), m >3$; \item[]
$\mc T_{1,1}=\mc T_1\times \{[0] \}, $ here $0\in \cplane$ and
$[0]$ denote the point in tori; \item[] $\mc T_{1,m}= \mc
T_{1,1}\times ((T_i-[0])^{m-1}-\Delta), m>1$ (refer this to
$\uni_1$ in \S\ref{section_2.3}); \item[] $\mc T_{g,m}= \mc
T_g\times (\Sigma^m-\Delta)$.
\end{enumerate}
In the first three terms, $S^2$ is identified with $\cplane\cup
\{\infty\}$; $\Delta$ is the big diagonal of the product $X^m$.

\subsection{Universal curves}\label{section_2.3}
By a universal curve, we mean a fibration
$$
\pi_{g,m}:\uni_{g,m}\to \mc T_{g,m}
$$
such that for $\mkj\in \mc T_{g,m}$ the fiber
$\pi_{g,m}\inv(\mkj)$ is the marked curve $\mkj$. We show the
existence of $\uni_{g,m}$.

\v\n{ \it Case 1, $g=0$}. Then
\begin{eqnarray*}
&&\uni_0=\uni_{0,0}=S^2; \uni_{0,1}=(S^2,\infty);\\
&&\uni_{0,2}=(S^2,0,\infty); \uni_{0,3}=(S^2,0,1,\infty)
\end{eqnarray*}
and
$$
\uni_{0,m}=\uni_{0,3}\times ((S^2-\{0,1,\infty\})^{m-3}-\Delta),
m>3.
$$
{\it Case 2, $g=1$.} We first construct $\uni_1$. Define an
$\integer\times \integer$ action on $\halfplane\times \cplane$ by
$$
(m,n)\cdot(\lambda, z)=(\lambda, z+m+n\lambda).
$$
Then
$$
\uni_1=\frac{\halfplane\times \cplane}{\integer\times \integer}
\to \halfplane
$$
is the universal curve.

$\uni_1$ can be topological trivialized to be $\halfplane\times
T_i, i=\sqrt{-1}$: define an $\integer\times \integer$ action on
$\halfplane\times \cplane$ by
$$
(m,n)\cdot(\lambda, z)=(\lambda, z+m+n\sqrt{-1}).
$$
Define the map
\begin{eqnarray*}
&&\tilde{\Phi}:\halfplane\times \cplane\to\halfplane\times \cplane\\
&&\tilde{\Phi}(\lambda,z)=(\lambda,\phi_\lambda(z)),
\end{eqnarray*}
where $\phi_\lambda:\cplane\to \cplane $ is the linear map
determined
 by
$\phi_\lambda(1)=1,\phi_\lambda(\sqrt{-1})=\lambda$. Then the map
is $\integer\times \integer$-equivariant. So $\tilde\Phi$ induces
an isomorphism
$$
\Phi: \halfplane\times T_{i}\to \uni_1.
$$
In particular $\Phi(\lambda,[0])=(\lambda,[0])$. We always assume
$\uni_1=\halfplane\times T_i$ from now on. Then
\begin{eqnarray*}
&& \uni_{1,1}=\uni_1\times \{[0]\}; \\
&& \uni_{1,m}=\uni_{1,1}\times ((T_i-\{[0]\})^{m-1}-\Delta), m>1.
\end{eqnarray*}
\v\n {\it Case 3, $g\geq 2$. } We set $\uni_g$ topologically to be
$$
\uni_g= \mc T_g\times \Sigma_g=BH\times \Sigma_g
$$
and the complex structure on the fiber over $\omega\in BH$ to be
$j_\omega$. This gives $\uni_g,g\geq 2$. Then set
$$
\uni_{g,m}=\uni_g\times (\Sigma^m-\Delta).
$$
We conclude that
\begin{theorem}\label{theorem_2.3.1}
The universal curve $\pi: \uni_{g,m}\to \mc T_{g,m}$ exists.
\end{theorem}


\section{Moduli space $M_{g,m}$}\label{section_3}

\subsection{Definitions}\label{section_3.1}
By definition,
\begin{eqnarray*}
&& \mc T_g=\frac{
J(\Sigma_g)}{Diff^+_0(\Sigma_g)};\\
&& M_g=\frac{ J(\Sigma_g)}{Diff^+(\Sigma_g)}.
\end{eqnarray*}
The group $ \Gamma_g= Diff^+(\Sigma_g)/Diff^+_0(\Sigma_g) $ is
called the {\em mapping class group}. Then
$$
M_g=\mc T_g/\Gamma_g.
$$
Similarly, one can define $M_{g,m}$, the moduli space of genus-$g$
curve with $m$-marked points.

When $g\leq 1$, $M_{g,m}=\mc T_{g,m}$. When $g\geq 2$, $M_{g,m}$
are orbifolds. We give  local descriptions for these orbifolds.

We recall some notions of orbifold. Let $X$ be an orbifold. For
any $x\in X$, there exists a neighborhood $U_x$ of $x$ such that
it is homeomorphic to $\rone^n/G_x$ for some finite group $G_x$.
The formal notion for these data is the so-called {\em
uniformization system} $(V, G_x,\phi)$: here $V$ is a smooth
manifold (usually is diffeomorphic to $\rone^n$), $G_x$ acts
smoothly on $V$ and
$$
\phi: V\xrightarrow{\pi} V/G_x \xrightarrow{\cong} U_x.
$$
$G_x$ is called the {\em isotropic group} of $x$. Clearly, such a
uniformization system describes the local of $x$.

Let
\begin{equation}\label{eqn_3.1.1}
\mkj_o=(\Sigma,j_o,x_{o1},\ldots, x_{om})\in M_{g,m}.
\end{equation}
We give a uniformization system for $\mkj_o$. The isotropic group
at $\mkj_o$ is
 $\aut(\mkj_o)$, the automorphism group of $\mkj_o$.
 By an automorphism of $\mkj_o$, we mean a bi-holomorphic
 maps of $(\Sigma,j_o)$ preserving
$x_{ok}$'s. Recall that
$$
\mc T_{g,m}= BH(j_o)\times (\Sigma^m-\Delta).
$$
$\aut(\mkj_o)$ acts on this space naturally as
$$
\sigma\cdot (\omega, y_1,\ldots, y_m) =(\sigma_\ast(\omega),
\sigma(y_1),\ldots, \sigma(y_m)).
$$
Then there exists a small $\aut(\mkj_o)$-invariant
 neighborhood $\tilde O$ of $\mkj_o$ such that $(\tilde O, \aut(\mkj_o),\pi)$ is
a uniformization system for $\mkj_o$.

Moreover, $\aut(\mkj_o)$ acts on $\uni_{g,m}$ similarly. Hence
\begin{equation}\label{eqn_3.1.2}
\frac{\uni_{g,m}|_{\tilde O}}{\aut(\mkj_o)}
\end{equation}
gives a universal curve for the neighborhood of $\mkj_o$ in
$M_{g,m}$.

\subsection{Hyperbolic metrics}\label{section_3.2}
A curve in $M_{g,m}$ is called {\em stable} if and only if
$2g+m\geq 3$. Let $\mkj_o\in \mc T_{g,m}$ be a stable curve given
as \eqref{eqn_3.1.1}. We may assign it a hyperbolic metric. This
is done as following: since $2g+m\geq 3$, the punctured surface
$(\Sigma-\{x_{o1},\ldots,x_{om}\},j_o)$ admits a universal
covering $\halfplane$ such that the punctured
 surface is  $\halfplane/\Gamma$ for some
Fuchsian group $\Gamma$. The hyperbolic metric on $\halfplane$
 induces a hyperbolic metric
on this punctured \rs.

Furthermore the neighborhood of $x_{oi}$ in $\Sigma$ has a nice
description.  Let
$$
\halfplane_{\geq 1}=\{x+y\sqrt{-1}\in \halfplane| y\geq 1\}.
$$
Then there exists a neighborhood of $x_{oi}$ (punctured at
$x_{oi}$) that is identified with
$$
\frac{\halfplane_{\geq 1}}{\langle z\to z+1\rangle }
$$
We call the area {\em the horocycle at $x_{oi}$}. This area can be
identified with the punctured disk
 $B^\ast (1)$ up to a rotation induced by:
$$
f(z)= c\exp(2\pi iz), c=\exp(2\pi).
$$
From now on, by horocycle, we always refer it as $B^\ast(1)$. We
note that $\aut(\mkj_o)$ acts on $B^\ast (1)$ as rotations.

We  deform the hyperbolic metric to be a {local-flat} metric such
that it is flat in $B^\ast (.5)$ and  is rotational invariant in
$B^\ast (1)$. So this new metric is still
$\aut(\mkj_o)$-invariant. We now assign the local-flat metrics
fiber-wisely to $\uni_{g,m}$. We denote the family metric to be
$h$ and the metric on the fiber over $\mkj$ as $h_\mkj$. Then
\begin{prop}\label{prop_3.2.1}
Let $\mkj\in \mc T_{g,m}, \sigma\in \aut(\mkj_o)$. Then
$$
h_{\sigma_\ast \mkj}= \sigma_\ast h_\mkj.
$$
\end{prop}
{\bf Proof. }We may assign hyperbolic metrics $\tilde h_{\mkj}$ on
each curve $\pi_{g,m}\inv(\mkj)$ and get a family metric $\tilde
h$. By the property of hyperbolic metrics,
\begin{equation}\label{eqn_3.2.1}
\tilde h_{\sigma_\ast \mkj}= \sigma_\ast\tilde h_\mkj.
\end{equation}
Moreover, $\sigma$ preserves horocycles: it maps horocycle of
$\mkj$ to $\sigma_\ast \mkj$ and behaves as rotations on $B^\ast
(1)$. Hence, \eqref{eqn_3.2.1} is still true for $h$. q.e.d.

\v\n We always assume that $\uni_{g,m}$ carries such a family of
metric $h$ if $2g+m\geq 3$.

\subsection{Local trivialization of universal curves}\label{section_3.3}
Recall that the local universal curve of $\mkj_o\in M_{g,m}$ is
given in the form \eqref{eqn_3.1.2}. Topologically, the fiber is
diffeomorphic to $(\Sigma, x_{o1},\ldots, x_{om})$, or $\Sigma
-\{x_{o1},\ldots, x_{om}\}$ if we use the language of puncture
curves. We now give a trivialization of $\uni_{g,m}|_{\tilde O}$
when $m\geq 1$.

\begin{prop}\label{prop_3.3.1}
There exists a smooth map
$$
\phi: \tilde O\times \Sigma\to\Sigma
$$
such that
\begin{enumerate}
\item for any $\mkj\in \tilde O$  with
$$
\pi_{g,m}\inv(\mkj)=(\Sigma, j, x_1,\ldots,x_m),
$$
 $\phi(\mkj,\ast): \Sigma\to\Sigma$
is diffeomorphic and $\phi(x_{ok})=x_{k}, 1\leq k\leq m$; \item
$\phi(\mkj,\ast)$ is holomorphic in  small neighborhoods of
$x_{ok}$ as a map
$$
\phi(\mkj,\ast): (\Sigma, j_o)\to (\Sigma,j);
$$
\item $\phi$ is $\aut(\mkj_0)$-equivariant in the sense
$$
\phi(\sigma\cdot \mkj, \sigma(z))= \sigma(\phi(\mkj,z)).
$$
\end{enumerate}
\end{prop}
{\bf Proof. }We explain for the case that $g\geq 2$. For $g=0,1$,
the proof is simpler, we leave it to readers. For simplicity, we
assume that $m=1$. We have horocycles on $\mkj_o$ and $\mkj$ given
by
$$
\zeta_{\mkj_o}: B^*(1)\to \Sigma-\{x_{o1}\}.
$$
and
$$
\zeta_{\mkj}:B^*(1)\to \Sigma-\{x_{1}\}.
$$
As $x_1$ is close to $x_{o1}$, we may assume that
$x_1\in\zeta_{\mkj_o}(B(.75))$ and
$$
\zeta_{\mkj_o}(B(.75))\subset \zeta_\mkj(B(1)).
$$
Set $x'_1= \zeta_{\mkj_o}\inv(x_1)$. We define a map $\phi_{\mkj}:
B(.75)\to B(.75) $ by
$$
\phi_{\mkj}(z)=\left\{
\begin{array}{lll}
z-x'_1, & \mbox{if} & |z|\leq .25;\\
z, &\mbox{if}& |z|\geq .5;\\
\eta(|\zeta_j(z)|)(z-x')+(1-\eta(|\zeta_j(z)|))z, & & \mbox{else}
\end{array}
\right.
$$
Here $\eta(t)$ is a cut-off function that is 1 when $t\leq .25$
and is 0 when $t\geq .5$. We now  set
 $$
\phi(\mkj,\ast)=\zeta_{\mkj}\phi_{\mkj} \zeta_{\mkj_o}\inv
$$
 on $\zeta_j(B(1))$ and  extended it over
$\Sigma$ by identity.

Conclusion (1) and (2) is obvious from the construction. It is
well known that the horocycle $\zeta_j$ depends smoothly with
respect to $\mkj$, hence $\phi$ is smooth.

We now explain (3). It is clear that
$$
\zeta_{\sigma\cdot j}= \sigma\circ\zeta_{j}.
$$
Also note that $\sigma$ is a rotation in the unit disk. Then (3)
can be verified directly. q.e.d.

\v\n We can now use $\phi$ to trivialize the universal curve over
the neighborhood of $\mkj_o$. Define
\begin{equation}\label{eqn_3.3.1}
\Phi_{\mkj_o}: \tilde O\times \mkj_o\to\uni_{g,m}|_{\tilde O}
\end{equation}
by
$$
\Phi_{\mkj_0}(\mkj, \Sigma)= \phi(\mkj, \Sigma).
$$
In fact, claim (1) in the proposition already serves the purpose.
Claim (3) implies that the trivialization is $\aut(\mkj_0)$
equivariant. Claim (2) is an additional property that is needed
later.


\section{Deligne Mumford moduli space $\bar M_{g,m}$}\label{section_4}

\subsection{Stable nodal curves}\label{sec_4.1}
Let $\bar M_{g,m}$ be the {\em Deligne-Mumford compactification}
of $M_{g,m}.$ We call it {\em Deligne-Mumford moduli space}. This
is a stratified space. The lower strata of $\bar M_{g,m}$ consists
of equivalence classes of stable nodal curves. A nodal curve is a
connected curve
$$
(\Sigma,j) = \bigcup_{k=1}^c (\Sigma_k,j_k), j=(j_1,\ldots,j_c),
$$
with normal crossing singularities
$$
\sing(\Sigma)=\{y_1,\ldots, y_s\}.
$$ 
We call these $y_i, 1\leq i\leq s $ the nodal points of $\Sigma$.  A marked point on
$\Sigma$ is a point $x\in \Sigma-\sing(\Sigma)$ Suppose we have
$m$-marked points $\{x_1,\ldots,x_m\}$. Let
\begin{equation}\label{eqn_4.1.1}
 \mkj=(\Sigma,j,
x_1,\ldots,x_m)
\end{equation}
be a nodal curve with $m$-marked points.

\def \comp{\mathrm{comp}}

Set
\begin{equation}\label{eqn_4.1.2}
D: \{1,\ldots, m\}\to \{1,\ldots,c\}
\end{equation}
to be the map that assigns marked point $x_i$ to component
$\Sigma_{D(i)}$. For each $y\in \sing(\Sigma)$, it is contained in
two components $\Sigma_{c_1}$ and $\Sigma_{c_2}$. Here $c_1$ may
equal to $c_2$. We define the set $\comp(y)=\{c_1,c_2\}$ (or,
$\comp(y)=\{c_1\}$ if $c_1=c_2$).

Each component of this curve is
$$
\mkj_k=(\Sigma_k, j_k, \{x_i\}_{i\in D\inv (k)}, \sing(\Sigma)\cap
\Sigma_k).
$$
This is an one-component-curve maybe with nodal points. The nodal
points are come from those singular $y$'s with $\comp(y)=\{k\}$.
Such a component admits a normalization  $\mk N(\mkj_k)$. Then the
normalization of $\mkj$ is defined to be the disjoint union
$$
\mk N(\mkj):= \coprod_{k=1}^c \mk N(\mkj_k).
$$
Recovering $\mkj$ from $\mk N(\mkj)$ is standard. It is given by a
proper quotient map
$$
\pi:\mk N(\mkj)\to \mk N(\mkj)/\sim \cong \mkj.
$$
\begin{defn}\label{def_1_2}
$\mkj$ is stable if $\mk N(\mkj_k)$ is stable for each $k$.
\end{defn}
Two  curves
$$
\mkj=(\Sigma, \mkj,x_1,\ldots,x_m) \mbox{ and } \mkj'=(\Sigma',
\mkj', x_1',\ldots,x_m')
$$
are equivalent if there exists a homeomorphism $\sigma: \Sigma\to
\Sigma'$ such that $\sigma(x_i)=x_i'$ and the natural induced map
$\mk N(\sigma):\mk N(\mkj)\to\mk N(\mkj')$ is bi-holomorphic.



\subsection{Data of stratum}
Let $\mkj$ be an $m$-marked stable nodal curve given by
\eqref{eqn_4.1.1}.
 We  assign  the following combinatoric data to this curve:
\begin{enumerate}
\item  a (weighted) connected graph (with tails) $T$ (refer to
item-2 for "weighted" and item-3 for "tail"):
 Let $V$ and $E$
be the set of vertices and edges of $T$ respectively, then  each
$k\in V$ stands for a component $\Sigma_k$ and $n\in E$ stands for
a nodal point $y_n\in \sing(\Sigma)$; \item the genus $g_k$ of
$\Sigma_k$ for each $k\in V$:$g_k$ is the weight of $k$ that is
mentioned in item-1; the data of genus is denoted by
$$
\mk g= (g_1,\ldots,g_c).
$$
Set
\begin{equation}\label{eqn_4.2.1}
g=g(T):=\sum_{v=1}^k g_k+ \rank H_1(T);
\end{equation}
\item a map
$$
D: \{1,\ldots,m\}\to \{1,\ldots,c\}
$$
mentioned in \eqref{eqn_4.1.2}: for each $1\leq j\leq m$ we assign
it a tail, that is,
 for
$D(j)=k$ we add $j$-th tail to vertex $k$.
\end{enumerate}
We denote the data by $S=(T, D,\mk g)$  and call it a {\em stratum
data} in $\bar M_{g,m}$. Such  a data is called {\em stable} if
for each vertex $k$,
$$
2g_k+ \mr{val}(k) \geq 3.
$$
Here $\mr{val}(k)$ is the valency of vertex $k$ (Tails are counted
for valency). It is easy to check that
\begin{claim}\label{claim_4.2.1}
$\mkj$ is stable if the data $S$ given by $\mkj$ is stable.
\end{claim}
On the other hand,
\begin{defn}\label{defn_4.2.1}
The genus of $\mkj$ is defined to be $g=g(T)$. $\mkj$ is called a
stable $(g,m)$-curve.
\end{defn}
We define $\bar M_{g,m}$ to be the set of equivalence classes of
stable $(g,m)$-curves. This space admits a natural stratification
given by data $S$'s: let $S=(T,D,\mk g)$ be a stable data with
$g=g(T)$, we define the stratum $ M_S\subset \bar M_{g,m} $ to be
the set of curves that give data $S$. The topology of $\bar
M_{g,m}$ is not clear at the moment. However this is studied
intensively (\cite{???}). It is well known that $\bar M_{g,m}$ is
a smooth orbifold of dimension $6g-6+2m$ if $2g+m\geq 3$.

In the rest of the section. We describe the strata and their
neighborhoods in $\bar M_{g,m}$ more carefully.

\subsection{Some facts of data $S$}\label{section_4.3}
 Let $S$ be a (stable) stratum data. There is an automorphism
group $\aut(S)$ of $S$ defined as following
\begin{defn}\label{defn_4.3.1}
We say $\gamma\in \aut(S)$ if $\gamma: T\to T$ is a graph
automorphism preserving weights and tails. Be precise, it induces
isomorphisms $\gamma: V\to V$ and $\gamma: E\to E$ such that
$$
\gamma(e(k_1,k_2))= e(\gamma(k_1),\gamma(k_2))
$$
and
$$
 g_{\gamma(k)}= g_k, \mbox{ \ } D(j) = \gamma(D(j)).
$$
\end{defn}

 Let $\mc D_{g,m}$ be the set of stable stratum data. It
can be  shown that
\begin{lemma}\label{lemma_1_2}
$|\mc D_{g,m}|<\infty$.
\end{lemma}
We skip the proof. The stability is crucial for the lemma.

For the set $\mc D_{g,m}$ we can assign a {\em partial order
$\prec$}. Let $S= (T,D,\mk g)$ be a data. Let $e=e(v,w)$ be an
edge of $T$. We can define a new data $S'$ by the following
modifications on $S$:
\begin{itemize}
\item  a new graph $T'$ is obtained by (i) erasing edge $e$, (ii)
identifying vertices $v$ and $w$ and  denote the new vertex by
$v'$; \item $g_{v'}$ is defined to be
$$
g_{v'}= g(T)- \sum_{k\not= v,w} g_k -\rank H^1(T');
$$
\item $D'(i)=v'$ if $D(i)=v$ or $w$. The attaching vertices of
tails are changed properly by new $D'$.
\end{itemize}
By this way, we say that $S'$ is a {\em contraction of $S$ at edge
$e$}. We write $S'=S(e)$. Similarly, we can define the contraction
$S'=S(e_1,\ldots, e_l)$ of $S$ at edge $e_1,\ldots, e_l$. Now, we
say that $S\prec S'$ if $S'$ is a contraction of $S$. This induces
a partial order on the strata of $\bar M_{g,m}$. In fact, this is
compatible with what we mean by "lower": $M_S $ is lower than
$M_{S'}$ if and only if $S\prec S'$.

Let $T$ be the simplest graph that consists of 1 vertice and no
edge. It defines an $S_0$ and the stratum is just
$M_{S_0}=M_{g,m}$.

\subsection{Strata $M_S$}\label{section_4.4}
Let $\mkj\in M_S$. The notions for $\mkj$ (cf. \eqref{eqn_4.1.1})
and $S$ are same as before. Recall that we have  normalizations
$\mk N(\mkj)$ and $\mk N(\mkj_k)$, with
$$
\pi: \mk N(\mkj)\to \mkj.
$$
We write
$$
\pi_S(\mk N(\mkj))=\mkj.
$$
These two maps are different! For $\pi$ the variable is a point on
\rs, while the variable for $\pi_S$ is a curve $\mk N(\mkj)$.
Suppose
\begin{equation}\label{eqn_4.4.1}
\mk N(\mkj_k) =(\tilde{\Sigma}_k, i_{k}, x_{k_1},\ldots,
x_{k_{m_k}}, \bar y_{k_1},\ldots,\bar y_{k_{s_k}}).
\end{equation}
Here $x$'s are marked points on $\Sigma_k$ and $\bar y$'s
correspond to nodal points. Be precisely, we may further assign an
edge $n=e(\bar y)\in E$ for $y$  that corresponds  to the nodal
point $\pi(\bar y)=y_n$.

Let
$$
\mc T_S: = \mc T_{g_1, m_1+s_1}\times \cdots\times \mc T_{g_c,
m_c+s_c}.
$$
This has a universal curve
$$
\pi_S:\uni_S=\uni_{g_1, m_s+s_s}\times \cdots\times \uni_{g_c,
m_c+s_c}\to \mc T_S.
$$
Since $\mk N(\mkj)\in \mc T_S$, we represent it by $\pi_S\inv(\mk
N(\mkj))$.

 We are now ready to describe the orbifold
structure of $M_S$ at $\mkj$. The isotropic group is $\aut(\mkj)$:
$\aut(\mkj)$  is a fibration
$$
\phi: \aut(\mkj)\to \aut(S)
$$
and acts on  $\mc T_S$. We explain this. Suppose that $\mk
N(\mkj)\in \mc T_S$ is
$$
\mk N(\mkj)=(\mk N(\mkj_1),\ldots, \mk N(\mkj_c)).
$$
Let $\gamma\in \aut(S)$. Then $\phi\inv (\gamma)$ is given by the
following elements: define $\lambda_k: N(\mkj_k)\to
N(\mkj_{\sigma(k)})$ such that
\begin{itemize}
\item the map
$$
\lambda_k: (\Sigma_k,i_k)\to (\Sigma_{\sigma(k)}, i_{\sigma(k)})
$$
is bi-holomorphic; \item $\lambda_k$ preserves marked points
$x_i$'s; \item $\lambda_k$ preserves $y$-set and
$$
\lambda_k(e(y))= \sigma(e(\lambda(y))).
$$
\end{itemize}
By this way, we define $\aut(\mkj)$. It acts naturally on the
neighborhood of $\mk N(\mkj)$ in $\mc T_S$. Let $\tilde O$ be  an
$\aut(\mkj)$-invariant neighborhood of $\mk N(\mkj)$. Then
$(\tilde O, \aut(\mkj), \phi)$ yields a uniformization system of
$\mkj$ in $M_S$ via $\pi_S$:
$$
\tilde O\xto{\phi} \tilde O/\aut(\mkj)\xto{\pi_S} M_S.
$$
All these charts form the orbifold $M_S$.

As before, $\aut(\mkj)$ also acts on  universal curve $\uni_S$.
Hence it induces an universal curve over $\tilde O/\aut(\mkj)$:
$$
\frac{\uni_S|_{\tilde O}}{\aut(\mkj)}.
$$
On the other hand, by using the trivialization constructed in
proposition \ref{prop_3.3.1}, there exists a trivialization of the
universal curve given by an $\aut(\mkj)$-equivariant map
\begin{equation}\label{eqn_4.4.2}
\Phi_\mkj: \tilde O\times \mk N(\mkj)\to \uni_S|_{\tilde O}.
\end{equation}

\subsection{Smoothing nodal curves at nodal points}\label{section_4.5}
Let $\mkj\in M_S$ be a nodal curve. For any nodal point $y$ and a
complex number $0\not= \rho$ with small radius, we can smoothen
$\mkj$ at $y$ and get a new curve $\mkj_{y,\rho}$. This is what we
mean by smoothing. We now explain this procedure.

Without loss of generality, we suppose that $\pi\inv(y)\subset \mk
N(\mkj)$ consists of  $v_1\in  \Sigma_1$ and $v_2\in \Sigma_2$.
 We treat $v_i$ as marked points on $\Sigma_i$.
 Since we are only concerned the local of $v_i$. We may assume
 that $\Sigma_i$ are smooth curves.
 The neighborhood
of $v_i$ can be canonically identified with  balls $B(1)$ up to
rotations in $z_i$-planes for $i=1,2$: if $\Sigma_i$ is stable, we
refer it as the horocycle at $v_i$; otherwise, $v_i$ is a special
point£¬ 0 or $\infty$, on $S^2$, we then refer the ball to be the
semi-sphere containing the point. We write the balls $B_{v_i}(1)$.
Furthermore we have
$$
\phi_i: B^\ast_{v_i}(1)\to (-\infty,0]_i\times S^1
$$
by
$$
\phi(r\ei)= (\log r, \theta).
$$
We write the punctured surfaces as
\begin{equation}\label{eqn_4.5.1}
\Sigma_i-\{v_i\}=\Sigma_{0i}\cup (-\infty,0]_i\times S^1,i=1,2.
\end{equation}

The neighborhood of $y$ can be put in $\cplane\times \cplane$ as
\begin{equation}\label{eqn_4.5.2}
\pi(B_{y_1}(1)\cup B_{y_2}(1)) =\{z_1z_2=0\}\cap B(1).
\end{equation}
For $\rho\in \cplane^\ast$ we  deform  \eqref{eqn_4.5.1} to
$$
\{z_1z_2=\rho\}\cap B(1).
$$
The new curve is denoted by
$$
\mkj_{y,\rho} =(\Sigma_\rho, i_\rho,\ldots).
$$
This smoothing procedure can be described explicitly. Set $\rho=
r_0e^{i\theta_0}$. Let us focus on $B_{v_i}(1)$. The remainder of
$\Sigma_{0i}$ remains unchanged in the whole process. We cut off
the cylinder ends of two cylinders at $\{1.25\log r_0\}\times
S^1$, namely we get
$$
[1.25\log r_0,0]_i\times S^1.
$$
Then we glue two tubes along a sub-tube of length $-0.5\log r_0$
with a twisted angle $\theta_0$. That is, we identify
$$
(\log r_0+t, \theta)_1\leftrightarrow (\log r_0-t, \theta
+\theta_0)_2, t\in [0.25\log r_0,- 0.25\log r_0].
$$
The resultant curve is then $\mkj_{y,\rho}$.

 We
note that the plane $\cplane$ of $\rho$ can be treated as
\begin{equation}\label{eqn_4.5.3}
T_{v_1}{\Sigma}_1\otimes T_{v_2}{\Sigma}_2.
\end{equation}
We denote the space as $\cplane_y$. Then we just construct a map
\begin{equation}\label{eqn_4.5.4}
\gs_\mkj: B_\epsilon \subset \cplane_y \to \bar M_{g,m}.
\end{equation}
Here $\epsilon$ is any small constant less than 1.

We remark that $\gs_\mkj$ is injective if and only if  $\mkj$ is
stable.

\subsection{Normal bundles of $M_S$ in $\bar M_{g,m}$}\label{section_4.6}
It is  natural to ask the neighborhood of  $\mkj\in M_S$ in $\bar
M_{g,m}$ when $S$ is stable. In particular, this is to ask what is
the normal direction of $M_S$. We assume that $S$ is stable.

Suppose $S=(\mk g,T,D)$. Given a point $\mkj\in M_S$, we define a
fiber
$$
\cplane_\mkj^{|E|}:= \oplus_{s=1}^{|E|}\cplane_{y_s}.
$$
Here $E$ is the set of edges of $T$. By this, we define an
orbifold bundle
$$
L_S\to M_S
$$
In fact, we have
$$
\tilde{L}_S\to \mc T_S.
$$
with fiber $\cplane_\mkj$. $\aut(\mkj)$ acts naturally on
$\tilde{L}_S$. Locally, the quotient gives the uniformization
system for $L_S$.

Let $O\subset M_S$ be a proper open subset. The gluing map
described earlier defines a neighborhood of $O$ in $\bar M_{g,m}$:
\begin{equation}\label{eqn_4.6.1}
\gs_S:  L_{S,\epsilon}|_O \to \bar M_{g,m},
\end{equation}
where $\gs_S(\mkj, \rho)= \gs_\mkj(\rho)$. Here by $
L_{S,\epsilon}$ we mean an $\epsilon$-neighborhood of 0 section.
It is known that $\gs_S$ is injective and locally diffeomorphic
when $S$ is stable.

Let $S'=S(e_1,\ldots, e_l)$. We can similarly define a sub-bundle
$L_{S,S'}$ of $ L_S$ by requiring the fiber to be
$$
{L}_{S,S'}|_{\mkj}=\bigoplus_{
\begin{array}{c}
y,
e(y)=e_j\\
1\leq j\leq l
\end{array}}
\cplane_{y}.
$$
We can then similarly define
\begin{equation}\label{eqn_4.6.2}
\gs_{S,S'}:  L_{S,S',\epsilon}^0|_O \to M_{S'}.
\end{equation}
Here, by $L^0_{S,S',\epsilon}$ we mean that the set of points
whose  all coordinates in fiber direction are not zero. For
example, by $V^0_\epsilon$ for a vector space $V=\cplane ^n$ we
have
$$
V^0_\epsilon= \{(z_1,\ldots,z_n)|z_i\not= 0, |z_i|\leq \epsilon
\}.
$$
Finally, we remark that the above discussion can be generalized to
unstable data $S$: $\tilde L_S$, $L_S$ and  $\gs$ are still
available. The only difference is that $\gs$ is neither injective
nor locally diffeomorphic. But this is  important when we consider
stable maps.

\vskip 0.2in

\begin{center}
{\bf Part II. Moduli spaces $\om_{g,m}(X,A)$ }
\end{center}

\section{Moduli spaces of stable maps }\label{section_5}

We review the fundamental facts of stable maps in this section, such as
notations of stable maps and related facts
(\cite{FO}).

\subsection{Stable maps without nodal points}\label{section_5.1}

 Let $(X,\omega)$ be a $C^{\infty}$ compact,  closed
symplectic manifold of dimension $2n$. We choose an $\omega$-tamed
almost complex structure $J$ on $X$, i.e, $\omega(\cdot, J\cdot)$
is positive. Define a $J$-compatible Riemannian metric
$\langle\cdot,\cdot\rangle$ by
$$ \langle V,W\rangle:= \frac{1}{2}\left(\omega(V,JW) + \omega(W,JV)\right).$$
Fix an element $A\in H_2(M,Z)$.

Let $\mkj=(\Sigma, j)$ be a \rs\ without nodal points. Consider
the space of smooth maps
$$
\mr{Map}_\mkj(X,A)= \{u: \Sigma \to X| [u(\Sigma)]=A \}.
$$
Here $[u(\Sigma)]$ represents the homology class of $u(\Sigma)$.

 Define
$$
\dbar(u)= \frac{1}{2}(du+ J\cdot du\cdot j).
$$
We say that the map $u$ is {\em $J$-holomorphic} if
\begin{equation}\label{eqn_5.1.1}
\dbar(u) =0.
\end{equation}
Let $\wtm_\mkj(X,A)\subset \mr{Map}_\mkj(X,A)$ be the space of all
$J$-holomorphic maps.

We explain $\dbar$. Given $u\in \mr{Map}_\mkj(X)$, we have a
bundle $T^\ast\Sigma\otimes u^\ast TX$ over $\Sigma$. The space of
sections of this bundle is denoted by
$$
\mr{End}(T\Sigma,u^\ast TX)= \Gamma(T^\ast\Sigma\otimes u^\ast
TX).
$$
According to $j\otimes u^\ast J$,  it has a decomposition
$$
\mr{End}(T\Sigma,u^\ast TX)= \mr{End}^{1,0}(T\Sigma,u^\ast TX)
\oplus\mr{End}^{0,1}(T\Sigma,u^\ast TX).
$$
To be consistent with conventions, we set
$$
\Omega^{0,1}_\mkj(u^\ast TX)= \mr{End}^{0,1}(T\Sigma,u^\ast TX).
$$
Now note that $du\in \mr{End}(T\Sigma,u^\ast TX)$.
 Then $\dbar u\in \Omega^{0,1}_\mkj(u^\ast TX)$
  is the $(0,1)$-component of $du$.

We summarize that we have a bundle
$$
\tilde{\mc E}_\mkj\to \mr{Map}_\mkj(X,A)
$$
with fiber $\Omega^{0,1}_\mkj(u^\ast TX)$ and $\dbar$ is a section
of $\tilde{\mc E}_\mkj$. Then
$$
\wtm_j(X,A) = \mr{zero \ section}\cap \dbar =\{u|\dbar u=0\}.
$$
We may replace $\mkj=(\Sigma, j)$ by a marked \rs\ $\mkj=(\Sigma,
j, x_1,\ldots, x_m)$.

Define
$$
\M_\mkj(X,A)= \frac{\wtm_\mkj(X,A)}{\aut(\mkj)}.
$$
Let $\aut(u,\mkj)$ be the stabilizer of the action for point $u\in
\wtm_\mkj(X,A)$. $u$ is called stable if $|\aut(u,\mkj)|\leq
\infty$.
\begin{remark}\label{rmk_5.1.1}
$(u,\mkj)$ is stable if $\mkj$ is stable. For this case, $2g+m\geq
3$ and we call the map is pre-stable. Otherwise, it is called
pre-unstable. We have four possibilities $(g,m)=(1,0), (0,0),
(0,1)$ and $(0,2)$ for pre-unstable maps.
\end{remark}
\begin{prop}\label{prop_5.1.1}
$u$ is stable if either $\mkj$ is stable or $u$ is not constant.
\end{prop}
\v We now allow $\mkj$ varies. Let $\mc T_{g,m}$ and $M_{g,m}$ be
spaces of curves described in \S \ref{section_3.1}. Define
$$
\wtm_{g,m}(X,A)=\coprod_{\mkj\in \mc T_{g,m}}\wtm_\mkj(X,A).
$$
Suppose that $(u,\mkj)$ and $(u',\mkj')$ are two maps. We say that
they are equivalent if there exists an isomorphism $\sigma:
\mkj\to \mkj'$ such that $u=u'\circ \sigma$. Define the moduli
space to be
$$
\M_{g,m}(X,A)= \wtm_{g,m}(X,A)/\sim.
$$
Roughly speaking, if the section $\dbar$ transverses to the 0-section,
$\M_{g,m}(X,A)$ is an orbifold. We postpone the complete
discussion to \S?. Here we describe its local structure. Let
$(u,\mkj)\in \M_{g,m}(X,A)$. Recall that a neighborhood $O_\mkj$
of $\mkj\in \M_{g,m}$ is given by
$$
O_\mkj= \frac{\tilde O_\mkj}{\aut(\mkj)}.
$$
Define
$$
\wtm_{\tilde O}(X,A) = \coprod_{\mkj'\in \tilde O_\mkj}
\wtm_{\mkj'}(X,A)
$$
and set
$$
\M_O(M,A)= \frac{\wtm_{\tilde O}(X,A)}{\aut(\mkj)}.
$$
This gives a neighborhood of $(u,\mkj)$. Furthermore, we can
choose a neighborhood $\tilde U$ of $(u,\mkj)$ in the numerator
that is invariant under the action of $\aut(u,\mkj)$. Then
$(\tilde U,\aut(u,\mkj),\pi)$ is a uniformization system of
$$
U:=\tilde U/\aut(u,\mkj).
$$

\subsection{Stable maps with nodal points}\label{section_5.2}
Let $\mkj=(\Sigma, i,x_1,\ldots, x_m)$ be a nodal curve (may  not
be stable). Most of the discussion in \S\ref{section_4.1.1} still
works for unstable curves. We still have notions $\mkj_k, \mk
N(\mkj)$ and etc. Similarly, we can define  data $S=(T,D,\mk g)$.
Suppose $g= g(T)$.

Let $u: \Sigma\to X$ be a continuous map such that
$[u(\Sigma)]=A$. We say that it is holomorphic if each restriction
$u_k:=u |_{\mkj_k}$ lifts to a $(J,\mk N(\mkj_k)$- holomorphic map
$$
\tilde u_k: \tilde{\Sigma}_k\to X.
$$
Here $\tilde\Sigma$ is the normalization of $\Sigma$.
\begin{defn}\label{defn_5.2.1}
We say that $(u, \mkj)$ is a $(g,m)$-stable map if $(\tilde{u}_k,
\mk N(\mkj_k))$ is stable for any $k$. The space of stable maps is
denoted by $\widetilde{\om}_{g,m}(X,A)$.
\end{defn}
Two stable maps  are equivalent, denoted by $(u,\mkj)\sim
(u',\mkj')$, if there exists an isomorphism $\sigma: \mkj\to
\mkj'$ such that $u=u'\circ \sigma$. Let
$$
\om_{g,m}(X,A)=\widetilde{\om}_{g,m}(X,A)/\sim.
$$

The moduli space $\om_{g,m}(X,A)$ has a similar stratification as
that of $\om_{g,m}$. Let $S=(T,D,\mk g)$ be a stratum data of
$\om_{g,m}$. We add an extra data
$$
\mk A= (A_1,\ldots, A_c)
$$
to it. Set
$$
A=A_1+\cdots+ A_c.
$$
We denote the new data to be $\strata=(T,D,\mk g,\mk A)$.
$\strata$ is a {\em stratum data of stable map}. Here $A_k$
represents the homology class of $u_k$. We then define the stratum
$\M_{\strata}(X,A)\subset \om_{g,m}(X,A)$ to be the set that
consists of equivalence classes of stable maps described above
with the property $[u_k(\Sigma_k)]=A_k$.

Let $\mc D_{g,m}^A$ be the set stratum data of stable map. As
before
\begin{lemma}\label{lemma_5.2.1}
$|\mc D_{g,m}^A|<\infty$.
\end{lemma}
We can also define a partial order $\prec$ on this set as in
\S\ref{section_4.3}. The only extra information we should add  is
that
$$
A_{v'}=A_v+A_w.
$$
By this, we can define $\prec$ in $\mc D_{g,m}^A$. So we can also
say that the stratum $\M_{\strata}(X,A)$ is lower than
$\M_{\strata'}(X,A)$ if $\strata\prec \strata'$.

With proper topology(\cite{PW}\cite{Ye}), one has
\begin{theorem}\label{theorem_5.2.1}
$\om_{g,m}(X,A)$ is compact. The closure of $\M_{\strata}(X,A)$ is
$$
\om_{\strata}(X,A)=\M_\strata(X,A)\cup \coprod_{\strata'\prec
\strata} \M_{\strata'}(X,A).
$$
\end{theorem}

Using the local description of $M_S$, we would like to give a
local description for $\M_{\strata}(X,A)$ as well. Let $(u,\mkj)$
be a $J$-holomorphic map in the stratum. Suppose $(\tilde
O,\aut(\mkj),\pi)$ is a uniformization system of a neighborhood
$O$ of $\mkj$ in $M_S$ (refer notations to \S\ref{section_4.4}.
 Let
$$
\wtm_{\tilde O}(X,A)=\coprod_{\mkj'\in \pi_S(\tilde O)}
\wtm_{\mkj'}(X,A).
$$
Then
$$
\M_O(X,A)= \frac{\wtm_{\tilde O}(X,A)}{\aut(\mkj)}
$$
is a neighborhood of $(u,\mkj)$ in $\M_\strata(X,A)$. As before,
we can choose a neighborhood $\tilde U$ of $(u,\mkj)$ in the
numerator that is invariant under the action of $\aut(u,\mkj)$.
Then $(\tilde U,\aut(u,\mkj),\pi)$ is a uniformization system of
$$
U:=\tilde U/\aut(u,\mkj).
$$


\section{Analytic set-up}\label{section_6}

\subsection{Analytic set-up for $\M_\mkj(X,A)$}\label{section_6.1}

Let $\mkj\in M_{g,m}$. We first assume that $\mkj$ is stable.
Recall that $\wtm_\mkj(X,A)$ is viewed as zeros of section $\dbar$
of bundle $\tilde{\mc E}_\mkj$. In order to show the smoothness of
$\wtm_{\mkj}(M,A)$, we need put Sobolev norms on  these spaces and
apply the transversality theorem for Banach manifolds.

Let $p>2$ be an even integer.  We denote by
 $\chi^{1,p}_\mkj(X,A)$ the space of
continuous map $u:\Sigma \rightarrow M$ of class $W^{1,p}$ such
that $[u(\Sigma)] = A$. We usually simplify the notation to be
$\chi_\mkj^{1,p}$.
 The space $\chi^{1,p}_\mkj$ is an infinite
dimensional Banach manifold. For any  map $u\in \chi^{1,p}_\mkj$,
its tangent space is the Banach space $W^{1,p}(u^*TM)$ of
$W^{1,p}$-vector fields $\zeta$ along u. The point-wise
exponential map
$$W^{1,p}(u^*TM)\rightarrow M: \zeta \rightarrow \exp_u\zeta$$
identifies a neighborhood $\mc U_u$ of 0 in $W^{1,p}(u^*TM)$
 with a neighborhood of u in $\chi^{1,p}_\mkj$. We have a coordinate chart $(\exp_u\mc U_u,\exp^{-1}_u)$.
Without loss of generality, we assume that $\mc U_u\subset
W^{1,p}(u^\ast TM)$ is a neighborhood of $u$.  The tangent bundle
is
$$
T\chi_\mkj^{1,p}\to \chi_\mkj^{1,p},
$$
a bundle with fiber $W^{1,p}(u^\ast TM)$. We denote the bundle by
$\tilde{\mc F}_{\mkj}$.

Similarly, we consider bundle $\tilde{\mc E}_\mkj$ over
$\chi_\mkj^{1,p}$. We put $L^p$ norm on the fiber. Hence the fiber
over $(u,\mkj)$ is
$$
L^p(\Lambda^{0,1}_\mkj(u^\ast TM)).
$$
 Recall that
$\dbar$ is a section of this bundle.

Now fix $u_o\in \chi^{1,p}_\mkj(X,A)$. We now trivialize
$\tilde{\mc E}_{\mkj}$ and $\tilde{\mc F}_{\mkj}$ over a small
neighborhood $\mc U_{u_o}$ of $u_o$. We trivialize $\tilde{\mc
F}_\mkj$ first. Let $\zeta\in W^{1,p}(u_o^\ast TM)$ and
$u=\exp_{u_o}\zeta$. Then the parallel transformation along path
$\exp_{u_o}s\zeta$
$$
P_\zeta: W^{1,p}(u^\ast TM)\to W^{1,p}(u_o^\ast TM)
$$
identifies two fibers. This defines a trivialization
$$
T\chi_\mkj^{1,p}|_{\mc U_{u_o}}\cong \mc U_{u_o} \times
W^{1,p}(u_o^\ast TM).
$$

Let $\Pi^{0,1}$ be the projection
$$
\Pi^{0,1}: L^p(T^\ast \Sigma\otimes u_o^\ast TM) \to
L^p(\Lambda^{0,1}_{\mkj}(u_o^\ast TM)).
$$
Then
$$
\Pi^{0,1}\circ P_\zeta: L^p(\Lambda^{0,1}_\mkj(u^\ast TM)) \to
L^p(\Lambda^{0,1}_{\mkj}(u_o^\ast TM))
$$
yields a trivialization of $\tilde{\mc E}_{\mkj}$.

We summarize the data we have:
\begin{itemize}
\item a base space $\mc U_{u_o}$; \item a bundle $ \tilde{\mc
E}_{\mkj}\cong  \mc U_u \times \Lambda^{0,1}_{\mkj}(u_o^\ast TM)$;
\item a  section $\bar\partial$ of the bundle $\tilde{\mc
E}_{\mkj}$; \item a tangent bundle $\tilde{\mc F}_{\mkj}$.
\end{itemize}
Let $(u,\mkj)\in \mc U_{u_o}$. The linearization of $\dbar$ at
$(u,\mkj)$ is
$$
D_{u,\mkj}: \tilde{\mc F}_{\mkj}|_{(u,\mkj)} \to \tilde{\mc
E}_{\mkj}|_{(u,\mkj)}.
$$
Be precise,  we have
$$
D_{u,\mkj}: W^{1,p}(u_o^\ast TM)\cong W^{1,p}(u^\ast TM)\to
L^p(\Lambda^{0,1}_\mkj(u^\ast TM)) \cong
L^p(\Lambda^{0,1}_{\mkj}(u_o^\ast TM)).
$$
Explicitly, by ignoring the identifications on the two ends given
by trivialization
$$
D_{u,\mkj}(\xi)= \frac{1}{2}(\nabla\xi+ J(u)\nabla \xi j +
\nabla_\xi J du j)
$$
\begin{prop}\label{prop_6.1.1}
The index of $D_{u,\mkj}$ is $n(2-2g)+ 2c_1(A)$.
\end{prop}
This follows from the Riemann-Roch theorem.
\begin{theorem}\label{theorem_6.1.1}
If $D_{u,\mkj}$ is surjective for all $(u,\mkj)\in
\wtm_\mkj(X,A)$, $\M_\mkj(X,A)$ is a smooth orbifold of dimension
$n(2-2g) + 2c_1(A)$.
\end{theorem}
The proof of smoothness is standard (\cite{?}). We will give the
proof in \S? using our terminology.

For $\mkj$ is unstable, the treatment is similar.

\subsection{Analytic set-up for $\M_{g,m}(X,A)$}\label{section_6.2}
For the analytic set-up, a general principle is to treat $M_{g,m}$
as a parameter space.  By this way,
 we can give a  family version of set-up. However, there
are some tedious issues.

Let us first assume $m=0$. For this case, there is no essential
change except that we replace $\mkj$ by $M_{g,m}$ (or $\tilde O$
if we emphasis the locality.) We summarize it:
\begin{itemize}
\item $\mc U_{u_o}$ is replaced by $\tilde O\times \mc U_{u_o}$;
\item $\tilde{\mc F}_\mkj$ is replaced by $\tilde O\times
\tilde{\mc F}_\mkj$; \item $\tilde{\mc E}_\mkj$ is replaced by the
parameterized bundle $\tilde{\mc E}_{\tilde O}$ which still can be
trivialized as $\tilde O\times \tilde{\mc E}_\mkj$.
\end{itemize}
We explain the last statement. The fiber of $\tilde{\mc E}_{\bar
O}$ over $(u,\mkj')$ is $L^p(\Lambda^{0,1}_{\mkj'}(u^\ast TX))$.
We explain the identification between
$$
L^p(\Lambda^{0,1}_{\mkj'}(u^\ast TX)) \leftrightarrow
L^p(\Lambda^{0,1}_{\mkj}(u_o^\ast TX)).
$$
First the identification between
$$
L^p(\Lambda^{0,1}_{\mkj'}(u^\ast TX)) \leftrightarrow
L^p(\Lambda^{0,1}_{\mkj'}(u_o^\ast TX))
$$
is given by $\Pi^{0,1}\circ P_\zeta$ as before; secondly, the
identification between
$$
L^p(\Lambda^{0,1}_{\mkj'}(u_o^\ast TM)) \leftrightarrow
L^p(\Lambda^{0,1}_{\mkj}(u_o^\ast TM))
$$
is induced by proposition \ref{prop_1.3.2}.

Next, we consider $m>0$. This case is subtle. On the one hand, we
can do the trivialization as what we do for $m=0$ case. But on the
other hand, we would like to trivialize bundles $\tilde{\mc E}$
and $\tilde{\mc F}$ in a different way. By using proposition
\ref{prop_3.3.1}, we may trivialize families $\tilde{\mc
E}_{\tilde O}$ and $\tilde{\mc F}_{\tilde O}$ locally. This is
necessary when we consider lower strata. However, these
trivialization
 causes  problems
technically at the first sight: as it is pointed out in \cite{R},
the family with such trivialization are not smooth. Namely, we
have  trivialization for both families $\tilde{\mc E}_{\tilde O}$
and $\tilde{\mc F}_{\tilde O}$ locally. But two different
trivialization $\tilde{\mc E}_{\tilde O}$ and $\tilde{\mc
E}_{\tilde O'}$ do not patch smoothly. However this trouble can be
solved by the following observation: first we note that trivialization is
patched well by restricting on smooth objects; secondly, by the
elliptic regularity property, all objects we are concerned are
smooth. Hence we may always assume that $\tilde{\mc E}_{\bar O}$
and $\tilde{\mc F}_{\bar O}$ are trivialized and study the theory
as if they are smooth families.

Hence,
\begin{theorem}\label{theorem_6.2.1}
If $D_{u,\mkj}$ is surjective for all $(u,\mkj)\in
\wtm_{g,m}(M,A)$, $\M_{g,m}(M,A)$ is a smooth orbifold of
dimension $n(2-2g) + 2c_1(A)+6g-6+2m$.
\end{theorem}
{\bf Proof. }Here $6g-6+2m$ is the dimension of parameter space
$M_{g,m}$. The theorem then follows from theorem
\ref{theorem_6.1.1}.

\subsection{Analytic set-up for $\M_{\strata}(X,A)$}\label{section_6.3}

Recall that
$$
\mkj_o:=(\Sigma,i_o, x_{o1}, \ldots, x_{om}).
$$
and
$$
\mk N(\mkj_o) =(\mk N(\mkj_{o1}),\ldots, \mk N(\mkj_{oc})).
$$
Let $\tilde\Sigma_{k}$ be the surface for $\mk N(\mkj_{ok})$.
Recall that
$$
\pi: \mk N(\mkj_o)\to \mkj_o.
$$
We define $ \chi^{1,p}_{\mk N(\mkj_o)}$ to be the set of elements
$$\tilde u:=(\tilde u_1,\ldots,\tilde u_c),
u_k\in \chi^{1,p}_{\mk N(\mkj_{ok})},
$$
such that it induces a continuous  map $u: \Sigma\to X$, i.e,
$u\circ \pi = \tilde u$. To avoid the complication of notations,
we simply use $\chi^{1,p}_{\mkj_o}$ for $ \chi^{1,p}_{\mk
N(\mkj_o)}$ and $u$ for $\tilde u$. This simplification only
causes a little ambiguity at nodal points. When this happens, we
always refer to the normalization of curves.

With $\mkj_o$ fixed, we still have $\tilde {\mc F}_{\mkj_o}$ and
$\tilde{\mc E}_{\mkj_o}$. Their fibers are given by the
followings. For $u_o\in \chi^{1,p}_{\mkj_o}$, its tangent space is
\begin{eqnarray*}
W^{1,p}(u_o^{\ast}TM)&=&\{(\zeta_1,...,\zeta_k)|
\zeta_v\in W^{1,p}(u_{ov}^\ast TM), \\
&&\zeta_v(y)=\zeta_w(y) \mbox{ for } y\in \Sigma_v \cap
\Sigma_w\}.
\end{eqnarray*}
This gives the fiber of $\tilde{\mc F}_{\mkj_o}$. Set
$$
L^p(\Lambda_\mkj^{0,1}(u_o^\ast TM)):=\bigoplus_vL^p
(\Lambda_{\mkj_{ov}}^{0,1}(u_{ov}^\ast TM)).
$$
This gives the fiber of $\tilde{\mc E}_{\mkj_o}$. We have the
linear operator
\begin{eqnarray*}
&&D_{\mkj_o,u_o}: W^{1,p}(u_o^*TM )\rightarrow
L^p(\Lambda^{0,1}_{\mkj_o}(u_o^\ast TM))\\
&&D_{\mkj_o,u_o}:=\left(D_{\mkj_{o1},u_{o1}},...,
D_{\mkj_{ok},u_{ok}}\right).
\end{eqnarray*}
\begin{lemma}\label{lemma_6.3.1}
 $D_{\mkj_o,u_o}$ is a Fredholm operator of
index $2c_1(A) + 2n(1-g)$.
\end{lemma}
\noindent For the proof see \cite{FO}.

To study $\M_{\strata}(X,A)$, we should allow that $\mkj$ varies.
Besides the similarities as above, there are parameters that
record the nodal points on each component. This is reflected in
the definition of $W^{1,p}(u_o^{\ast}TM)$. Therefore, we should
use the trivialization method mentioned at the end of last
subsection.

\begin{prop}\label{prop_6.3.1}
The stratum $\M_{\strata}(X,A)$ is a smooth orbifold of dimension
$n(2-2g) + 2c_1(A) + 6g - 6 +2m - 2|\sing|$, if $D_{i,u}$ is
surjective for any $(u,\mkj)\in \M_{\strata}$.
\end{prop}
\n{\bf Proof.} We verify the claim of dimension. For a stable
component, the moduli space of the component has dimension
$$
2c_1(A_k) + 2n(1-g_k) + 6g_k-6+2m_k +2s_k,
$$
where $m_k$ is the number of marked points and $s_k$ is the number
of nodal points (on the normalized surface); for an unstable
component, (only when $g=0$ and $m_k+ s_k\leq 2$), the dimension
is
$$
2c_1(A_k) + 2n -6+2m_k +2s_k.
$$
Totally we have
$$
2c_1(A) + 2n(c-\sum g) + 6\sum g-6c + 2m + 4|\sing|- 2n |\sing|
$$
Note that
$$
g= \sum g_v + \mr{rank} H^1(T); c- |\sing| = 1- \mr{rank} H^1(T).
$$
We have the formula of dimension.

For the smoothness, the proof is same as that of theorem
\ref{theorem_6.1.1}. We omit it.
 q.e.d.


\section{Coordinate charts for $\M_{g,m}(X,A)$}\label{section_7}

\subsection{Data of coordinate charts}\label{section_7.1}

We consider $\M_\mkj(X,A)$ with a fixed $\mkj\in M_{g,m}$. Let
$u\in \chi^{1,p}_\mkj$ and $\mc U_u$ be a neighborhood of $u$. We
may identify $\mc U_u$ with with an open set $W\subset
W^{1,p}(u^\ast TM)$ via $\exp_{u}$. Set
$$
L=L^p(\Lambda_\mkj^{1,0}(u^\ast TM)), \ M=\M_\mkj(X,A)\cap W.
$$

\begin{defn}\label{defn_7.1.1}
Let $W,L,M$ be as above. Suppose that we have
\begin{enumerate}
\item a smooth sub-manifold $U$ in $W$, \item a small open ball
$B_\delta\subset L$, a neighborhood $V$ of $u$ in $W$ and a
diffeomorphism
$$
\Phi:U\times B_\delta\to V,
$$
\item a smooth section
$$
f:U\to B_\delta
$$
\end{enumerate}
such that the map given by
$$
F:U\xrightarrow{(1,f)} U\times B_\delta\xrightarrow{\Phi} V
$$
maps $U$ onto $V\cap M$ and the map is diffeomorphic, we then call
$(U,\phi,F)$ (or $(U,\Phi,f)$, if no confusion may be caused,) a
data of coordinate chart.
\end{defn}
Obviously, by the definition $(U,F)$ gives a coordinate chart for
$M\cap V$.

\begin{prop}\label{prop_7.1.1}
Any two coordinate charts given by two different data are
$C^\infty$ compatible.
\end{prop}
{\bf Proof. } Suppose that we have two data of coordinate charts.
Be precise: we have $u_i,i=1,2$ and $\mc U_{u_i}$ which are
identified with $W_i$; then furthermore, we have
$(U_i,\Phi_i,F_i)$ which give coordinate charts $(U_i,F_i)$ for
$V_i\cap M$. So we have a transition map
$$
F_1\inv(F_1(U_1)\cap F_2(U_2))\xrightarrow{F_1} F_1(U_1)\cap
F_2(U_2)\xrightarrow{F_2\inv} F_2\inv(F_1(U_1)\cap F_2(U_2)).
$$
This map is the composition of the following chain:
$$
U_1\xrightarrow{(1,f_1)} U_1\times B_{\delta_1}
\xrightarrow{\Phi_1} V_1\xrightarrow{\Psi} V_2\xrightarrow{\Phi_2}
U_2\times B_{\delta_2}\xrightarrow{projection} U_2
$$
Here $\Psi=\exp_{u_2}\inv\exp_{u_1}$. Since each map in the chain
is smooth,  the transition map is smooth. q.e.d.

\v Note that in the proof, we use the fact that
$\chi^{1,p}_\mkj(X,A)$ is smooth. This is needed for the
smoothness of the map $\Psi$. However, if we consider
$\chi^{1,p}_{g,m}(X,A)$ and $\M_{g,m}(X,A)$, the fact is not true.
The problem can be solved by a small modification:
\begin{remark}\label{rmk_7.1.1}
We modify   the definition by requiring that $U$ consists of
smooth maps. Then we may repeat the argument of proposition
\ref{prop_7.1.1} for the $\M_{g,m}(X,A)$. The only problem is
$\Psi$. Although $\Psi$ is not smooth in general, it is smooth
when restricted on smooth maps.
\end{remark}

\subsection{Proof of theorem \ref{theorem_6.2.1}}\label{section_7.2}
We only prove the smooth structure of $\M_\mkj(X,A)$. The proof
for that of $\M_{g,m}(X,A)$ is similar.

The goal is to construct a data of coordinate chart for each point
$u\in \M_{\mkj}(X,A)$. Set
$$
W=W^{1,p}(u^\ast TM), L=L^p(\Lambda_\mkj^{0,1}(u^\ast TM)).
$$
By our assumption,
$$
D_{u,\mkj}: W\to L
$$
is surjective. Hence we may construct a right inverse $Q_{u,\mkj}$
to $D_{u,\mkj}$ such that $Q_{u,\mkj}$ is $\aut(u,\mkj)$
equivariant: note that a right inverse gives a splitting
$$
W=\ker D_{u,\mkj} \oplus \mathrm{range} Q_{u,\mkj}
$$
and vice versa. We choose $\mathrm{Q_{u,\mkj}}$ to be $(\ker
D_{u,\mkj})^\perp$ with respect to $L^2$-norm. Since $\ker
D_{u,\mkj}$ and $L^2$-norm are $\aut(u,\mkj)$ invariant, the
splitting is $\aut(u,\mkj)$-equivariant.

Now we define
\begin{eqnarray*}
&&
\Phi: \ker D_{u,\mkj}\times L\to W; \\
&& \Phi(\xi, \eta)=\xi +Q_{u,\mkj} \eta.
\end{eqnarray*}
Then, there exists a small neighborhood $U$ of $0\in \ker
D_{u,\mkj}$, a small ball $B_\delta\subset L$ and a neighborhood
$V$ of $u$ in $W$ such that
$$
\Phi: U\times B_\delta\to V
$$
is diffeomorphic.

It remains to construct a section
$$
f: U\to B_\delta
$$
such that
\begin{equation}\label{eqn_7.2.1}
\dbar(\Phi(\xi,f(\xi)))=0
\end{equation}
for any $\xi\in U$. For this purpose, we consider the map
\begin{eqnarray*}
&& H: U\times B_\delta\to U\times L; \\
&& H(\xi, \eta)= (\xi,\dbar(\Phi(\xi,\eta)).
\end{eqnarray*}
Then
$$
H(0,0)=(0,0); dH_{(0,0)}=id.
$$
By the inverse function theorem, there is a smooth section solving
\eqref{eqn_7.2.1}. This completes the proof. q.e.d.

\subsection{Constructing data of coordinate charts}\label{section_7.3}

Again, we only consider $\M_\mkj(X,A)$. The situation is: let $U$
be a smooth sub-manifold of $\chi^{1,p}_\mkj(X,A)$; fix a point
 $u_o\in
U$; set $W$ to be a small neighborhood of $u_o\in W^{1,p}(u_o^\ast
TM)$ and $L=L^{1,p}(\Lambda^{0,1}_{\mkj}(u_o^\ast TM)$; let
$$
\mc Q=\{Q_{u,\mkj}|u\in U\}
$$
be a smooth family of right inverses for $u\in U$. Then we define
$$
\Phi(u,\eta) = u+ Q_{u,\mkj}\eta.
$$
Furthermore, we have the following assumption on $(U,\mc Q)$:
\begin{assumption}\label{assumption_7.3.1}
Let $(U,\mc Q)$ be as above with  properties
\begin{enumerate}
\item $\|\nabla u\|_{L^p}\leq C$ for any $u\in U$; \item for any
$u\in U$
$$
\|\dbar u\|_{L^p}\leq \epsilon;
$$
\item for any $\zeta\in T_uU$
$$
\|\frac{d\dbar u}{d\zeta}\|_{L^p}\leq \epsilon \|\zeta\|;
$$
\item for right inverses
$$
\|Q_{u,\mkj}\|\leq C
$$
and
$$
\|Q_{u_1,\mkj}-Q_{u_2,\mkj}\|\leq C\|u_1-u_2\|_{L^{1,p}}.
$$
\end{enumerate}
Here $C$ is a constant and $\epsilon$ is a small constant such
that  $C\epsilon \ll 1$.
\end{assumption}
For any $(U,\mc Q)$ satisfying the assumption, we explain that we
may produce a data of coordinate chart from it for a neighborhood
of $u_o$.

Applying the famous Taubes argument, we have
\begin{prop}\label{prop_7.3.1}
There exists a smooth map
$$
f: U\to B_\delta
$$
such that $u+ Q_u f(u)$ is holomorphic. Any holomorphic curve in
the form $u+ Q_u\xi, \xi\in B_\delta$ is given by $\xi=f(u)$. Here
$\delta$ is a small number that depends only on $C$. Moreover
\begin{equation}\label{eqn_7.3.1}
\|f(u)\|_{L^p}\leq 2\epsilon.
\end{equation}
\end{prop}
We remark that we may assume that $\epsilon \ll \delta\ll C$.

\v\n {\bf Proof. } Composing with
$$
\Phi: U\times L\to W,
$$
we have a family of operators parameterized by $u\in U$:
$$
\bar\partial: U\times L\xto{\Phi} W\xto{\dbar } L.
$$
Be precise, for each $u$, we have
$$
\bar\partial(u,\cdot) : L\to L; \bar\partial(u,\eta)=\dbar(u+
Q_u\eta).
$$
We now solve $\eta$ for the equation
$$
\dbar(u+Q_u\eta)=0.
$$
Expand  the equation we have
$$
\dbar(u+Q_u\eta)= \dbar u + D_uQ_u\eta + N_u(Q_u\eta) = \cdots
+\eta+ \cdots.
$$
Here $N_u(Q_u\eta)$ is a term with second or higher order. We use
the fact
$$
\|N_u(\xi_1)-N_u(\xi_2)\|_{L^p} \leq
C_0(\|\xi_1\|+\|\xi_2\|)(\|\xi_1-\xi_2\|).
$$
Here $C_0$ depends only on $\|\nabla u\|_{L^p}$.

The equation to solve is
$$
\eta = -\dbar u -N_u(Q_u\eta).
$$
Let $H: B_\delta\to B_\delta$ be a map defined by
$$
H\eta = -\dbar u-N_u(Q_u\eta).
$$
By choosing proper $\delta$, $H$ is a contraction map. This
follows by two simple estimates.
\begin{eqnarray*}
\|H\eta\| &\leq & \|\dbar u\| + \|N_u(Q_u\eta)\|\\
&\leq & \epsilon + C_0C^2 \|\eta\|^2\\
&\leq & \epsilon + C_0C^2\delta^2\\
&\leq & \epsilon + \delta/4 \leq \delta;
\end{eqnarray*}
here we require that $C_0C^2\delta<1/4$ and $\epsilon\ll \delta$;
\begin{eqnarray*}
\|H\eta_1-H\eta_2\|_{L^p} &=&
\|N_u(Q_u\eta_1)-N_u(Q_u\eta_2)\|\\
&\leq& C_0(\|Q_u\eta_1\|+\|Q_u\eta_2\|)\|Q_u(\eta_1-\eta_2)\|\\
&\leq& 2C_0C^2\delta \|\eta_1-\eta_2\|\\
&\leq& 0.5 \|\eta_1-\eta_2\|.
\end{eqnarray*}
We conclude that $H$ is a contraction map. On the other hand, we
can also show that $H:B_{2\epsilon}\to B_{2\epsilon}$ is a
contraction map. This implies the estimate for $f(u)=\eta$. q.e.d.
\v\n In this proposition, we essentially only use the property
(1) in Assumption \ref{assumption_7.3.1}.
\begin{theorem}\label{theorem_7.3.1}
There exists a small neighborhood $U'\subset U$ of $u_o$,
$\delta_1\leq \delta$ and $V\subset W$ such that
$$
\Phi: U'\times B_{\delta_1}\to V
$$
is diffeomorphic. Here $\delta_1$ depends only on $C$.
\end{theorem}
{\bf Proof. } We may identify $ W^{1,p}(u_o^\ast TM)$ with $ \ker
D_{u_o,\mkj}\oplus L^p$ via
$$
\xi+Q_{u_o,\mkj}\eta\leftrightarrow (\xi,\eta).
$$
We rewrite map $\Phi$ as
\begin{eqnarray*}
&&\Phi: U\times B_\delta \to W^{1,p}(u_o^\ast TM)= \ker
D_{u_o,i}\oplus L^p;\\
&& \Phi(u, \eta)= (\bar u+Q_u\eta-Q_{u_o}D_{u_o}(\bar u+Q_u\eta),
D_{u_o}(\bar u+ Q_u\eta)),
\end{eqnarray*}
Here $\bar u=u-u_o$. The tangent map  of $\Phi$ at $u,\eta$ is
$$
D\Phi_{u,\eta}(\xi,\zeta) =\left(
\begin{array}{ll}
\xi + I_{11} & I_{12}\\
I_{21} & \zeta +I_{22}
\end{array}
\right),
$$
where
\begin{eqnarray*}
&& I_{11}= \frac{dQ_u}{d\xi}\eta -Q_{u_o}D_{u_o}\xi-
Q_{u_o}D_{u_o}\frac{dQ_u}{d\xi}\eta=: I_{111}+I_{112}+I_{113};\\
&& I_{12}=Q_u\zeta -Q_{u_o}D_{u_o}Q_u\zeta;\\
&& I_{21}=D_{u_o}(\xi + \frac{dQ_u}{d\xi}\eta);\\
&& I_{22}=D_{u_o}Q_{u_o}\zeta-\zeta.
\end{eqnarray*}
By direct estimates, we have that for proper chosen $U'\subset
U,\delta'<\delta$ and $(u,\eta)\in U'\times B_{\delta'}$,
$$
\|I_{ij}\|\leq \frac{1}{100}\|(\xi,\zeta)\|
$$
Hence $D\Phi_{u,\eta}$ is invertible and
$$
\|D\Phi_{u,\eta}\|\leq 2
$$
for $(u,\eta)\in U'\times B_{\delta'}$.

Finally, we show that $\Phi$ is injective. Suppose that
$$
\Phi(u_1,\eta_1)= \Phi(u_2,\eta_2).
$$
In general, we have
$$
\Phi(u,\eta)= \Phi(u_o,0) + D\Phi_{u_o,o}(\bar u,\eta) + N(\bar
u,\eta).
$$
Here
$$
N(\bar u,\eta)= (Q_u\eta- Q_{u_o}D_{u_o}Q_u\eta,
D_{u_o}Q_u\eta-\eta).
$$
It is not hard to get
\begin{equation}\label{eqn_7.3.2}
\|N(\bar u_1,\eta_1)-N(\bar u_2,\eta_2)\| \leq C(\|(\bar
u_1,\eta_1)\|+\|(\bar u_2,\eta_2)\|)(\|(u_1-u_2,
\eta_1-\eta_2)\|).
\end{equation}
We have
$$
D\Phi_{u_o,0}((\bar u_1,\eta_1)-(\bar u_2,\eta_2)) =-(N(\bar
u_1,\eta_1)-N(\bar u_2,\eta_2)).
$$
Set $h=\|(u_1,\eta_1)-(u_2,\eta_2)\|$, then
$$
\frac{h}{2}\leq C(\|(\bar u_1,\eta_1)\|+\|(\bar u_2,\eta_2)\|)h.
$$
This is impossible if $\|(\bar u_v,\eta_v)\|,v=1,2,$ are small.
Here $\Phi$ is injective. q.e.d.

\v\n As a corollary, $(U,\Phi)$ yields a data of coordinate chart
$(U',\Phi,F)$ .

\subsection{Estimates of $df/d\xi$}\label{section_7.4}
Finally, we discuss the derivative
$$
\frac{df}{d\xi},\xi\in T_uU.
$$
We show that
\begin{theorem}\label{theorem_7.4.1}
Let $f$ be constructed in proposition \ref{prop_7.3.1}. Then
\begin{equation}\label{eqn_7.4.1}
\|\frac{df}{d\xi}\|\leq C\epsilon \|\xi\|.
\end{equation}
\end{theorem}
{\bf Proof. } The proof is rather long although it is
straightforward.

Let $u_t$ be a path with $u_0=u_o$ and representing $\xi\in
T_{u_o}U$. We differentiate the equation
$$
\dbar u_t + f(u_t) + N_{u_t}(Q_{u_t}f(u_t))=0
$$
and get
\begin{eqnarray*}
0 &=&\frac{d}{d \xi}\dbar u
+ \frac{d f(u) }{d\xi} + \frac{d(N_uQ_uf(u))}{d\xi}\\
&=&\frac{d}{d \xi}\dbar u
+ \frac{d f(u) }{d\xi} \\
&&+\frac{dN_u(Q_{u_o}f(u_o))}{d\xi}
+N_{u_o}(\frac{d(Q_uf(u))}{d\xi})
\\
&=& I_1+I_2+I_4+I_4.
\end{eqnarray*}
We have
$$
\|I_1\|\leq \epsilon \|\xi\|
$$
by property (2) in assumption \ref{assumption_7.3.1}.

To get the estimate for $I_4$  we consider
\begin{eqnarray*}
&&N_{u_o}(Q_{u_t}f(u_t)-Q_{u_o}f(u_o))\\
 &&\leq  C_0(\|(Q_{u_t}f(u_t)\|+\|Q_{u_o}f(u_o)\|)
(\|Q_{u_t}f(u_t)-Q_{u_o}f(u_o)\|) \\
&&\leq  2C_0 \|Q_{u_o}f(u_o)\|(C\|u_t-u_o\|\|f(u_o)\|
+C\|f(u_t)-f(u_o)\|)
\end{eqnarray*}
which says that
$$
I_4\leq 2C_0C^2\|f(u_o)\|^2\|\xi\|+ 2C_0C^2\|f(u_o)\| \|I_2\|\leq
C\epsilon^2\|\xi\|+ 0.5\|I_2\|.
$$
The estimate
$$
\|I_3\|\leq \epsilon\|\xi\|
$$
is given in the next lemma. Combine all these together, we have
$$
\|I_2\|\leq C\epsilon \|\xi\|.
$$
q.e.d.

\v\n
\begin{prop}Let $u_t$, $N_u$ and $\xi$ be as above, then
$$
\|\frac{dN_u(\eta)}{d\xi}\|\leq C\|\xi\|_{L^{1,p}}(\epsilon
+\|\eta\|_{L^{1,p}}).
$$
\end{prop}
{\bf Proof. } As we know
$$
\dbar(u_t+\eta)= \dbar u_t +D_{u_t}\eta +N_{u_t}\eta.
$$
On the other hand,
$$
\dbar(u_t+\eta)= \dbar u_o+D_{u_o}(\bar u_t+\eta) + N_{u_o}(\bar
u_t+\eta),
$$
where $\bar u=u-u_o$. Set two right hand sides equal. Then
\begin{eqnarray*}
\frac{\bar\partial_{J,i}u_t-\bar\partial_{J,i}u_o}{t}
&+&\frac{D_{u_t}\eta-
D_{u_o}(\bar u_t+\eta)}{t}\\
&+& \frac{N_{u_t}(\eta)- N_{u_o}(\bar u_t+\eta)}{t}=0
\end{eqnarray*}
 By taking $t\to 0$, we have
\begin{eqnarray*}
\frac{\bar\partial_{J,i}u_t-\bar\partial_{J,i}u_o}{t}
&\to& \frac{d}{d\xi}(\bar\partial_{J,i}u);\\
 \frac{D_{u_t}\eta-
D_{u_o}(\bar u_t+\eta)}{t}&\to&
\frac{d}{d\xi}(D_{u})\eta-D_{u_o}\xi;
\end{eqnarray*}
while for
\begin{eqnarray*}
\frac{N_{u_t}(\eta)- N_{u_o}(\bar u_t+\eta)}{t} &=&
\frac{N_{u_t}(\eta)- N_{u_o}(\eta)}{t} \\
&&+\frac{N_{u_o}(\eta)- N_{u_o}(\bar u_t+\eta)}{t},
\end{eqnarray*}
its limit is
\begin{equation}\label{eqn_7.4.2}
\frac{d}{d\xi}N_{u_o}(\eta)+\lim_{t\to 0}\frac{N_{u_o}(\eta)-
N_{u_o}(\bar u_t+\eta)}{t}.
\end{equation}
Therefore
\begin{eqnarray*}
\frac{d}{d\xi}N_{u_o}(\eta)&=&
-\frac{d}{d\xi}(\bar\partial_{J,i}u)
-\frac{d}{d\xi}D_{u}\eta+D_{u_o}\xi\\
&&\lim_{t\to 0}\frac{N_{u_o}(\eta)- N_{u_o}(\bar u_t+\eta)}{t}.
\\
&=:& I_1+I_2+I_3+I_4.
\end{eqnarray*}
For each term we have
\begin{eqnarray*}
\|I_1\|_{L^p} &\leq & \epsilon\|\xi\|_{L^{1,p}}, \\
\|I_2\|_{L^p} &\leq & C\|\xi\|_{L^{1,p}}\|\eta\|_{L^{1,p}},\\
\|I_3\|_{L^p} &\leq & \epsilon\|\xi\|_{L^{1,p}},\\
\|I_4\|_{L^p} &\leq & C\|\xi\|_{L^{1,p}}\|\eta\|_{L^{1,p}}.
\end{eqnarray*}
The estimate of $I_4$ follows from lemma \ref{ap1_b}. q.e.d.


\section{Balanced $J$-holomorphic curves}\label{section_8}
We consider the moduli space $\M_{g,m}(X,A)$ with $2g+m\leq 2$.
There are 4 cases: $(g,m)=(0,0),(0,1),(0,2)$ and $(1,0)$. In this
section, we focus on $(g,m)=(0,1)$ and $(0,2)$ since we need them
when consider gluing.

Let
$$
\mkj_m =(S^2, j,x_1,\ldots, x_m),1\leq m\leq 2.
$$
The moduli spaces are
$$
\M_{0,m}(X,A)=\frac{\wtm_{0,m}(X,A)}{\aut(\mkj_m)},
$$
where $\wtm_{0,m}(X,A)$ is defined below.

Since $\aut(\mkj_m)$ is a non-compact finite dimensional Lie
group, it is useful to construct the slice for the quotient space,
or  reduce the quotient group to be compact. For this purpose, we
introduce balanced holomorphic maps.

\v\n{\bf Case 1, $(g,m)=(0,1)$.}

\v\n $\M_{0,1}$ consists of only one element $\mkj_1=(S^2,
\infty)$. Here $S^2-\infty = \cplane$. We use $\cplane$ in our
discussion in this subsection. Let $\mk{t}=\cplane$ be the group
of translations of $\cplane$ and $\mk{m}=\cplane^\ast$ that acts
on $\cplane$ by multiplications. The semi-product
$\mk{B}=\mk{t}\ltimes \mk{m}$ acts on $\cplane$ as
$$
(t,m)\cdot z= m(z-t).
$$
It is well known that
$$
\aut(\mkj_1)= \mk{B}.
$$
Let
$$
\wtm_{0,1}(X,A):=\wtm_{0,0}(X,A):=\{ u:S^2\to X| \dbar u=0,
[u(S^2)] =A \}.
$$
Then
$$
\M_{0,1}(X,A)=\frac{ \wtm_{0,1}(X,A)}{\mk B}.
$$
For $u\in \wtm_{0,0}(X,A)$ we usually call $|du|^2$ the energy
density.  Note that the energy of $u$ is $\omega(A)$. Let $\hbar=
\omega(A)/2$.
\begin{defn}\label{defn_8.1}
A $J$-curve  $u\in \wtm_{0,1}(X,A)$
 is called balanced if
\begin{itemize}
\item the energy center of $u$ is $0\in \cplane$; \item the energy
on the unit disk is $\hbar$.
\end{itemize}
Let $\M^b_{0,1}(X,A)$ be the space of balanced $J$-curves.
\end{defn}
We remark that for any $u\in \wtm_{0,1}(X,A)$ there is a canonical
balanced curve $\mk b_1(u)$ constructed
\begin{itemize}
\item by translating the energy center of $u$ to 0; \item by
proper dilation (i.e, multiplying a proper real number) such that
the energy on the unit disk is $\hbar$.
\end{itemize}
It is then easy to see that
\begin{equation}\label{eqn_8.1}
\M_{0,1}(X,A)= \frac{\M^b_{0,1}(X,A)}{S^1}.
\end{equation}
Here $S^1$ acts on $\cplane $ by rotations and therefore has an
induced action on $\M^b_{0,1}(X,A)$. When we consider
$\M_{0,1}(X,A)$ we always use \eqref{eqn_8.1}.

\v\n{\bf Case 2, $(g,m)=(0,2)$. }

\v\n This case is similar but easier. $\M_{0,2}$ consists only an
element $\mkj_2=(S^2, 0,\infty)$. Then
$$
\aut(\mkj_2) =\mk{m}=\cplane^\ast.
$$
Set $\wtm_{0,2}(M,A)=\wtm_{0,0}(M,A)$.
\begin{defn}\label{defn_8.2}
A $J$-curve $u \in \wtm_{0,2}(M,A)$
 is called balanced if the energy of $u$ on the unit disk is
 $\hbar$. Let $\M^b_{0,2}(X,A)$ be the space of balanced $J$-curves.
\end{defn}
We also have
\begin{equation}\label{eqn_8.2}
\M_{0,2}(X,A)= \frac{\M^b_{0,2}(X,A)}{S^1}.
\end{equation}

\vskip 0.2in
\begin{center}
{\bf Part III. The Gluing Theory}
\end{center}


\section{Gluing maps }\label{section_9}
In \S\ref{section_9}--\S\ref{section_12}, we discuss the basic
case, i.e, the gluing theory for 1-nodal strata. Then we
generalize it to general  strata in \S\ref{section_13}.

\subsection{Pre-gluing}\label{section_9.1}
Let $\strata=(\mk g,\mk A, T,D)$ be a data of stratum in
$\om_{g,m}(X,A)$. For simplicity, we assume $m=0$.  Here
$$
\mk g=\{ g_1,g_2\}, \mk A=\{A_1,A_2\}
$$
and $T$ consists of two vertices $v_1,v_2$ and one edge $e$.
 $D$ is trivial since $m=0$.

Set $S=(\mk g,T,D)$. Let $\mkj_o\in M_S$ and $(u_o,\mkj_o)\in
\M_{\mkj_o}(X,A)$. Suppose that $\mkj_o$ consists of
$$
\mkj_{ov}=(\Sigma_v, j_{ov}, y_{ov}), v=1,2.
$$
By identifying $y_{o1}$ and $y_{o2}$, we get
$\mkj_o=(\Sigma,i_o)$. We  write
$$
\Sigma =\Sigma_1\cup_{y_{o1}=y_{o2}} \Sigma_2.
$$
We denote the singular point by $y_o$. $u_o$ consists of
$J$-holomorphic curves
$$
u_{ov}: \Sigma_v\to M, \ [u_v(\Sigma_v)]=A_v
$$
with $u_{o1}(y_{o1})=u_{o2}(y_{o2})$.

Recall that we have an (orbi-)line bundle
$$
 L_S\to M_S.
$$
The forgetting-map map
\begin{eqnarray*}
&& \mk f: \M_\strata(X,A)\to M_S; \\
&& \mk f(u,\mkj)=\mkj
\end{eqnarray*}
induces an orbi-line bundle
$$
\mc L_\strata=\mk f^\ast L_S\to \M_\strata(X,A).
$$
Given a point $\mk p\in \mc L_\strata$, our goal is to construct a
holomorphic map $Gl(\mk p)\in \M_{g,m}(X,A)$. Put in the local
coordinate, we write $\mk p=(u_o,\mkj_o,\rho), \rho=re^{i\theta}$,
we construct $Gl(u_o,\mkj_o,\rho)$. The first step of the
construction is pre-gluing, which gives an approximation
holomorphic map $\pgl(u_o,\mkj_o,\rho)$.

Recall that we have a gluing map for surfaces:
$$
\gs: L_S\to M_{S_0}.
$$
In local coordinates, we write
$$
\mkj_{o\rho}= \gs(\mkj_o,\rho)=(\Sigma_{\rho,y_o},j_{o\rho}).
$$
Geometrically, $\Sigma_{\rho,y_o}$ is obtained as the following.
We use the holomorphic cylindrical coordinates $(\log s_i,t_i)$ on
$\Sigma_i$ near $y$, and write
$$\Sigma_2-\{y_{o2}\}=\Sigma_{20}\bigcup\{[0,\infty)\times S^1\},$$
$$\Sigma_1-\{y_{o1}\}=\Sigma_{10}\bigcup\{(-\infty,0]\times S^1\}.$$
We cut off the part of $\Sigma_i$ with cylindrical coordinate
 glue the remainders by identifying the
$|\log r|$-long ends of the cylinders with a twist of angle
$\theta$. The new curve is $\mkj_{o\rho}$. $\pgl(u_o,\mkj_o,\rho)$
is  expected to be a map on $\mkj_{o\rho}$.

More generally, we may replace holomorphic map $u_o$ by $u\in
\chi^{1,p}_{\mkj_o}(X,A)$. Write $\phi= \pgl(u,\mkj_o,\rho)$ where
$u=(u_1,u_2)$. $\phi$ is supposed to be a map on surface
$\Sigma_{\rho,y_o}$. Define
$$
\phi(x)=\left\{
\begin{array}{ll}
u_1(x) & \textup{ if } x\in \Sigma_1 - D_{y_{o1}}(2r^{1/4})
\\
p=u_1(y_{o1})=u_2(y_{o2})  & \mbox{ if } x\in  D_{y_{o1}}(r^{1/4})- D_{y_{o1}}(r^{3/4})\\
u_2(x) & \textup{ if } x\in \Sigma_2 - D_{y_{o2}}(2 r^{1/4})
\end{array}\right.
$$
To define the map in the rest part we fix a smooth cutoff function
cutoff function $\beta : \rone\rightarrow [0,1]$ such that
\[
\beta (s)=\left\{
\begin{array}{ll}
1 & if\;\; s \geq 2 \\
0 & if\;\; s \leq 1
\end{array}
\right.
\]
and $|\beta^{\prime}(s)|\leq 2.$ We assume that $r$ is small
enough such that $u_i$ maps the disk $D_{y_i}(4r^{1/4})$ into a
normal coordinate domain of $p$. We can define $\phi$ by
$$
\phi(x) =\exp_p\left(\beta\left(\frac{x}{r^{1/4}}\right)
           \exp_p^{-1}u_1(x)+ \beta\left(\frac{r^{1/4}}{z}\right)
           \exp_p^{-1}u_2\left(\frac{\rho}{x}\right) \right).
$$

\begin{lemma}\label{lemma_9.1.1}
Suppose $\phi=\pgl(u,\mkj_o,\rho)$, then
$$
\|\bar{\partial}_{J, \mkj_{o\rho}}\phi\|_{L^p} \leq
\|\bar{\partial}_{J,\mkj_o} u\|_{L^p} +Cr^{\frac{1}{2p}},
$$
where $C$ is independent of $\rho$. In particular,
$$
\|\bar{\partial}_{J, \mkj_{o\rho}}\phi\|_{L^p} \leq
Cr^{\frac{1}{2p}}
$$
if $u$ is holomorphic.
\end{lemma}
The proof is given in \S\ref{section_10}.

For $u=u_o$, set
$$
\phi_o= \pgl(u_o,\mkj_o,\rho).
$$

\subsection{Right inverses}\label{section_9.2}

Let $u\in \chi^{1,p}_{\mkj_o}$. We assume that $D_{u,\mkj_o}$ is
surjective.
 Therefore, there is a right inverse
$$
Q_{u,\mkj_o}: L^p(\Lambda^{0,1}_{\mkj_o}(u^\ast TM)) \to
W^{1,p}(u^\ast TM).
$$
with $\|Q_{u,\mkj_o}\|\leq C$. Let $\phi=\pgl(u,\mkj_o,\rho)$. We
construct the right inverse to $D_{\phi,\mkj_{o\rho}}$.

We identify
$$
\Sigma''_v:= \Sigma_{v0}\cup \{(\log r,0)\times S^1\}, v=1,2
$$
with $\Sigma_{\rho}$ in an obvious way. We introduce two pairs
$\lambda_{v,\rho}$ and $\gamma_{v,\rho},v=1,2,$ of cut-off
functions
 on them.  We only describe these functions on the
cylinder ends only since they are 1 on $\Sigma_{v0}$. Let
$$
\lambda_{v,\rho}(t,\theta)= \left\{
\begin{array}{ll}
1,& \mbox{ if } t>\log r/2 +1;\\
0,& \mbox{ if } t<\log r/2 -1,
\end{array}
\right.
$$
with  $\lambda_{1,\rho}+\lambda_{2,\rho}=1$. Let
$$
\gamma_{v,\rho}(t,\theta)= \left\{
\begin{array}{ll}
1,& \mbox{ if } t>\log r/2-1;\\
0,& \mbox{ if } t<\log r,
\end{array}
\right.
$$
Note that $\gamma_{v,\rho}$ is 1 on the support of
$\lambda_{v,\rho}$. Also
$$
|\nabla \gamma_\ast|\leq \frac{C}{|\log r|}.
$$

Suppose that  $\eta$  is a function (or a form) on $\Sigma_\rho$.
We define
$$
\Lambda(\eta)=\lambda_{1,\rho}\eta\dotplus\lambda_{2,\rho}\eta
$$
to be a function (or a form) on $\Sigma$. Note that
$\lambda_{v,\rho}\eta, v=1,2$ are functions (or forms) on
$\Sigma''_{v}\subset \Sigma_v\subset \Sigma$. By $\dotplus$, we
mean the sum is taken over $\Sigma$.

Conversely, suppose $\sigma$ is a continuous function (or form) on
$\Sigma$. Define
$$
\sigma_1(x)= \gamma_{1,\rho}(\sigma(x)-\sigma(y))+\sigma(y), x\in
D_{y_{o1}}(\sqrt r)
$$
and equals to $\sigma$ outside the disk. $\sigma_1$ is a function
on $\Sigma_1''\subset \Sigma_\rho$. Similarly, we have a function
$\sigma_2$ on $\Sigma_2''\subset \Sigma_\rho$. Define
$$
\Gamma(\sigma)= \sigma_1\dotplus\sigma_2
$$
to be a function on $\Sigma_\rho$. Here by $\dotplus$, we mean the
sum is taken over $\Sigma_\rho$.

\begin{lemma}\label{lemma_9.2.1}
For $\eta\in L^p(\Lambda_{\mkj_{o\rho}}^{0,1}\phi^\ast TM)$
$$
\|D_{ \phi,\mkj_{o\rho}}R\eta- \eta\|_{L^p}\leq \frac{C}{|\log r|}
\|\eta\|_{L^{p}},
$$
where $R=\Gamma Q_{u,\mkj_o}\Lambda(\eta)$.
\end{lemma}
The proof is given in \S\ref{section_10}. The lemma says that
$D_{\phi,\mkj_{o\rho}}R$ is invertible. Set
$$
Q_{\phi,\mkj_{o\rho}}= R(D_{\phi,\mkj_{o\rho}}R)\inv.
$$
\begin{prop}\label{prop_9.2.1}
$Q_{\phi,\mkj_{o\rho}}$ is a right inverse to
$D_{\phi,\mkj_{o\rho}}$. Moreover
$$
\|Q_{\phi,\mkj_{o\rho}}\|\leq C
$$
where $C$ is independent of $\rho$.
\end{prop}
In particular, for $\phi_o$ we construct the right inverse
$Q_{\phi_o,\mkj_{o\rho}}$.

\subsection{Gluing maps}\label{section_9.3}

With $\phi_o$ and $Q_{\phi_o,\mkj_{o\rho}}$, we can construct a
holomorphic curve as in proposition \ref{prop_7.3.1}.

We need the lemma
\begin{lemma}\label{lemma_9.3.1}
Let $\phi=\pgl(u, \mkj_o,\rho)$.
$$
\|N_{\phi}(\zeta_1)-N_{\phi}(\zeta_2)\|_{L^p} \leq
C(\|\zeta_1\|_{L^{1,p}}+\|\zeta_2\|_{L^{1,p}})(\|\zeta_1-\zeta_2\|_{L^{1,p}}),
$$
where $C$ depends only on $\|u\|_{L^{1,p}}$.
\end{lemma}
{\bf Proof. } By theorem \ref{est_N}, we have this inequality with
some constant $C'$ depending on $\|\phi\|_{L^{1,p}}$. By the
construction of $\phi$, we know that
$$
\|\phi\|_{L^{1,p}}\leq C''\|u\|_{L^{1,p}}
$$
So the claim follows. q.e.d.

\v\n
\begin{theorem}\label{theorem_9.3.1}
Suppose that  $\phi$ is as above and let $C_0$ be the constant
given in lemma \ref{lemma_9.3.1}. Suppose that
$$
\|\bar\partial_{J,\mkj_{o\rho}}\phi\|_{L^p}\leq \epsilon
$$
for some $\epsilon\ll  C_0\inv$. Then  in the $\delta$-ball of
$L^p(\Lambda_{\mkj_{o\rho}}^{0,1}\phi^\ast TM)$ with $\delta C_0
<1/2$, there exists a unique element, denoted by $f(u,
\mkj_o,\rho)$,  such that
$$
\exp_{\phi} Q_{\phi,\mkj_{o\rho}}f(u,\mkj_o,\rho)
$$
is $J$-holomorphic and
$$
\|f(u, \mkj_o,\rho)\|\leq C\epsilon,
$$
where $C$ can be any constant such that $CC_0\epsilon <1/2$.
\end{theorem}
The proof is a repeat of that in proposition \ref{prop_7.3.1}.

\begin{remark}\label{rmk_9.3.1}
Since we are working on some spaces with orbifold structure, we
should require that the gluing maps are equivariant with respect
to isotropic groups.

Let $(u_o,\mkj_o)\in \M_\strata(X,A)$, the local uniformization
system for a neighborhood $O$ of $(u_o,\mkj_o)$ in
$\M_\strata(X,A)$
 and bundle $\mc L_\strata|_O$
are in the form
$$
(\tilde O,\aut(u_o,\mkj_o),\pi) \mbox{ and } (\tilde{\mc
L}_\strata|_{\tilde O}, \aut(u_o,\mkj_o),\pi).
$$
The gluing map is, at the moment, defined on $\tilde{\mc
L}|_{\tilde O}$ other than on $\mc L|_O$. Then  we note that
\begin{enumerate}
\item when $\mkj_o$ is pre-stable, the gluing is
$\aut(u_o,\mkj_o)$-equivariant. Hence the gluing is defined on
$\mc L_\strata$; \item when $\mkj_o$ is not pre-stable, there is
at least one non-pre-stable component $\mkj_{ov},v=1,2$. The
component is of $g=0,1\leq m\leq 2$. For this case, we have to
use  moduli spaces of balanced curves. Then it is easy to see that
the gluing is well defined on $\mc L_\strata$.
\end{enumerate}
\end{remark}

Let $U\subset \M_\strata(X,A)$ be any proper open subset. Define
the gluing map to be
\begin{eqnarray*}
&& Gl: \mc L_{\strata,\epsilon}^0|_U\to \M_{g,m}(X,A)\\
&& Gl(u,\mkj,\rho) = \pgl(u,\mkj,\rho)+ f(u,\mkj,\rho).
\end{eqnarray*}
Here $\epsilon$ depends only on $U$. To stress the process of
gluing, we set
$$
\pert(u,\mkj,\rho)=f(u,\mkj,\rho).
$$
Here $\pert$ stands for perturbation which is exactly what we are
doing in the second step.

\subsection{Gluing maps for general strata}\label{section_9.4}

Now suppose that  $\strata=(\mk g,\mk A,T,D)$ is any stratum and
$S=(\mk g, T,D)$. For simplicity we assume that $\M_\strata(X,A)$
is compact, otherwise we always restrict our discussion on a
proper open subset in the stratum.

Recall that for any $S\prec S'$ (and correspondingly
$\strata\prec\strata'$) there exists a gluing bundle $\mc
L_{\strata,\strata'}$. Repeat the process in
\S\ref{section_9.1}-\S\ref{section_9.3}, we have a gluing map
$$
Gl_{\strata,\strata'}: \mc L_{\strata,\strata,\epsilon}^0 \to
\M_{\strata'}(X,A).
$$
Now consider a point $\mk p\in \mc L_{\strata}$
$$
\mk p= (u, \mkj, \rho_1,\rho_2),
$$
where $\rho_1$ denotes the coordinate corresponding to the fiber
in $\mc L_{\strata,\strata'}$ and $\rho_2$ is the rest. Then
applying the gluing map $Gl_{\strata,\strata'}$ we have
$$
(u,\mkj,\rho_1,\rho_2) \to
(Gl_{\strata,\strata'}(u,\mkj,\rho_1),\rho_2) \in \mc
L_{\strata'}.
$$
We denote this gluing map on the bundle level by $BGl$. It is
clear that
\begin{lemma}\label{lemma_9.4.1}
$Gl_{\strata,\strata'}^\ast (\mc L_{\strata'})= \mc L_{\strata}$.
\end{lemma}
Suppose $\strata''$ is any stratum that is bigger than $\strata'$.
Set $W= Gl_{\strata,\strata'}(\mc
L_{\strata,\strata',\epsilon}^0)$. Had we proved that
$Gl_{\strata,\strata'}$ is a homeomorphic,
$Gl_{\strata,\strata''}$ would induce a gluing map
$$
Gl'_{\strata',\strata''}: \mc
L_{\strata',\strata'',\epsilon}^0|_W\to \M_{\strata''}(X,A)
$$
given by
$$
Gl'_{\strata',\strata''} =Gl_{\strata,\strata''}\circ
Gl_{\strata,\strata'}\inv.
$$
The homeomorphism (in fact, diffeomorphism) of
$Gl_{\strata,\strata'}$ will be proved in \S\ref{section_12}.


\section{Estimates}\label{section_10}

\subsection{Estimates for pre-gluing maps}\label{section_10.1}

We first prove lemma \ref{lemma_9.1.1}.

\v\n {Proof of lemma \ref{lemma_9.1.1}: }Let $\Sigma_1'= \Sigma_1-
D_{y_{o1}}(r^{1/2})$. We have
\begin{eqnarray*}
\|\bar{\partial}_{J,\mkj_{o\rho}}\phi\|_{L^p (\Sigma_1')} &\leq&
\|\bar{\partial}_{J,\mkj_o} u\|_{L^p(\Sigma_1)}+
C\left(\int_{D_{y_{o1}}(2r^{1/4} )}|\nabla
\beta(\frac{x}{r^{1/4}}) (u-p)|^p\right)^{1/p} \\
&&+ C\left(\int_{D_{y_{o1}}(2r^{1/4} )}|\nabla
J\cdot(u-p)|^p\right)^{1/p}.
\end{eqnarray*}
Note that in $D_{y_{o1}}(2r^{1/4})$
$$
|\nabla \beta(\frac{x}{r^{1/4}}) (u-p)|\leq C|u|_{C^1};
$$
and
$$
|\nabla J\cdot(u-p)|\leq C|J|_{C^1}r^{1/4}.
$$
So  on $\Sigma_1'\subset \Sigma_\rho$
$$
\|\bar{\partial}_{J, \mkj_{o\rho}}\phi\|_{L^p(\Sigma_1')} \leq
\|\bar{\partial}_{J,\mkj_o} u\|_{L^p} +Cr^{\frac{1}{2p}},
$$
One can compute the other side on $\Sigma_2$ similarly, so the
lemma follows.
 q.e.d.

\v We are also interested in the derivative of pre-gluing maps.
Let
$$
u_t=(u_{1t}, u_{2t}), t\in [0,\delta)
$$
be a path in $\chi^{1,p}_{\mkj_o}$ with
$$
\dot{u_t}|_{t=0}=\left.\frac{d}{dt}\right|_{t=0}u_t= \zeta:=
(\zeta_1,\zeta_2),
$$
Let $\phi_t= \pgl(u_t,\mkj_o, \rho)$ we study $
\dot{\phi}_t|_{t=0}. $ Similar to the computations for previous
lemma, we have
\begin{lemma}\label{lemma_10.1.1}
 Let $\zeta=(\zeta_1,\zeta_2)$ be as above,
Then
\begin{eqnarray*}
\left\|\left.\frac{d}{dt}\right|_{t=0}\phi_t\right\|_{L^{1,p}}
&\leq& C \|\zeta\|_{L^{1,p}};\\
\left\|\left.\frac{d}{dt}\right|_{t=0}\bar{\partial}_{J,\mkj_{o\rho}}\phi_t\right\|_{L^p}
&\leq&
\left\|\left.\frac{d}{dt}\right|_{t=0}\bar{\partial}_{J,\mkj_{o}}u_t\right\|_{L^p}
+
 Cr^{1/2p}\|\zeta\|_{C^1},
\end{eqnarray*}
where $C$ is independent of $\rho$. In particular,
$$
\left\|\left.\frac{d}{dt}\right|_{t=0}\bar{\partial}_{J,\mkj_{o\rho}}\phi_t\right\|_{L^p}
\leq
 Cr^{1/2p}\|\zeta\|_{C^1},
$$
if $u_t$ is a holomorphic path.
\end{lemma}
We leave the proof to readers. Note that the last term in the
estimates is $\zeta$ with respect to $C^1$-norm rather than
$L^{1,p}$-norm. Also for the last statement, it is clear that it
holds as long as $\zeta\in \ker D_{i_o,u_0}$.

\subsection{Estimates for right inverses}\label{section_10.2}

{\bf Proof of lemma \ref{lemma_9.2.1}: }Suppose $\sigma_1$ is
constructed from  $\sigma= Q_{u,\mkj_o}(\Lambda\eta)$ as explained
in \S\ref{section_9.2}. It is supported in $\Sigma_1''$. We
compute
$$
I:=D_{\phi,\mkj_{o\rho}}(\sigma_1)=
D_{\phi,\mkj_{o\rho}}(\sigma_1)(\gamma_{1,\rho}(\sigma-
\sigma(y_{o1}))+ \sigma(y_{o1})).
$$
We find that
\begin{eqnarray*}
|I|&\leq &  |\lambda_{1,\rho}\eta| + |\nabla \gamma_{1,\rho} (\sigma-\sigma(y_{o1}))|\\
&& + |J(\phi)\nabla\gamma_{1,\rho}(\sigma-\sigma(y_{o1}))| +
|(J(\phi)-J(u))\gamma_{1,\rho}d\sigma|
\\
&& + |\langle \nabla J, \gamma_{1,\rho}\sigma \rangle d(\phi-u)| +
|\langle \nabla J, \sigma(y_{o1})-\gamma_{1,\rho}\sigma(y_{o1}) \rangle d\phi|\\
&=:& I_1+I_2+I_3+I_4+I_5+I_6.
\end{eqnarray*}
The difficult terms are $I_2$ and $I_3$. They behave similarly:
for example,
$$
\|I_2\|_{L^p}^p \leq \int_{D_{y_{o1}}(r^{1/4})}
\left(\frac{C}{|\log r|}\frac{1}{r^{1/4}}
(r^{1/4})^{1-2/p}|\sigma|_{C^\alpha} \right)^p = \frac{C}{|\log
r|^p}|\sigma|_{C^\alpha}^p
$$
where $\alpha= 1-2/p$. So
$$
\|I_2\|_{L^p}\leq \frac{C}{|\log r|}\|\eta\|_{L^p}.
$$
estimates for $I_4$ to $I_6$ are trivial, so the claim is true.
q.e.d.

\v Let $u_t, \phi_t$ be as before.
\begin{prop}
Let $Q_{\phi_t,\mkj_{o\rho}}$ be the right inverse to $D_{
\phi_t,\mkj_{o\rho}}$ constructed as before. Then
\begin{eqnarray*}
\|Q_{\phi_t,\mkj_{o\rho}}\| &\leq & C;\\
\|\frac{\partial}{\partial \zeta} Q_{\phi_t,\mkj_{o\rho}}\| &\leq&
C\|\zeta\|_{L^{1,p}},
\end{eqnarray*}
where $C$ are constants depending only on $u$.
\end{prop}
{\bf Proof. } All statements are standard except the last
estimate. We explain this.
$$
\frac{\partial}{\partial \zeta} Q_{\phi_t,\mkj_{o\rho}}
=(\frac{\partial}{\partial \zeta}R)(D_{\phi_t,\mkj_{o\rho}}R)\inv
+ R\frac{\partial}{\partial \zeta}(D_{\phi_t,\mkj_o}R)\inv.
$$
For the first term it is sufficient to estimate
$$
\frac{\partial}{\partial \zeta}R =\Gamma (\frac{\partial}{\partial
\zeta}Q_{u_t,\mkj_o})\Lambda.
$$
It is standard to have
$$
\|\frac{\partial}{\partial \zeta}Q_{u_t,\mkj_o}\| \leq
C\|\zeta\|_{L^{1,p}}.
$$
and therefore
$$
\|\frac{\partial}{\partial \zeta}R\|\leq C\|\zeta\|_{L^{1,p}}.
$$
For the second term, we use the identity
$$
\frac{\partial}{\partial \zeta}(D_{\phi_t,\mkj_{o\rho}}R)\inv=
 -(D_{\phi_t,\mkj_o}R)\inv
\frac{\partial}{\partial \zeta}(D_{\phi_t,\mkj_{o\rho}}R)
(D_{\phi_t,\mkj_{o\rho}}R)\inv.
$$
Then the rest of estimates is standard. q.e.d.

\subsection{Estimates of $f(u,\mkj,\rho)$}\label{section_10.3}

As a consequence of theorem \ref{theorem_7.4.1}, we have
\begin{theorem}\label{theorem_10.3.1}
Let $\zeta\in \ker D_{u,\mkj}$. Then
$$
\|\frac{\partial}{\partial \zeta}f(u,\mkj_o,\rho)\|_{L^p}\leq
Cr^{1/p}\|\zeta\|_{L^{1,p}},
$$
where $C$ depends only on $u$.
\end{theorem}


\section{$C^0$-compatibility of gluing maps}\label{section_11}

\subsection{Admissible  gluing maps}\label{section_11.1}

As we have seen,  gluing maps consist of two parts: pre-gluing and
perturbation, i.e, map $\pgl$ and $\pert$ described in
\S\ref{section_9}. Hence, they depends on cut-off functions and
right inverses used in the constructions. Since cut-off functions
only depend on the coordinates of horocycles, we may assume that
cut-off functions are fixed. This kills the ambiguities caused by
cut-off functions.

On the other hand, there are more general gluing maps realized by
the following data (again, we only explain for the 1-nodal stratum
case): let $\Lambda=(V,\mc Q)$ be a pair satisfying assumption
\ref{assumption_7.3.1}. Suppose that it generates  a data of
coordinate chart $(V,\Phi,F)$ of a  proper open  subset $U$ of
$\M_\strata(X,A)$. We may define a gluing map $Gl_\Lambda$ based
on these data:
\begin{itemize}
\item for $\mk p=(u,\mkj_o,\rho)\in \mc L_\strata|_{U}$ we define
$$
\pgl_\Lambda(\mk p)= \pgl(F\inv(u), \mkj_o,\rho),
$$
set $\phi=\pgl_\Lambda(\mk p)$; \item construct right inverse for
$Q_{\phi,\mkj_{o\rho}}$ by using $Q_{F\inv(u),\mkj_o}$; \item
construct $\pert_\Lambda$ by using $\phi$ and
$Q_{\phi,\mkj_{o\rho}}.$
\end{itemize}
More explicit, $Gl_\Lambda$ is the composition
$$
\mc L_{\strata}|_{U}\to F^\ast \mc L_{\strata}|_{V}
\xrightarrow{Gl} \M_{g,m}(X,A).
$$
We call a gluing map $Gl_\Lambda$ constructed as above is an {\em
admissible gluing map}. Clearly, the original gluing maps are
admissible.
\begin{defn}\label{defn_11.1.1}
$Gl_\Lambda$ is called type-1 if $V\subset \M_\strata(X,A)$,
otherwise, it is called type-2.
\end{defn}

\subsection{$C^0$-compatibility}\label{section_11.2}

Suppose that we have two different gluing maps $Gl_\Lambda$,
$\Gamma=(V,\mc Q)$ and $Gl$. Later, we will prove that both of
them are compatible with the smooth structure of top stratum. How
they compatible with each other? Note that all gluing maps are
identity when $\rho=0$. We want to understand how much difference
between $Gl(u,\mkj_o,\rho)$ and $Gl_\Gamma(u,\mkj_o,\rho)$ when
$\rho\to 0$. The expected result should be
\begin{theorem}
$\lim_{\rho\to 0} Gl(u,,\mkj_o,\rho)=\lim_{\rho\to 0}
 Gl_\Gamma(u,\mkj_o,\rho)$
\end{theorem}
{\bf Proof. } We show that
$$
\|Gl(u,\mkj_o,\rho)- Gl_\Gamma(u,\mkj_o,\rho)\|\leq C(\rho)
$$
where $C(\rho)\to 0$ for $\rho\to 0$.

For simplicity, we introduce notations. Let
\begin{eqnarray*}
&&u'_o= F\inv(u_o),\\
&& \eta_o=f(u'_o), \\
&&\sigma_o= u_o-u'_o= Q_{u'_o}\eta_o.
\end{eqnarray*}
Let
$$
\phi_o'=\pgl(u'_o,\mkj_o,\rho), \phi_o=\pgl(u_o,\mkj_o,\rho).
$$
Let $\Lambda$ and $R$ be those terms in \S \ref{section_9.2}. We
compare
$$
\phi_o'':= \phi_o'+ Q_{\phi_o',\mkj_{o\rho}}(\Lambda \eta_o)
$$
with $\phi_o$. We claim that
\begin{eqnarray}
&&\|\phi_o''-\phi_o\|\leq C_1(\rho);\label{eqn_11.2.1}\\
&& \|\bar\partial_{J,\mkj_\rho} \phi_o''\|\leq C_2(\rho),
\label{eqn_11.2.2}
\end{eqnarray}
where $C_j(\rho)\to 0, j=1,2,$ when $\rho\to 0$. These two
equations imply this theorem by theorem \ref{theorem_9.3.1}.
 The proof of these two equations is given below.
q.e.d.

\v\n
\begin{prop}\label{prop_11.2.1}
Equation \ref{eqn_11.2.1} is true.
\end{prop}
{\bf Proof. } Step 1,
\begin{equation}
\|Q_{\phi_o'}(\Lambda\eta_o) -R(\Lambda\eta_o)
 \|\leq \frac{C}{|\log \rho|}
\|\eta_o\|.
\end{equation}
This follows directly by the definition of $Q_{\phi_o'}$.

It remains to compare $ R(\Lambda\eta_o)$ with $\beta\cdot
\sigma_o$. By definition
$$
R(\Lambda\eta_o)= \Gamma Q_{u_o'}(\lambda \eta_o).
$$
Note that $\sigma_o= Q_{u_o'}\eta_o$. We can easily verify that
$$
\|\Gamma Q_{u_o'}(\lambda \eta_o)-\beta\cdot\sigma_o\|\leq
Cr^{1/2p}.
$$
Combine these, we get equation \ref{eqn_11.2.1}. q.e.d.

\v\n
\begin{prop}
Equation \ref{eqn_11.2.2} is true.
\end{prop}
{\bf Proof. }We have that
$$
\|\bar\partial \phi_o\|\leq Cr^{1/2p}.
$$
and
$$
\|\bar\partial (\phi_o''-\phi_o)\|_{L^p} \leq
C\|\nabla(\phi_o''-\phi_o)\|_{L^p} \leq C(\rho).
$$
So \ref{eqn_11.2.2} follows. q.e.d.

\v\n As a corollary, we have that
\begin{corollary}
$Gl\inv Gl_\Gamma$ and its inverse are continuous.
\end{corollary}


\section{Coordinate charts from gluing maps}\label{section_12}

We explain that how the differential structure on $\M_{g,m}(X,A)$
induced by gluing maps fits with the one given in
\S\ref{section_7.2}.

We discuss these case by case:
\begin{itemize}
\item Case I: $2g_1+ m_1\geq 3$ and $2g_2+ m-m_1\geq 3$; \item
Case II: $2g_1+ m_1\geq 3$ and $2g_2+ m-m_1 < 3$; \item Case III:
$2g_1+ m_1 < 3$ and $2g_2+ m-m_1 < 3$;
\end{itemize}

\subsection{Case I}\label{section_12.1}
We study the gluing maps near $(u_o,\mkj_o)$. By assumption
$\mkj_{ov},v=1,2$ are stable, so is $\mkj_o$. For simplicity, we
will ignore finite groups $\Gamma_{\mkj_o}, \Gamma_{\mkj_{ov}}$
and etc. unless it is stressed.

Let $M_S$ be the stratum containing $\mkj_o$. For simplicity, we
assume that $M_S$ and $\M_\strata(X,A)$ are compact. It is known
that
$$
\gs: L_{S,\epsilon}^0  \to M_{g,m}
$$
is a local diffeomorphism.

Let $O$ be any neighborhood of $\mkj_o$ in $M_S$. Let
$$
\mc O= \mk f\inv(O),
$$
where $\mk f$ is the forgetting-map map.

Set
$$
U=\pgl(\mc L_\strata|_{\mc O})\subset \chi^{1,p}_{g,m}(X,A).
$$
For each $\phi=\pgl(u,\mkj,\rho)\in U$ we have right inverse
$Q_{\phi,\mkj_\rho}$. If we fix $(\mkj,\rho)$, set
\begin{eqnarray*}
U_{\mkj,\rho} &=& \{\pgl(*,\mkj,\rho)\}, \\
\mc Q_{\mkj,\rho} &=& \{Q_{\phi,\mkj_\rho}|\phi\in
U_{\mkj,\rho}\}.
\end{eqnarray*}
By estimates in \S\ref{section_9.1} and \S \ref{section_9.2}, we
have
\begin{theorem}\label{theorem_12.1.1}
$(U_{\mkj,\rho},\mc Q_{\mkj,\rho})$ is a  pair satisfying
assumption \ref{assumption_7.3.1}. Hence the gluing map generates
a coordinate chart of $\M_{\mkj_\rho}(X,A)$. Since
$\{(\mkj,\rho)\}$ may be treated as parameters, $(U,\mc Q)$
generates a coordinate chart of $\M_{g,m}(X,A)$ which is given by
gluing maps.
\end{theorem}
In the other word,
 $Gl_{\mkj_o,\rho}$ is diffeomorphic
automatically.

\subsection{Gluing maps: case II}\label{section_12.2}
We now discuss the gluing for case II. That is:
$\mkj_{o1}$ is stable and $\mkj_{o2}$ is
unstable. In particular, we note that $\Sigma_2=S^2$.

We will further divide case II into four subcases:
\begin{enumerate}
\item[IIa.] $m=m_1$ and $(\Sigma_1, i_1,x_1,\ldots ,x_{m})$ is
stable;

\item[IIb.] $m_1=m-1$ and $(\Sigma_1, i_1,x_1,\ldots ,x_{m-1})$ is
stable;

\item[IIc.] $m=m_1$ and $(\Sigma_1, i_1,x_1,\ldots ,x_{m})$ is
unstable;

\item[IId.] $m_1=m-1$ and $(\Sigma_1, i_1,x_1,\ldots ,x_{m-1})$ is
unstable;
\end{enumerate}

We start with  case IIa  which is one of the most complicated
cases. Before we proceed, let us remark what is new comparing with
case 1. The problem is that $\gs$ is no longer local
diffeomorphic. Hence, we are not able to treat $L_S$ as
parameters.

\v \n {\em Case IIa}. We specify the notations for this case.
$\mkj_o$ consists of
$$
\mkj_{o1}=(\Sigma_1,j_{o1},x_{o1},\ldots,x_{om}, y_{o1})
$$
and
$$
\mkj_{o2}=(S^2, y_{o2}=\infty).
$$
We describe $M_S$. For simplicity, we assume $m=0$ and
$(\Sigma_1,j)$ is stable. For
$$
\mkj_{1}=(\Sigma_1,j_{1}, y_{1})
$$
set
$$
\mkj_{1}'=(\Sigma_1,j_{1}).
$$
Then
$$
M_S\cong M_{g_1}\times \Sigma_1 \times \{\mkj_{o2}\}
$$
where the isomorphism is given by
$$
(\mkj_1,\mkj_{o2})\leftrightarrow (\mkj_1', y_1,\mkj_{o2}).
$$
By the construction of
$$
\mr{gs}: M_S\rtimes \cplane^\ast_{\epsilon} \to \M_{g}
$$
we know it is an fibration with fiber
$$
\Sigma_1\rtimes \cplane^\ast.
$$
Geometrically, this says: with fixed surface
$$\mkj_1'=(\Sigma_1,j_1),$$
for any $ y\in \Sigma_1 $ and  $0\not=\rho\in
\cplane_\epsilon^\ast$,
$$
gs(\mkj_1,\rho)=\mkj_1'
$$

Let $u=(u_1,u_2)\in \M_{\strata}(X,A)$ be a map. We may assume
that $u_2$ is balanced.  Be precise, we define
\begin{eqnarray*}
\M_{\strata}^b(X,A) &=&\{(u_1,u_2)|u_1\in \M_{g,1}(X,A_1),\\
&&u_2\in\M^b_{0,1}(X,A_2), u_1(y)=u_2(y)\}.
\end{eqnarray*}
Then
$$
\M_{\strata}(X,A)=\frac{\M_{\strata}^b(X,A)}{S^1}.
$$

By this exposition, we know that: $y$ and $\rho$ can not be
treated as parameters, however $\mkj_1'$ can be. So we will fixed
$\mkj_1'$ in the rest of argument for this subcase. This is
equivalent to fixing $j_1$.

We summarize the notations again: $\mkj_o$ consists of
$$
\mkj_{o1}=(\Sigma_1, j_{o1}, y_{o1}) \mbox{ and } \mkj_{o2}=(S^2,
\infty)
$$
Set $\mkj_{o1}'=(\Sigma_1,j_{o1})$. Define
$$
M_{S,{\mkj_{o1}'}}=\{(\mkj_1:= (\Sigma_1,j_{o1},y_1), \mkj_{o2})\}
$$
Clearly, $M_{S,\mkj_{o1}'}\cong \Sigma_1$. Then
$$
\gs: M_{S,\mkj'_{o1}} \times \cplane^\ast_\epsilon\to\mkj_{o1}'.
$$
Correspondingly, for moduli spaces,  we have $\M^b_{\mkj}(X,A)$
and $\M_\mkj(X,A)$ for $\mkj\in M_{S,\mkj_{o1}'}$. Set
$$
\M_{\strata,\mkj_{o1}'}^b(X,A) =\coprod_{\mkj\in
M_{S,\mkj_{o1}'}}\M_{\mkj}^b(X,A)
$$
and
$$
\M_{\strata,\mkj'_{o1}}(X,A):
=\frac{\M_{\strata,\mkj'_{o1}}^b(X,A)}{S^1}.
$$

For any
$$
u\in \M_{\mkj}^b(X,A)\subset \M_{\strata,\mkj'_{o1}}^b(X,A),
$$ we assume that $D_{u,\mkj}$ is
surjective. Then we get a gluing map
$$
Gl: \M_{\strata,\mkj_{o1}'}^b(X,A)\times_{S^1}
\cplane^\ast_\epsilon\to \M_{\mkj_{o1}}(X,A).
$$
The map is well defined: since it is easy to see that the gluing
map defined on
$$
\M_{\mkj_{o1}'}^b(X,A)\times \cplane^\ast_\epsilon
$$ is
$S^1$-equivariant.  The balanced moduli spaces are
 necessary for the equivariance.
We move on to discuss the diffeomorphic issue.

Fix a map $u_o=(u_{o1},u_{o2})\in \M_{\strata,\mkj_{o}}^b(X,A)$.
Since $D_{u_o,\mkj_o}$ is surjective,
$\M_{\strata,\mkj_{o1},y_{o1}}^b(X,A)$ is a smooth manifold. Let
$N_o$ be a slice (with respect to the $S^1$-action) through $u_o$.
Let $V$ be a neighborhood of $y_{o1}\in \Sigma_1$. Then the
neighborhood $U_{u_o}$ of $u_o\in \M_{\strata,\mkj_{o1}'}(X,A)$
 can be identified with
$$
U_{u_o}\cong V\times N_o.
$$
Fix $\rho_o=r_o$ and its neighborhood $Glu(\rho_o)\in
\cplane^\ast_{\epsilon}$. Then
 the gluing map is
locally rewritten as
\begin{equation}\label{eqn_12.2.1}
Gl: V\times N\times Glu(\rho_o)\to \M_{\mkj_{o1}'}(X,A).
\end{equation}
We want to show that this is local diffeomorphic. The new point is
to compute differentiation with respect to new variables in
$V\times Glu(\rho_o)$. To treat them properly, we compare this map
with
 with another well-studied map,
 which has been shown to be diffeomorphic
 by case 1.

 By
adding two marked points $\{0,1\}$ to $S^2$, we get a stable curve
$$
\bar\mkj_{o2}=(S^2, 0,1,\infty).
$$
Let $\bar\mkj_{o}=(\mkj_{o1},\bar\mkj_{o2})$. Set
$$
\bar\mkj'=\gs(\bar\mkj_o,\rho_o)\in M_{g,2}.
$$
We regard $u_o$ as an element in $\M_{\bar\mkj_o}(X,A)$ in an
obvious way. Let $\bar N$ be a neighborhood of $u_o$ in this
moduli space. We have a gluing map
$$
Gl_{\bar \mkj_o}: \bar N\times \{\rho_o\}\to \M_{\bar \mkj'}(X,A)
$$
which is diffeomorphic according to case I. We rewrite the map as
$$
Gl_{\bar \mkj_o,\rho_o}: \bar N\to \M_{\bar \mkj'}(X,A).
$$
Since $y_{o1}$ and $\rho_o$ are fixed, $\bar\mkj'$
 can be identified with $\mkj_{o1}'$
by forgetting the two extra marking points. This induces an
isomorphism
$$
\M_{\bar\mkj'}(X,A)\leftrightarrow \M_{\mkj_{o1}'}(X,A)
$$
via forgetting-marked-point. So we have
$$
Gl_{\bar \mkj_o', \rho_o}: \bar N\to \M_{\mkj_{o1}'}(X,A).
$$

Next, we explain that there is a natural isomorphism
$$
B: V\times N\times Glu(\rho_o)\to \bar N.
$$
Had $Gl=Gl_{\bar \mkj_o', \rho_o}\circ B$, we would prove that the
former one is diffeomorphic. Though this is not case, they are
rather close. This is what we do next.

We know that $\bar N=\mk B_0N$. By $\mk B_0$ we mean a
neighborhood of identity $(0,1)\in \mk B$. So it is sufficient to
define a map $b: V\times Glu(\rho_o)\to \mk B_0$. This is given by
$$
b(y,\rho)= (r_o\inv(y-y_0),r_o\inv \rho).
$$
So
$$
B(y,u, \rho)= b(y,\rho)\cdot u.
$$
Set
$$
Gl_{\bar \mkj_o,\rho_o}'= Gl\circ B\inv.
$$
And think of $Gl_{\bar \mkj_o,\rho_o}$ and $Gl_{\bar
\mkj_o,\rho_o}'$ are both maps from $\bar N=\mk B_0N $ to
$\M_{\mkj_{o1}'}(M,A)$.

Let
$$
\phi_o= \pgl(u_o,\mkj_o,\rho_o), \hat u_o= Gl(u_o,\mkj_o,\rho_o).
$$
Set
\begin{equation}
W=W^{1,p}(\phi_o^\ast TM), L=L^p(\Lambda^{0,1}
_{\mkj_{o1}'}\phi_o^\ast TM).
\end{equation}
We know that $Gl_{\bar \mkj_o',\rho_o}$ induces the following data
of a coordinate chart:
\begin{enumerate}
\item $O_{\bar\mkj_o',\rho_o}= \pgl(\mk B_0N,\bar\mkj_o,\rho_o)$;
\item $\Phi_{\bar\mkj_o',\rho_o}: O_{\bar\mkj_o',\rho_o} \times
L\to W$ given by
\begin{equation}\label{eqn_12.2.2}
\Phi_{\bar\mkj_o',\rho_o}(\phi,\eta)=\phi + Q_{\phi,
\mkj_{o1}'}\eta;
\end{equation}
\item $\pert_{\bar\mkj'_o,\rho_o}: O\to L$  the map that yields
the gluing map $Gl_{\mkj_o',\rho_o}$.
\end{enumerate}

\begin{theorem}\label{theorem_12.2.1}
The following is the data of  a coordinate chart induced by
$Gl_{\bar \mkj_o',\rho_o}'$:
\begin{enumerate}
\item $O:= \pgl(V\times N\times Glu(\rho_o)\times \{\mkj_{o1}\})
=\pgl(B\inv(\mk B_0N),\mkj_{o1})$; \item $\Psi: O\times L\to W$
given by
\begin{equation}\label{eqn_12.2.3}
\Psi(\phi,\eta)= \psi+ Q_{\psi,\mkj_{o1}'}\eta;
\end{equation}
\item $\pert: O\to L$  the map that yields the gluing map $Gl$.
\end{enumerate}
\end{theorem}
{\bf Proof. } We only need to verify that $\Psi$ is an
isomorphism. Since $O\cong \mk B_0N$, it is equivalent to show
that
$$
\tilde\Psi= \Psi\circ(\pgl\circ B\inv, 1):
 \mk B_0N\times
L\to W
$$
is an isomorphism.

On the other hand, by \eqref{eqn_12.2.1} we know that
$$
\tilde\Phi=\Phi_{\mkj_o',\rho_o}\circ \pgl:\mk B_0N\times L\to W
$$
is an isomorphism. Now both maps $\tilde\Psi$ and $\tilde\Phi$
 have same domain and range. We
claim that when $r_o$ is small,
\begin{eqnarray}
&&\|\tilde\Psi-\tilde\Psi\|\leq \epsilon, \label{eqn_12.2.4} \\
&&\|d\tilde\Psi-d\tilde\Psi\|\leq \epsilon \label{eqn_12.2.5}
\end{eqnarray}
for some small $\epsilon$. Then that the isomorphism of
$\tilde\Psi$ implies that of $\tilde\Psi.$ The proof of
\eqref{eqn_12.2.4} and  \eqref{eqn_12.2.5} is rather
straightforward but tedious. The proof of them is explained below.
q.e.d.

\begin{defn}\label{defn_12.2.1}
Let $Gl_1$ and $Gl_2$ be two gluing maps defined on same (local)
domain. Let $\Psi_1$ and $\Psi_2$ are corresponding maps defined
in the form as \eqref{eqn_12.2.3}.
 We say
$$
Gl_1\approx Gl_2
$$
if $\Psi_1-\Psi_2$ satisfies \eqref{eqn_12.2.4} and
\eqref{eqn_12.2.5}.
\end{defn}
\v\n To compare $\tilde\Psi$ and $\tilde\Phi$, we should go
through the process of gluing and compare them in each step.

\v\n {\em Pre-gluing maps}. We first compare the pre-gluing maps
for two different gluing processes.

\v\n Suppose $(t, z)\in \mk B$ is given. Let $(y,\rho)=
b\inv(t,z)$. Let $u=(u_1,u_2)\in N$ and
$$
\tilde u=(u_1, (t,z)\cdot u_2)=:(u_1,\tilde u_2).
$$
Set
$$
\mkj_y=((\Sigma_1,j_{o1},y), \mkj_{o2}).
$$
The pre-gluing for $Gl'_{\bar \mkj_{o},\rho_o}$ is
$\pgl(u,\mkj_y,\rho)$ and that for $Gl_{\bar\mkj_o, \rho_o}$ is
$\pgl(\tilde u, \bar\mkj_o,\rho_o)$. We denote them by
$$
\pgl'_{\bar\mkj_o}, \pgl_{\bar\mkj_o}: \mk B_0N\to W
$$
respectively. We have
\begin{prop}\label{prop_12.2.1}
Let $(t,z)=b(y,\rho),u\in N$, i.e,
$$
\rho=r_oz,y=r_ot.
$$
Then
$$
\|\pgl'_{\bar\mkj_o}((t,z)\cdot u), \pgl_{\bar\mkj_o}((t,z)\cdot
u)\|_{L^{1,p}} \leq C|t|\sqrt{r_o}.
$$
Here $C$ is a constant independent of $r_o$.
\end{prop}
{\bf Proof. } Set
$$
\phi'=\pgl'_{\bar\mkj_o}((t,z)\cdot u);
\phi=\pgl_{\bar\mkj_o}((t,z)\cdot u).
$$
By the construction of pre-gluing, $\phi'$ and $\phi$ are maps on
$\Sigma_{y,\rho}$ and $\Sigma_{y_o,\rho_o}$.
 We should identify them properly: in fact,
both of them  are identified with $\Sigma_1$ in a canonical way
and so they are identified. In particular, we explain how two
sphere components identified. We name the spheres $S^2_y$ and
$S^2_{y_o}$. Let
$$
\cplane_y=S^2_y-\{\infty\}; \cplane_{y_o}=S^2_{y_o}-\{\infty\}.
$$
We write down  the identification map
\begin{eqnarray*}
&&w: \cplane_y\to \cplane_{y_o};\\
&&w(z)= \left(r_0\inv(y+\rho z\inv) \right)\inv,
\end{eqnarray*}
The inverse of $w$ is
$$
w\inv(z)=\left(\rho\inv(r_0x\inv -y ) \right)\inv.
$$
Explicitly, we write down $\phi$ and $\phi'$ on
$\Sigma_{y_o,\rho_o}$. We separate $\Sigma_{y_o,\rho_o}$ into
three pieces:
\begin{itemize}
\item $P_1:=\Sigma_1'= \Sigma_1-D_{y_o}(2\sqrt {r_o})$;

\item $P_2:=\Sigma_2'=\Sigma_2-D_{y_o}(2\sqrt {r_o})$;

\item $P_3= D_{y_o}(2\sqrt{r_o})- D_{y_o}(\sqrt {r_o}/2)\subset
\Sigma_1$.
\end{itemize}
On $P_3$,
$$
\phi= \phi'= u(y_o).
$$
On $P_1$
\begin{eqnarray}
&&\phi(z) =u(y_o)+
\beta(\frac{z}{\sqrt {r_o}})(u_1(z)-u(y_o)),\\
&&\phi'(z)= u(y_0)+ \beta(\frac{z-y}{\sqrt r})(u_1(z)-u(y_o)).
\end{eqnarray}
On $P_2$
\begin{eqnarray}
&&\phi(z) =u(y_o)+
\beta(\frac{z}{\sqrt {r_o}})(\tilde u_2(z)-u(y_o)),\\
&&\phi'(z)= u(y_o)+ \beta(\frac{w\inv(z)}{\sqrt r})(\tilde
u_2(z)-u(y_o)).
\end{eqnarray}
Clearly, to prove the proposition, the computation of cut-off
functions is involved. We need the results from appendix
\ref{appendix-cut-off}.

We explain the computation on $P_1$.
$$
\phi(z)-\phi'(z) =(\beta(\frac{z}{\sqrt
{r_o}})-\beta(\frac{z-y}{\sqrt r}))(u_1(z)-u(y_o)).
$$
This is supported  in $D_{y_o}(3\sqrt
{r_o}N)-D_{y_o}(\sqrt{r_o}N/2)$ And we have estimates in this
area:
$$
|(\beta(\frac{z}{\sqrt {r_o}})-\beta(\frac{z-y}{\sqrt
r}))(u_1(z)-u(y_o))| \leq
C|\frac{y}{\sqrt{r_o}}|\sqrt{r_o}^{1-2/p}
$$
and
$$
\left|\nabla\left((\beta(\frac{z}{\sqrt
{r_o}})-\beta(\frac{z-y}{\sqrt
r}))\left(u_1(z)-u(y_o\right)\right)\right| \leq
C|\frac{y}{r_o}|\sqrt{r_o}^{1-2/p} +C|\frac{y}{\sqrt{r_o}}|.
$$
Then their $L^p$-norms  are bounded by
$$
C|y| +C|y|r_o^{-1/2}+ C|y|r_o^{1-p/2}.
$$
Plug in $y=r_ot$, we have
$$
\|\phi-\phi'\|_{L^{1,p}(P_1)}\leq C|t|\sqrt{r_o}.
$$
The computation on $P_2$ is same. Then the claim of proposition
follows. q.e.d.

\begin{remark}\label{rmk_12.2.1}
The key to the whole process is that we use coordinate $(t,z)$
rather than $(y,\rho)$: we note that the computation of cut-off
functions with respect to $(y,\rho)$ does not preform friendly,
while there is no problem when it is with respect to $(t,z)$. This
is  due to the factor $r_0$. On the other hand, we know that it is
$(t,z)\in \mk B$ that is essential inspired by the map
$Gl_{\bar\mkj_o,\rho_o}$. So it is not surprise that the
computation behaves well. We will skip  the computations of the
rest of these type results. It is just a matter of recycling the
above computations and those in appendix.
\end{remark}
Similar computations  imply
\begin{prop}\label{prop_12.2.2}
Given $u\in N$ and  a path $(t(s),z(s))\in \mk B, s\in [0,1)$ with
$(t,z)=(t(0),z(0))$ and
$$
(v_1,v_2)=\left.\frac{\partial}{\partial
s}\right|_{s=0}(t(s),z(s)),
$$
 Then
$$
\|\left.\frac{\partial}{\partial
s}\right|_{s=0}\pgl_{\bar\mkj_o}((t(s),z(s))\cdot u)-
\left.\frac{\partial}{\partial
s}\right|_{s=0}\pgl_{\bar\mkj_o}'((t(s),z(s))\cdot u) \|_{L^{1,p}}
\leq C\sqrt{r_o}|(v_1,v_2)|.
$$
In particular, this implies that at $(t,s)\cdot u$
$$
\|d( \pgl'_{\bar\mkj_o}- d \pgl_{\bar\mkj_o}\|\leq C\sqrt{r_o}.
$$
when $(t,s)$ is bounded.
\end{prop}


\v\n{\em Right inverses.} Recall that in the construction of right
inverse $Q$, we first define $R=\Gamma Q_{u}\Lambda$ and then set
$Q_\phi= R(DR)\inv$.
 Here $\Gamma$ and $\Lambda$ involves cut-off
functions. Hence we should deal with the derivatives of cut-off
functions as well.

For gluing maps $Gl_{\bar\mkj_o.\rho_o}$ and $Gl_{\bar\mkj_o,
\rho_o}'$ we have two families of right inverses $\mc Q$ and $\mc
Q'$:
\begin{eqnarray*}
&&
\mc Q=\{Q_{\pgl_{\bar\mkj_o}(x),\mkj_{o1}}|x\in \mk B_0N\};\\
&& \mc Q'=\{Q_{\pgl'_{\bar\mkj_o}(x),\mkj_{o1}}|x\in \mk B_0N\}.
\end{eqnarray*}
We may treat them as maps
$$
\mc Q,\mc Q': \mk B_0N\times L\to W.
$$
Then
\begin{prop}\label{prop_12.2.3}
Let $\mc Q$ and $\mc Q'$ be as above.
\begin{eqnarray*}
 && \|Q-Q'\|\leq C\sqrt{r_0};\\
&&\|d Q-dQ'\|\leq C\sqrt{r_0}.
\end{eqnarray*}
\end{prop}
Combine these results, we prove theorem \eqref{eqn_12.2.4} and
\eqref{eqn_12.2.5}.

\begin{remark}\label{rmk_12.2.2}
We explain the idea that guides us in the above proof. Let
$$
(u_o,\rho_o)\in \tilde{M}_\strata(X,A)\times \cplane^\ast_\epsilon
$$
and $\tilde U_{u_o}\times Glu(\rho_o)$ be a neighborhood of this
point.
 We may be
expecting a gluing map
$$
\widetilde{Gl}: \frac{\tilde U_{u_o}\times Glu(\rho_o)}{\mk B} \to
\M_{g,m}(M,A).
$$
We may construct a gluing map defined on a proper chosen slice.
This is essentially what $Gl$ does. Another reasonable approach
would be $\widetilde{Gl}\Rightarrow \widetilde{Gl}_1\Rightarrow
Gl'$. Here
$$
\widetilde{Gl}_1: \frac{\tilde U_{u_o}}{\mk t} \times
\{\rho_o\}\to \M_{g,m}(M,A).
$$
We use $Glu(r_o)\cong \mk m$.

Set
$$
\mkj_{o,y}=((\Sigma_1,j_{o1},y),\mkj_o2); \mkj_o=\mkj_{o,y_o}.
$$
Let
$$
\wtm_V(X,A)=\coprod_{y\in V}\M_{\mkj_{o,y}}(X,A).
$$
Suppose that
$$
\wtm_{\mkj_{o,y}}(X,A)\cong \wtm_{\mkj_o}(X,A)
$$
and
$$
\wtm_{V}(X,A)\cong \wtm_{\mkj_o}(X,A)\times V.
$$
Set
$$
\tilde U_{u_o,y_o}=\tilde U_{u_o}\cap \wtm_{\mkj_{o}}(X,A).
$$
Then using a natural identification of $V$ with $\mk t$ we reduce
$\widetilde{Gl}_1$ to
$$
Gl': \tilde U_{u_o,y_o}(X,A)\times \{\rho_o\}\to \M_{g,m}(M,A).
$$
Elements in $\tilde U_{u_o,y_o}$ are treated pre-stable by adding
two marked points on $\mkj_{o2}$. Hence $Gl'$ is exactly
$Gl_{\bar\mkj_{o},\rho_o}$. So it is not surprise to have a
natural comparison between $Gl$ and $Gl'$. Although the
computation is tedious, it is quite straightforward.

From this, we also see that in this local comparison  $Glu(r_0)$
can always compare with $\mk m$. So we will always cancel
$Glu(r_0)$ with $\mk m$ when the similar issue occurs.
\end{remark}

\v\n{\em Case 2b and 2d.} These two cases are simpler. The group
$\mk B$ in case 2a is replaced by $\mk m< \mk B$. We skip them.

\v\n{\em Case 2c.} This is a relatively new case. The point is
that the resultant curves after gluing are pre-unstable. The
treatment of this case is same as case 3. We discuss case 3
directly.

\subsection{Gluing maps: case 3}\label{section_12.3}

Both
$$
(\Sigma_1,i_1,x_1,\ldots, x_{m_1},y) \;\;\;\mbox{and} \;\;\;
(\Sigma_2, i_2,x_{m_1+1},\ldots, x_m,y)
$$ are  unstable,
$\Sigma_j=S^2,j=1,2$ and $m_1\leq 1, m-m_1\leq 1$. We take the
most complicated case: $m=0$. To tell the difference between two
components, we mark spheres by $S^2_j$. Namely
\begin{equation}\label{eqn_12.3.1}
\Sigma= (S^2_1,y= \infty_1)\cup (S^2_2,y=\infty_2).
\end{equation}
Then
$$
\aut(\Sigma) =\mk B_1\times \mk B_2.
$$
We put the subscripts to tell the difference. We define a normal
subgroup of $\aut(\Sigma)$
\begin{eqnarray*}
\aut_y(\Sigma)&=&
\{(\psi_1,\psi_2)|\psi_i\in \aut(\Sigma_i),\\
&& d\psi_1(y_1)\otimes d\psi_2(y_2)|_{T_{y_1}(\Sigma_1) \otimes
T_{y_2}(\Sigma_2)}=1\}.
\end{eqnarray*}
Set
$$
\Delta^\ast (\mk m)=\{(m_1,m_2)\in \mk m_1\times \mk
m_2|m_1m_2=1\}.
$$
By direct computation, we have
\begin{lemma}\label{lemma_12.3.1}
$\aut_y(\Sigma)=\mk t_1\times \mk t_2\ltimes \Delta^\ast (\mk m)$.
\end{lemma}
Then
$$
\frac{\aut(\Sigma)}{\aut_y(\sigma)}\cong \mk m\cong \mk m_2.
$$
Now we consider the gluing. For each component, we use balanced
curves, i.e,
$$
\M_{0,1}(X,A_j)= \frac{M^b_{0,1}(M,A_j)}{S^1}.
$$
Therefore
$$
\M_{\strata}(X,A)= \frac{M_{\strata}^b(X,A)} {S^1\times S^1},
$$
where $A=A_1+A_2$. For simplicity, we assume that the stratum is
compact. The gluing is
\begin{equation}\label{eqn_12.3.2}
Gl:{M_{\strata}^b(X,A)}\times_{S^1\times S^1} \times
\cplane^*_\epsilon\to \wtm_{0,0}(X,A).
\end{equation}
or in a more precise  form, the right hand side is  treated as a
subset
$$
{M_{\strata}^b(X,A)}\times_{S^1\times S^1} \times
\cplane^*_\epsilon \subset
M^b_{0,0}(M,A_1)\times_{S^1}\cplane\times_{ S^1} M^b_{0,0}(M,A_2).
$$
Note that
$$
\M_{0,0}(X,A)=\frac{\wtm_{0,0}(X,A)}{\aut(S^2)}.
$$
In order to show that $Gl$ defined in \eqref{eqn_12.3.2} induces a
local-diffeomorphic gluing map
$$
{M_{\strata}^b(X,A)}\times_{S^1\times S^1} \times
\cplane^*_\epsilon\to \M_{0,0}(X,A).
$$
we should conclude that
\begin{theorem}\label{theorem_12.3.1}
The image of $Gl$ represents a slice in $\wtm_{0,0}(X,A)$ with
respect to the action of $\aut(S^2)$.
\end{theorem}
We now explain the idea following the guide line given in remark
\ref{rmk_12.2.2} to speculate the proof. The key is proposition
\ref{prop_12.3.1}.

Locally, an expecting map is
\begin{equation}\label{eqn_12.3.3}
\widetilde{Gl}: \frac{U}{\mk B_1\times \mk B_2} \times Glu(r_o)
\to \frac{\wtm_{0,0}(X,A)}{\aut(S^2).}
\end{equation}
Here $U\subset \wtm_\strata(X,A)$ is a small open subset in the
stratum that is $\mk B_1\times \mk B_2$ invariant. Again, $U$
should be thought as a subset  of
$$
\wtm_{0,1}(X,A_1)\times \wtm_{0,1}(X,A_2)
$$
and
$$
\frac{U}{\mk B_1\times \mk B_2} \subset
\frac{\wtm_{0,1}(X,A_1)}{\mk B_1}\times
\frac{\wtm_{0,1}(X,A_2)}{\mk B_2}.
$$
The left hand side of \eqref{eqn_12.3.3} can be written as
$$
\left(\left.\frac{U}{\aut_y(\Sigma)}\right/\mk m_1\right) \times
Glu(r_o).
$$
As before, locally $\mk m_1$ is cancelled by $Glu(\rho_0)$, and we
have
$$
\widetilde{Gl}_1: \frac{U}{\aut_y(\Sigma)}\times \{\rho_0\} \to
\frac{\tilde{\M}_{0,0}(M,A)}{\aut(S^2)}.
$$
Next we need an important fact for this kind of gluing.
\begin{prop}\label{prop_12.3.1}
For any $\rho\leq r_o$, there exists
 a neighborhood $V$ of $id$ in $\aut_y(\Sigma)$, a neighborhood
 $V'$ of $id$ in $\aut(S^2)$ and
a diffeomorphism map
$$
gl: V\to V'.
$$
Here $gl$ is construct via gluing process.
\end{prop}
We skip the proof.

Using this fact:
$$
\aut_y(\Sigma)\cong \aut(S^2)
$$
locally, we would show that the image of $Gl$ is a slice.

\v\n {\bf Sketch the proof of theorem \ref{theorem_12.3.1}: }
First we introduce a slice of
$$
\frac{\wtm_\strata(X,A)}{\aut_y(\Sigma)}.
$$
We say an element $(u_1,u_2)\in \wtm_\strata(X,A)$ is balanced
with respect to $\aut_y(\Sigma)$ if $u_1$ is balanced and $u_2$ is
centered. Let $\wtm_{\strata}^b(X,A)$ denote the set of such
elements. It is not hard to see that
$$
\frac{\wtm_\strata(X,A)}{\aut_y(\Sigma)}=
\frac{\wtm_\strata^b(X,A)}{S^1_1\times S^1_2},
$$
and
$$
\M_\strata^b(X,A)=\frac{\wtm_\strata^b(X,A)}{\mk m_2}.
$$
Set
$$
Gl_1: \wtm_\strata^b(X,A)\times \{\rho_0\} \to \wtm_0(X,A)
$$
to be  a gluing map. Then locally
$$
Gl\approx Gl_1
$$
in the sense of definition \ref{defn_12.2.1}. The problem is now
translated to show that the image of $Gl_1$ is a slice in
$\wtm_0(X,A)$.

Take a slice $N$ in $\wtm_\strata^b(X,A)$. A neighborhood of $N$
in $\wtm_{\strata}(X,A)$ is $V\cdot N$. Define a map
$$
Gl_2:V\cdot N\times \{\rho_0\}\to \wtm_0(X,A)
$$
by
$$
Gl_2(v,n)= gl(v)\cdot Gl(n).
$$
We compare it with the original gluing map, extending $Gl_1$,
$$
Gl:V\cdot N\times \{\rho_0\}\to \wtm_0(X,A).
$$
By using the property of $gl$, it is straightforward to show that
$$
Gl\approx Gl_2.
$$
Since $Gl$ is local diffeomorphic, we conclude that
$$
Gl_2(1\cdot N)=Gl(1\cdot N)
$$
represents a slice in $\wtm_0(X,A)$. This proves the theorem.

\subsection{On gluing maps for lower strata}\label{section_12.4}
We generalize our results from 1-nodal case to general strata.
\begin{corollary}\label{corollary_12.4.1}
The gluing map
$$
Gl_{\strata,\strata'}: \mc L_{\strata,\strata,\epsilon}^0 \to
\M_{\strata'}(X,A).
$$
gives a coordinate chart for $\M_{\strata'}(X,A)$.
\end{corollary}
Similarly,
\begin{corollary}\label{corollary_12.4.2}
The isomorphism
$$
Gl_{\strata,\strata'}^\ast (\mc L_{\strata'})= \mc L_{\strata}
$$
is diffeomorphic.
\end{corollary}

Set $W= Gl_{\strata,\strata'}(\mc
L_{\strata,\strata',\epsilon}^0)$. Since $Gl_{\strata,\strata'}$
is  diffeomorphic, the gluing map
$$
Gl'_{\strata',\strata''}: \mc
L_{\strata',\strata'',\epsilon}|_W\to \M_{\strata''}(X,A)
$$
given by
$$
Gl'_{\strata',\strata''} =Gl_{\strata,\strata''}\circ
Gl_{\strata,\strata'}\inv.
$$
is a diffeomorphism.

Moreover,
\begin{corollary}\label{corollary_12.4.3}
$Gl'_{\strata',\strata''}$ is admissible, so it is
$C^0$-compatible with $Gl_{\strata',\strata''}$.
\end{corollary}
{\bf Proof. }$Gl'_{\strata',\strata''}$ is admissible by its
construction and  definitions. The second assertion follows from
\S\ref{section_11}. q.e.d.

\section{Smooth structures on $\om_{g,m}(X,A)$}\label{section_13}

\subsection{Topology on $\om_{g,m}(X,A)$}\label{section_13.1}
By far, $\om_{g,m}(X,A)$ is a union of strata, each of which is a
smooth orbifold. We have not defined the topology on the whole
set. This is provided by gluing maps.

Recall that for any $(u,\mkj)\in \M_\strata(X,A)$, there exists a
neighborhood $U\subset \M_\strata(X,A)$ of $(u,\mkj)$ and
$\epsilon$ such that the gluing map
$$
Gl_\strata: \mc L_{\strata,\epsilon}|_{U} \to \om_{g,m}(X,A)
$$
exists. We define the image of $Gl_\strata$ to be a neighborhood
of $(u,\mkj)\in \om_{g,m}(X,A)$. By this way, we may define a
topological base at $(u,\mkj)$: to see we form a topological base,
we use the property of $C^0$-compatibility between gluing maps,
which says that any two such open sets are compatible. Therefore,
we have a topology on $\om_{g,m}(X,A)$. In fact, we have
\begin{theorem}\label{theorem_13.1.1}
$\om_{g,m}(X,A)$ is a topological orbifold.
\end{theorem}
{\bf Proof. }For each point $(u,\mkj)\in \M_{\strata}(X,A)$, a
neighborhood described above has a coordinate chart:
$$
(\mc L_{\strata,\epsilon}|_{U}, Gl_\strata).
$$
The transition maps between any two charts are $C^0$. Hence it is
an orbifold. q.e.d.

\subsection{Smooth structures on $\om_{g,m}(X,A)$ }\label{section_13.2}
In this subsection, we explain that there exists an atlas such
that $\om_{g,m}(X,A)$ is smooth. However, we do not show any two
atlas are compatible.

\begin{defn}\label{defn_13.2.1}
A stratum-covering of $\om_{g,m}(X,A)$ consists of
$U_\strata,\epsilon_\strata$ for each stratum such that
\begin{itemize}
\item $U_\strata$ is a proper subset of $\M_\strata(X,A)$; \item
there exists a $Gl_\strata$  on $\mc
L_{\strata,\epsilon_\strata}|_{U_\strata}$; \item for
$$
W_\strata= Gl_\strata(\mc
L_{\strata,\epsilon_\strata}|_{U_\strata}),
$$
$W_\strata\cap W_{\strata'}\not=\emptyset$ if and only if
$\strata\prec \strata'$ (or, $\strata'\prec \strata$); \item
$\{W_\strata\}_{\strata\in \mc D_{g,m}^A}$ is a covering of
$\om_{g,m}(X,A)$;
\end{itemize}
\end{defn}
The following lemma shows that stratum-coverings are abundance.
\begin{lemma}\label{lemma_13.2.1}
There are many stratum-coverings.
\end{lemma}
{\bf Proof. } Set $\mc D=\mc D^A_{g,m}$ Let $\mc S_0$ be the set
of smallest strata $\strata\in \mc D$. Choose
$$
U_\strata= \M_\strata(X,A).
$$
They are compact. By the gluing theory, there exists
$\epsilon_\strata$ such that the gluing map exists on $\mc
L_{\strata,\epsilon_\strata}$. If we choose $\epsilon_\strata$
small, we may have
$$
W_\strata\cap W_\strata'=\emptyset.
$$
Inductively, let $\mc S_k$ be the set of smallest strata
$\strata\in \mc D-\mc S_{k-1}$.

Suppose that $U_\strata,\epsilon_\strata$ are chosen for all
$\strata\in \mc S_{l},l\leq k-1$. Set
$$
W_{\strata,\strata'}=Gl(\mc
L_{\strata,\strata',\epsilon_\strata}|_{ U_\strata}).
$$
For any $\strata\in \mc S_k$ we choose a proper open set
$U_\strata$ such that
$$
\{W_{\strata',\strata}|\strata'\prec \strata\}\cup \{U_\strata\}
$$
covers $\M_\strata(X,A)$.  Moreover, we choose $\epsilon_\strata$
such that there exists a gluing map $Gl_\strata$ defined on $\mc
L_{\strata,\epsilon_\strata}|_{U_\strata}$ and $W_\strata$ is
disjoint with other $W_{\strata'}$ unless $\strata'\prec \strata$.
Inductively, this construct a stratum-covering. Since we are free
to choose $\epsilon_\strata$, $U_\strata$ (except $\strata\in \mc
S_0$) and $Gl_\strata$, hence there are many choices of
stratum-coverings. q.e.d.

\v\n Note that for  a given stratum-covering, we have  an atlas on
$\om_{g,m}(X,A)$ given by
$$
(\mc L_{\strata,\epsilon_\strata}|_{U_\strata}, Gl_\strata).
$$
Given such an atlas, we ask if the transition maps between any two
charts
$$
Gl_\strata\circ Gl_{\strata'}\inv
$$
are smooth. If so, we have shown the smoothness of
$\om_{g,m}(X,A)$. However, this may be too tedious and not true.
Instead, we show that there exists certain $Gl_\strata$ for each
$\strata$ such that
$$
Gl_\strata\circ Gl_{\strata'}\inv
$$
are smooth for any pair $(\strata\prec\strata')$.

The main idea is given by the following. Let $\strata_i\in \mc
D,i=1,2,$ with  $\strata_1\prec \strata_2$. Let $U_{\strata_i}$ be
proper open subsets of $\M_{\strata_i},i=1,2$. Suppose that  we
have a gluing map
$$
Gl_{\strata_1}: \mc L_{\strata_1,\epsilon_1}|_{U_{\strata_1}} \to
\om_{g,m}(X,A).
$$
Set
\begin{eqnarray*}
W_{\strata_1,\strata_2} &=&Gl_{\strata_1}(\mc
L_{\strata_1,\strata_2,\epsilon}
|_{U_{\strata_1}})\cap U_{\strata_2}; \\
W_{\strata_1,\strata_2} '&=&Gl_{\strata_1}(\mc
L_{\strata_1,\strata_2,0.5\epsilon}
|_{U_{\strata_1}})\cap U_{\strata_2};\\
W_{\strata_1,\strata_2}'' &=&Gl_{\strata_1}(\mc
L_{\strata_1,\strata_2,0.75\epsilon} |_{U_{\strata_1}})\cap
U_{\strata_2};
\end{eqnarray*}
As explained in \S\ref{section_9.4}, $Gl_{\strata_1}$ induces a
gluing map $Gl'_{\strata_2}$ on $W_{\strata_1,\strata_2}$. We show
that
\begin{prop}\label{prop_13.2.1}
There exists $\epsilon_2$ and gluing map
$$
Gl_{\strata_2}: \mc L_{\strata_2,\epsilon_2}|_{U_{\strata_2}} \to
\om_{g,m}(X,A)
$$
such that $Gl_{\strata_2}=Gl_{\strata_2}'$ on
$$
\mc L_{\strata_2,\epsilon_2}|_{W_{\strata_1,\strata_2}'}
$$
\end{prop}
{\bf Proof. }By the gluing theory, there exist $\epsilon$ and a
type-1 gluing map
$$
Gl_{\strata_2}'': \mc L_{\strata_2,\epsilon}|_{U_{\strata_2}} \to
\om_{g,m}(X,A).
$$
Note that $Gl_{\strata_2}'$ is of type-2. We use a cut-off
function on gluing parameter to patch these two gluing maps. To be
precise, let us introduce coordinates: by local coordinates, a
point in $W_{\strata_1,\strata_2}$ is denoted by
$$
Gl(u,\mkj, \rho), (u,\mkj,\rho)\in \mc L_{\strata_1,\strata_2}.
$$
Let $\beta$ be a cut-off function such that
$$
\beta(t)=\left\{
\begin{array}{ll}
1, & t\leq 0.5\epsilon \\
0, & t\geq 0.75\epsilon.
\end{array}\right.
$$
For an admissible gluing, we start with a  coordinate data
$(V,\Phi,F)$. Suppose this is the data used for $Gl_{\strata_2}'$.
Namely,
$$
V=pgl_{\strata_1,\strata_2}(\mc
L_{\strata_1,\strata_2,\epsilon_1})
$$
and
$$
F: V\to \M_{\strata_2}(X,A)
$$
realizes the gluing map. In terms of formula, it says
$$
Gl_{\strata_1,\strata_2}(u,\mkj,\rho)
=pgl_{\strata_1,\strata_2}(u,\mkj,\rho) +
\pert(pgl_{\strata_1,\strata_2}(u,\mkj,\rho)).
$$
Now we define $V'$ to be
$$
V'=\{pgl_{\strata_1,\strata_2}(u,\mkj,\rho) +
\beta(\rho)\pert(pgl_{\strata_1,\strata_2}(u,\mkj,\rho))
 \}.
$$
Start with $(V,\mc Q)$, it is easy to generate a new pair $(V',\mc
Q')$. Therefore, we define a new admissible gluing map
$Gl_{\strata_2}$ based on this coordinate data. Since the part of
$V'$ is in $\M_{\strata_2}(X,A)$ when $|\rho| \geq 0.75\epsilon$,
we may extend $Gl_{\strata_2}$ over $U_{\strata_2}$. q.e.d.

\v\n We remark that the cut-off function used to patch two gluing
maps is a function on $L_S$. We call the method to be patching
gluing maps.

\begin{theorem}\label{theorem_13.2.1}
There exists a stratum-covering $(U_\strata,\epsilon_\strata)$ and
gluing maps $Gl_{\strata}$ such that for any $\strata$
$Gl_{\strata}$ agrees with any gluing map $Gl'_{\strata}$ induced
from $Gl_{\strata'},\strata'\prec \strata,$ on the overlapping
domain.
\end{theorem}
{\bf Proof. } We use the same process as in lemma
\ref{lemma_13.2.1}. For $\strata\in \mc S_0$, no modification is
needed. Suppose the construction is done for all $\strata\in \mc
S_{l},l\leq k-1$.

Let $\strata\in \mc S_k$. For any $\strata'\prec \strata$, set
$W_{\strata', \strata}$ as before. Let
$$
W_{\strata}=\cup_{\strata'\prec \strata}W_{\strata',\strata}.
$$
For the moment, we denote $Gl_{\strata}(\strata')$ for the gluing
map defined over $W_{\strata', \strata}$ induced by
$Gl_{\strata'}$. We assert that
$$
Gl_{\strata}(\strata')=Gl_{\strata}(\strata'')
$$
over $W_{\strata',\strata}\cap W_{\strata'',\strata}$. First of
all, by the definition of stratum-covering, the intersection is
non-empty if and only if $\strata''\prec \strata'$. Since
$$
Gl_{\strata'}(\strata'')=Gl_{\strata'}
$$
over $W_{\strata'',\strata'}\cap U_{\strata'}$, hence they induce
same gluing maps on stratum $\M_{\strata}$. So totally, we have a
gluing map $Gl'_{\strata}$ over $W'_{\strata}$ induced by all
gluing maps from lower strata. For any gluing map $Gl''_{\strata}$
defined over $U_\strata$, we may apply proposition
\ref{prop_13.2.1} and get a new gluing map $Gl_{\strata}$ that is
a patching of $Gl_{\strata}'$ and $Gl_{\strata}''$. Then by
induction, we complete the construction. q.e.d.

\v\n As a corollary, we have
\begin{theorem}\label{theorem_13.2.2}
$\om_{g,m}(X,A)$ admits smooth structure.
\end{theorem}

\vskip 0.2in
\begin{center}
{\bf Part IV. Virtual theory on $\om_{g,m}(X,A)$}
\end{center}
\vskip 0.1in
In \cite{CT}, we introduce a new concept "virtual
manifolds/orbifolds". Furthermore, we develop the integration
theory  on them, which including the equivariant integration and
localization formulae. The background of the concept is to define
invariants on the moduli
spaces from Fredholm systems. In this part, our goal is to construct
a (smooth) virtual orbifold from $\om_{g,m}(X,A)$. Then all the theory
on virtual orbifolds can be applied to this particular moduli space. Therefore
the virtual localization formulae of Gromov-Witten invariants follow.

In \S\ref{section_14}, we review the material of virtual
orbifolds in \cite{CT}. Then we construct the virtual orbifold
structure on $\om_{g,m}(X,A)$ in \S\ref{section_19}-\S\ref{section_21}.
An application is given in \S\ref{section_22}.

\section{Virtual orbifolds}\label{section_14}

\subsection{Basic concepts}
Let $N=\{1,\ldots,n\}$ and $\mc N=2^N$ be the set of all subsets
of $N$. Let
$$
\mc X= \{X_I|I\in \mc N\}
$$
be a collection of sets indexed by $\mc N$.
For any $I\subset J$ there exist $ X_{I,J}\subset X_I,
X_{J,I}\subset X_J $ and a surjective map
$$
\phi_{J,I}: X_{J,I}\to X_{I,J}.
$$
Set ${\Phi}=\{\phi_{J,I}|I\subset J\}$. We always assume that $X_\emptyset \not=\emptyset$.
\begin{defn}\label{defn_2.2.1}
A pair $(\mc X,{\Phi})$ is called {\em patchable} if for any
$I,J\in \mc N$ we have
\begin{itemize}
\item[P1.] $X_{I\cup J,I\cap J}= X_{I\cup J,I}\cap X_{I\cup J,J}$;
\item[P2.] $X_{I\cap J,I\cup J}= X_{I\cap J,I}\cap X_{I\cap J,J}$;
\item[P3.] $\phi_{I\cup J,I\cap J}= \phi_{I,I\cap J}\circ
\phi_{I\cup J,I} =\phi_{J,I\cap J}\circ \phi_{I\cup J,J}$;
\item[P4.] $\phi_{I\cup J,I}(X_{I\cup J,I\cap
J})=\phi\inv_{I,I\cap J}(X_{I\cap J,I\cup J})$; \item[P5.]
$\phi_{I\cup J,J}(X_{I\cup J,I\cap J})=\phi\inv_{J,I\cap
J}(X_{I\cap J,I\cup J})$.
\end{itemize}
Set
\begin{eqnarray*}
X_{I,J}&=&\phi_{I\cup J,I}(X_{I\cup J,I\cap J})=\phi\inv_{I,I\cap J}(X_{I\cap J,I\cup J}),\\
X_{J,I}&=&\phi_{I\cup J,J}(X_{I\cup J,I\cap J})=\phi\inv_{J,I\cap
J}(X_{I\cap J,I\cup J}).
\end{eqnarray*}
\end{defn}
There is an equivalence relation for points in $\cup X_I$:
For $x\in X_I$ and $y\in X_J$ we say that $x\sim y$ if and only if
there exists a $K\subseteq I\cap J$ such that
$$
\phi_{I,K}(x)=\phi_{J,K}(y).
$$
We "patch" $X_I$ together and get a set
$$
\mathbf{X}= \bigcup_{I\in \mc N}X_I/\sim.
$$

A virtual manifold is a patchable pair $(\mc X,\Phi)$  with
specified properties.
\begin{defn}\label{defn_2.3.1}
Let $(\mc X,\Phi)$ be a patchable pair. Suppose that
\begin{itemize}
\item $X_I\in \mc X$ are  smooth orbifolds; \item $X_{I,J}$ and
$X_{J,I}$ are  {\em open} suborbifolds in  $X_I$ and $X_J$
respectively; \item $\Phi_{J,I}: X_{J,I}\to X_{I,J}$ is an
orbifold {\em vector bundle}.
\end{itemize}
Then $(\mc X,\Phi)$ is called a {\em virtual orbifold} if for any
$I$ and $J$,
\begin{eqnarray*}
&&\phi_{I,I\cap J}: X_{I,J}\to X_{I\cap J,I\cup J},\\
&&\phi_{J,I\cap J}: X_{J,I}\to X_{I\cap J,I\cup J}
\end{eqnarray*}
are orbifold vector bundles and
\begin{equation}\label{eqn_2.3.1}
X_{I\cup J,I\cap J}= X_{I,J}\times_{X_{I\cap J,I\cup J}} X_{J,I}.
\end{equation}
We call
$$
{\mathbf X}=\bigcup_{I\in \mc N} X_I/\sim
$$
 the {\em virtual space} of $(\mc X,\Phi)$.
We denote the projection map $X_I\to {\mathbf X}$ by $\phi_I$.

Let $d_I$ be the dimension of $X_I$. We call $d_\emptyset $ the
virtual dimension of $(\mc X,\Phi)$.
\end{defn}
One can also define virtual manifolds/orbifolds with boundary. From
now on, for simplicity, we forget the orbifold singularities and focus on
manifolds only.

The following example gives a typical method to construct  virtual manifolds.
\begin{example}\label{ex_2.3.1}
 Let $X$ be a manifold. Let $\{U_0,U_1, \ldots,U_n\}$ be an open
cover of $X$. Let $U_i^\circ = \frac{3U_i}{4}, i\geq 1$. Here
$\frac{3U_i}{4}$ just means an open subset whose closure is in
$U_i$. We use $\frac{3}{4}$ to make the notations more suggestive.

Let $N=\{1,\ldots, n\}$ and $I,J,K$ be as before. Define
\begin{eqnarray*}
X_\emptyset &=& U_0 - \bigcup_{i=1}^n U_i^\circ \\
X_I &=& \bigcap_{i\in I} U_i - \bigcup_{j\not\in I}U_j^\circ.
\end{eqnarray*}
Let $\mc X= \{X_I|I\in \mc N\}$. Define
$$
X_{I,J}=X_{J,I}= X_I\cap X_J.
$$
All possible $\psi_{J,I}$ are taken to be identities and let
$\mathbf \Phi=\{\phi_{J,I}\}$. Then $(\mc X,\mathbf\Phi)$ is a
virtual manifold (cf. Proposition \ref{prop_2.3.1}). Moreover, the
virtual space $\mathbf X$ is $X$.
\end{example}

We can define differential forms on virtual manifolds. There are two types.
The first type is nature.
Let $(\mc X,\Phi)$ be a virtual manifold.
\begin{defn}\label{defn_3.1.1}
A $k$-form on $(\mc X,\Phi)$ is
$$
\alpha=\{\alpha_I\in \Omega^k(X_I)|I\in \mc N\}
$$
such that
$$
\alpha_J=\phi^\ast_{J,I}\alpha_I
$$
on $X_{J,I}$.
\end{defn}
This is called a pre-$k$-form in \cite{CT}. It, in fact, induces a
$k$-form on the virtual manifold in the sense of \cite{CT}.

In order to consider the second type of forms, we need
Thom forms $\Theta_{J,I}$   of the bundle
$\Psi_{J,I}: X_{J,I}\to X_{I,J}$. To avoid the unnecessary
complication caused by the degree of forms, we always assume that
the degree of $\Theta_{J,I}$ is even.
\begin{defn}\label{defn_3.1.2}
A set of  forms $\Theta= \{\Theta_{J,I}\}_{I\subseteq J} $ is
called a {\em transition data} of $\mc X$ if it satisfies the
following compatibilities: for any $I$ and $J$,
$$
\Theta_{I\cup J, I\cap J} = \Psi^*_{I\cup J,I}\Theta_{I,I\cap
J}\wedge \Psi^*_{I\cup J,J} \Theta_{J,I\cap J}
$$
on $X_{I\cup J,I\cap J}$.
\end{defn}
\begin{defn}\label{defn_3.1.3}
A virtual form on $(\mc X,\Phi)$ is
$$
\mfk z=\{z_I\in \Omega^\ast(X_I)|I\in \mc N\}
$$
such that
$$
z_J=\phi^\ast_{J,I}z_I\wedge \Theta_{J,I}
$$
on $X_{J,I}$ for some transition data $\Theta$. $\mfk z$ is called
a $\Theta$-form on $\mc X$.
\end{defn}
For either forms or virtual forms, one can define close and compact
supported forms. Let $\mfk z$ be a compact supported $\Theta$-form, one can
define integration
$$
\int_\mathbf{X} \mfk z.
$$
The Stokes' theorem holds for this type integration.

The discussion given above can be generalized to
the equivariant case. Let $G$ be a compact Lie group.
\begin{defn}\label{defn_4.1.1}
By a {\em $G$-virtual manifold} ($\mc{X},\Phi)$,  we mean that
(a.) $(\mc X,\Phi)$ is a virtual manifold, (b.) each $X_I$ is
$G$-manifold and (c.) $\Psi_{J,I}: X_{J,I}\to X_{I,J}$ are
$G$-equivariant bundles for any $I\subset J$.
\end{defn}
To study the $G$-equivariant integration theory on $\mc X$, we
may consider  {\it $G$-equivariant transition data}
$\Theta_G=\{\Theta^G_{J,I}\}_{I\subseteq J}$. Then similarly, we
may define:
 $G$-equivariant
forms, $G$-equivariant $\Theta_G$ forms, and etc. For
a compact supported $\Theta_G$ form $\zeta=(\zeta_I)$, we can define
$$
\int^G_X \zeta.
$$
The virtual localization formula is stated as
\begin{theorem}\label{thm_4.3.1}
Let $\mc X$ be a finite dimensional virtual manifold with $G=S^1$
action. Let $X$ be its virtual space. Let $\zeta\in
\Omega_{\Theta_G,c} (\mc X^\circ)$ and $\alpha\in
\Omega^*_G(\mc X)$, then
$$
\mu_{\zeta}(\alpha)=\int_{X^G}\frac{ i_{X^G}^*(\alpha\wedge
\zeta)}{ e_{G}(X^G)}.
$$
\end{theorem}
We explain the notations. $\zeta$ is a compact supported $\Theta_G$
forms in the interior of $\mc X$;
$$
\mu_{\zeta}(\alpha):=\int_{X}^G \zeta\wedge \alpha;
$$
$X^G$ is the fix locus of the action, which itself is a virtual manifold;
and $e_G(X^G)$ is the  $G$-equivariant Euler class of
the virtual normal bundle of $X^G$ in $X$.

\subsection{From Fredholm systems to virtual manifolds}
We start with the following set-up.
\begin{defn}\label{defn_5.1.1}
A Fredholm system consists of following data:
\begin{enumerate}
\item[(B1)] let $\pi:\mcf\to \mcb$ be a Banach orbifold bundle
over a Banach orbifold $\mcb$; \item[(B2)] let $S: \mcb\to \mcf$
be a proper smooth section.
 In particular, the
properness implies that $M=S^{-1}(0)$ is compact; \item[(B3)] for
any $x\in M$, let $L_x$ be the linearlization of $S$ at $x$
$$
L_x: T_x\mcb\to \mcf_x.
$$
We assume that $L_x$ is a  Fredholm operator. Let $d$ be the index
of the operator.
\end{enumerate}
We refer the triple $(\mcb,\mcf,S)$ as a {\em Fredholm system}.
$M$ is called  the {\em moduli space} of the system.
\end{defn}
A core topic in studying moduli problems is to define  invariants
on such a system. This is  based on the study of $M$. It is well
known that if $L_x$ is surjective for all $x\in M$, $M$ is a
compact smooth orbifold. Then $M$ can be thought as a cycle in
$H_d(\mcb)$ representing the Euler class of  bundle $\mcf\to\mcb$.
Let $a\in H^d(\mcb,\mathbb{R})$, define
$$
\Phi(a)= \int_M a.
$$
The challenging problem is  to define invariants when  the
surjectivity of $L_x$
fails. The virtual technique is introduced
to deal with this  situation.
There are several different versions
of this  technique, however the main idea is the \s, which has
become popular since 60's. Our method follows \cite{R} closely.

We recall stabilization for a Fredholm system. Let $U$ be an open subset of $\mcb$, let
\def \mfo{\mathfrak{o}}
\def \mco{\mathcal{O}}
$$
\mfo: \mco_U \to U
$$
be a rank-$k$ vector bundle, let
$$
s: \mco_U\to \mcf_U
$$
be a  bundle map. Define a map
$$
\hat{S}: \mco_U\to \mcf_U; \hat{S}(u,o)= (u, S(u) + s(o)),
$$
where the expression is given in the form of local coordinates and
$S(u)+s(o)$ is the sum on fibers. By abusing the notations, we
usually use $S+s$ for $\hat{S}$ to emphasis that $S$ is stabilized
by $s$.

Let $\hat{L}_{(u,o)}$ be the linearization of $\hat{S}$ as a map
$$
\hat{L}_{(u,o)}: T_{(u,o)}\mco_U\to \mcf_u.
$$
We say that the pair $(\mco_U,s)$ {\it stabilizes} the system
$(\mcb,\mcf,S)$ at $U$ if $\hat{L}_{(u,o)}$ are surjective for all
$(u,o)\in \mco_U$. Set
$$
V_U = \hat{S}^{-1}(0)\subseteq \mco_U.
$$
This is now a smooth manifold of dimension $d+k$. Clearly, $M\cap
U\subseteq V_U$ and
$$
(u,o)\in M \iff o=0.
$$
We now explain  the existence of local \s s.

Suppose $L_x$ is not surjective for some $x\in M$. Let $O^x$ be a
finite dimensional subspace of $\mcf_x$ such that
\def \Im{\mathrm{Image}}
$$
\Im(L_x) + O^x = \mcf_x.
$$
For example, we may take $O^x$ to be the "cokernel" of $L_x$.

Let $U^x$ be a neighborhood of $x$ in $\mcb$. In order to make
notations more suggestive, we assume that $U^x=B_r(x)$ is the
radius-$r$ disk centered at  $x$ and $cU^x=B_{cr}(x)$ for $c\in
\mathbb{R}^+$.

We can restate this construction by using the concept of Fredholm
system. Let $\mfo^*\mcf\to \mco_U$ be the pull-back  bundle over
$\mco_U$. $\hat{S}$  then gives a canonical section of this bundle
in an obvious way. For simplicity, we still denote the section by
$\hat{S}$. Therefore, we have a  Fredholm system $(\mco_U,
\mfo^*\mcf, \hat{S})$. If $(\mco_U,s)$ stabilzes the system at
$U$, we say that $(\mco_U, \mfo^*\mcf, \hat{S})$ {\it stabilizes}
 $(\mcb,\mcf,S)$
at $U$. $V_U\subseteq \mco_U$ is the moduli space of  the new
system.

We may construct a canonical bundle $\mfo^*\mco_U\to V_U$, then
there is a canonical section $\sigma: V_U\to \mfo^*\mco_U$ given
by $(u,o)\to (u,o,o)$ with respect to the local coordinates. Then
$M\cap U= \sigma^{-1}(0)$. This reduces the infinite dimensional
system $(U,\mcf_U,S)$ to a finite dimensional system $(V_U,
\mfo^*\mco_U, \sigma)$. We call $(V_U, \mfo^*\mco_U, \sigma)$, or
simply $V_U$,  to be the {\it virtual neighborhood} of $M$ at $U$.
Bundles $\mco_U$ and $\mfo^*\mco_U$ are called the {\em
obstruction bundles}.

Suppose that
 $\mcf_{U^x}$ is trivialized as $\mcf_{U^x}= U^x\times \mcf_x$. We now
describe the \s\ using the notations given above by setting
$U=U^x$:
\begin{itemize}
\item[(C1)] the obstruction bundle is
$$
\mco_{U^x}= U^x\times O^x;
$$
\item[(C2')] the bundle map $s=I^x: \mco_{U^x}\to \mcf_{U^x}$ is
the standard embedding via the trivialization of $\mcf_{U^x}$
given above.
\end{itemize}

We may assume that the pair  $(\mco_{U^x}, I^x)$ stabilizes the
system at $U^x$ if $U^x$ is chosen small. This explains the
existence of local \s.
The trivialization of $\mcf_{U^x}$ prevents us to extend the
construction outside $U^x$. This is "taken care" by modifying the
bundle map $s$ as the following. Let $\eta^x$ be a cut-off
function on $U^x$ such that $\eta^x= 1$ in $\frac{U^x}{2}$ and
$=0$ outside $\frac{3U^x}{4}$. (C2') is then replaced by \v\n
\begin{itemize}
\item[(C2)] the bundle map is given by $s^x= \eta^xI^x.$
\end{itemize}
Clearly, $(\mco_{U^x}, s^x)$ stabilizes the system at
 $\frac{U^x}{2}$.
In this paper, we always use (C2) to construct  virtual
neighborhoods. It turns out that (C2) is the key towards the
construction of virtual orbifolds from a Fredholm system.

Since $M$ is compact by our assumption, there exists finite points
$\{x_i \}_{i=1}^n$ in $M$ such that
$$
M\subseteq \bigcup_{i=1}^n\frac{1}{2}U^{x_i}=:U,
$$
where $U^{x_i}$ are as above.

For simplicity, we set
$$
U_i= U^{x_i}, \mco_{i}=\mco_{U^{x_i}}, s_i=s^{x_i}.
$$
 We call the data $\{
(U_i, \mc O_i, s_i) \}$ a {\em local \s\ system of $U$}.
From such a local \s\ system, one is able to construct a virtual
manifold and other data that yield integrations on it. This is stated
as
\begin{prop}\label{prop_6.2.1}
Let $(\mcb,\mcf,S)$ be a Fredholm system.
\begin{enumerate}
\item
there exists a local \s\ system $\{U_i,s_i,\mco_i\}$.
\item
Let $\mc{X}$ be the natural virtual manifold for $\mcb$ generated
by the covering $\{U_i\}$. Using the \s\  data given above, one is
able to define a virtual manifold $\mc{W}=\{W_I\}$,
where $(W_I, \mfo^*_I\mco_I,\sigma)$ is a virtual neighborhood
over $U_I$.  Let $\mathbf W$ be the virtual space of $\mc{W}$.
\item $\mc O$ is a virtual bundle over $\mc W$. $\sigma$ is a section
of the bundle;
\item 
Let $\Theta_i$ be Thom form of $\mco_i$. All Thom forms $\Theta_I$
of $\mco_I$ restricting on $W_I$ form a $\Theta$-form. Denote the
form by $\theta$. If the moduli space $M$ is compact, $\theta\in
\Omega_{\Theta,c}(\mc{W})$. $\theta$ is an Euler class of $\mc O$.
\item For any $a\in \Omega^*(\mcb)$, let $ a_I= \pi_I^*a $ on
$W_I$. Then $(a_I)_{I\subseteq N}\in \Omega^*(\mc{V})$. To abuse
the notations, we still denote the form by $a$.
\end{enumerate}
\end{prop}

By the proposition, we have $ \mu_{\theta}(a). $ Also we know that
this is well defined not only on $\Omega^*(\mcb)$, but also on
$H^*(\mcb)$. If a global \s\ as in \S\ref{sect_5.2} exists, it is
easy to see that
$$
\Phi(a) = \mu_{ \theta}(a).
$$
This leads to the following definition.
\begin{defn}\label{defn_6.2.1}
Let $(\mcb,\mcf, S)$ be a Fredholm system. Let
$\{(U_i,\mco_i,s_i)\}$ be a local \s\  system  constructed in
Proposition \ref{prop_6.2.1}.
Let $\mc{W},\theta$ be the virtual manifold and obstruction form
given  above. For $a\in H^*(\mcb)$, define the invariants
$\Phi(a)$ to be $\mu_{\theta}(a)$.
\end{defn}
One can prove that the invariants is independent of the choice of
local \s\ systems.

One can further assume that the Fredholm system admits an $S^1$-action.
Then  we can construct a $G$-virtual manifold
$\mc{V}$ from a local $G$-\s \ system. Then we replace $\Theta_i$
by equivariant Thom forms $\Theta_i^G$. So we have $\Theta_I^G$'s
and $\Theta_{J,I}^G$'s. Clearly, $\theta_G=\{\Theta_I^G\}_I$ is a
$\Theta^G=\{\Theta_{J,I}^G\}$ form. For any $\alpha\in
\Omega^*_G(\mcb)$, define
$$
\Phi_G(\alpha)= \mu_{\mc{V}, \Theta_G}(\alpha).
$$

Now we can state the virtual localization formula for Fredholm
systems. Again, let $G=S^1$. We consider the Fredholm system
$(\mcb,\mcf,S)$ with $G$-action. Let $\mc V$ be the virtual
orbifold for the moduli space $M$. Let $V$ denote the virtual
space. Then $\mc V^G$ is the virtual orbifold for $M^G$ and its
virtual space is $V^G$.
 We have
\begin{theorem}\label{thm_6.4.1}
Let $(\mcb,\mcf,S)$ be an $S^1$-Fredholm system.
 For $\alpha\in \Omega_G^\ast(\mcb)$,
$$
\Phi_G(\alpha) = \int_{V^G} \frac{i^*_{V^G}\alpha\wedge \theta_G}
{e_{G}(V^G)}=\mu_{e_{\theta_G}(V^G)}(i^\ast_{V^G}\alpha).
$$
\end{theorem}

\section{Local Stabilizations}\label{section_19}

\subsection{Neighborhoods in $\bar\chi_{g,m}(X,A)$}
\label{section_19.1}

Let $u_o\in \bar\chi_{g,m}(X,A)$.
For simplicity, we will drop $(X,A)$ and write $\bar\chi_{g,m}$.
We describe neighborhoods
of $u_o$. Suppose that
$$
u_o\in \chi_{\mkj_o}\subset \chi_{\strata_o},
$$
where $\mkj_o\in M_{S_o}$. Within the stratum
$\chi_{\strata_o}$ the neighborhoods of $u_o$ is
well defined. Here we give an explicit construction
of neighborhoods which may be generalized to $\bar\chi_{g,m}$.
\v
Neighborhoods of $u_o$ within the stratum:
let $V$ be a neighborhood of $\mkj_o$ in $M_{S_o}$, there is a
trivialization of
$$
\chi_V:=\bigcup_{\mkj\in V}\chi_\mkj
$$
given by
$$
\phi: \chi_{\mkj_o}\times V\to \chi_V.
$$
Set $u_j=\phi(u,j)$.

We now consider two cases with respect to
whether  $S$ is stable or not. First,
suppose that $S$ is stable.
Let $U_{\mkj}(u_\mkj,\delta)$
be an $\delta$-neighborhood of $u_{\mkj}$
in $\chi_\mkj$. We define a neighborhood of $u_o$
in the stratum to be
$$
U_{\strata_o}(u_o,\delta, V)
=\bigcup_{\mkj\in V} U_{\mkj}(u_{\mkj},\delta).
$$
Now if $S$ is unstable. $u_{o}$ may have nontrivial
isotropic group $\aut(u_o)$. Set
$$
\Lambda_{\mkj_o}=\aut(u_o)\cdot u_o.
$$
We define a normal bundle of $\Lambda_{\mkj_o}$:
at $u_o$ we use $L^2$-norm to define a normal
tangent space $N_{u_o}$ that is normal to $\Lambda_{\mkj_o}$,
then define $N_u=\alpha\cdot N_{u_o}$ for
$u=\alpha\cdot u_o$. This automatically
define a normal bundle $N$ over $\Lambda_{\mkj_o}$
with fiber $N_{u}$. Then take a $\delta$-disk bundle $N_\delta$
of $N$ and use  $\exp$ mapping it to $\chi_{\mkj_o}$
to get a neighborhood
$$
U_{\mkj_o}(u_o,\delta)=\exp_{\Lambda_{\mkj_o}} N_\delta.
$$
Similarly, we do this for all $u_\mkj$ in $\chi_\mkj$.
We put them together and get $U_{\strata_o}(u_o,\delta, V)$.
The method provide here is standard to treat
nontrivial isotropic groups.
\v
Next we consider  "neighborhoods" of $u_o$ in
$\bar\chi_{g,m}$.
Recall that there is a gluing bundle $L_{S_o}$
over $M_{S_o}$ and
$$
\mc L_{\strata_o}\to \chi_{\strata_o}
$$
which is $\pi^\ast L_{S_o}$ via projection
$\pi: \chi_{\strata_o}\to M_{S_o}$.
Let $L_{S_o,\epsilon_o}$ be the $\epsilon_o$-disk bundle
of $L_{S_o}$. There is a gluing-surface map
$$
\gs: L_{S_o,\epsilon_o}\to \bar M_{g,m}.
$$
Set
$$
\mkj_{\rho}=\gs(\mkj,\rho).
$$
{\em Case 1.}
We now assume that
$\mkj_o$ is stable.
Then
$\gs$ is injective and
$\mkj_\rho$ is stable.
On $\mkj_\rho$ there is a map
$$
u_{\mkj_\rho}=\pgl(u_\mkj,\rho).
$$
We define a neighborhood of $u_o$ in $\bar\chi_{g,m}$ given
by
$$
U_{u_o}( \delta,\epsilon_o, V)
=\bigcup_{|\rho|<\epsilon_o}\bigcup_{\mkj\in V}
U_{\mkj_\rho}(u_{\mkj_\rho},\delta).
$$
{\em Case 2.} Suppose that $\mkj_o$ is not stable
but $\gs(\mkj_o,\rho)$ is stable. The typical example
is case (IIa) in \S\ref{section_12.2}. We take it as an example. To
avoid too much complication caused by notations.
We follow notations in \S\ref{section_12.2}.
As it is explained, $\gs$ is no longer injective.

The neighborhood of $\mkj_o$ in $M_{S_o}$ can be
parameterized by
$$
V'\times D_{y_o1},
$$
where $V'$ is a neighborhood of $\mkj_{o1}'$ in its
stratum, denoted by $M_{S_o'}$, and $D_{y_{o1}}$
is a neighborhood of $y_{o1}$ in $\Sigma_1$.
For $\mkj=(\mkj_1',y)$ we write
$u_{\mkj}$ to be $u_{\mkj_1',y}$.
Note that $\gs(\mkj,\rho)=\mkj_1'$.
Set
$$
u_{\mkj_1',y,\rho}=\pgl(u_{\mkj_1',y}, \rho).
$$
Hence on $\chi_{\mkj_1'}$ we get a slice
$$
\Lambda_{\mkj_1'}=\{u_{\mkj_1',y,\rho}|y\in D_{y_{o1}},\rho<\epsilon_o\}.
$$
We then give a $\delta$ neighborhood of this slice:
as usual we use $L^2$-norm to get its normal bundle
in $\chi_{\mkj_1'}$ and then using $\exp$ we map
a $\delta$-disk on $\chi_{\mkj_1'}$ to get a neighborhood.
{\em But} there is a tricky point:
the $L^2$-norm is induced from the metrics on $\Sigma_1$
and $X$, here we require  the metric on
$\Sigma_1$ varies as parameters $y,\rho$ vary.
In fact, for fixing $y$ and $\rho$ the metric we use is
the metric on the connected sum
$$
\Sigma_1\sharp_{y,\rho}S^2.
$$
We can arrange the metric varies smoothly with respect to
$y$ and $\rho$.
By this way, we get a $\delta$-neighborhood
$$
U_{\mkj_1'}(\Lambda_{\mkj_1'},\delta)
$$
of $\Lambda_{\mkj_1'}$.
Then the neighborhood of $u_o$
 is defined to be
 $$
U_{u_o}(\delta,\epsilon_o,V)= U_{\strata_o}(u_o,\delta,V)
\cup\bigcup_{\mkj_1'\in V'}U_{\mkj_1'}(\Lambda_{\mkj_1'},\delta).
 $$
\v\n
{\em Case 3. }We consider the case that $\mkj_o$
and $\gs(\mkj_o,\rho)$ are both un-stable.
The typical example is \S\ref{section_12.3}. The idea is a combination of
case 2 and case 1 with non-trivial isotropic groups.
We leave the construction to readers.

\subsection{Cut-off functions}\label{section_19.2}

On neighborhoods $U_{u_o}(\delta,\epsilon_o, V)$, we can construct
smooth cut-off functions easily: let $\beta_1$ be a cut-off
function such that
$\beta_1(t)=1, t\leq \delta/4$ and 0 when $t\geq \delta/2$;
let $\beta_2$ be a cut-off
function such that
$\beta_2(t)=1, t\leq \epsilon_o/4$ and 0 when $t\geq \epsilon_o/2$;
let $\beta_3$ be a cut-off function on $V$ which is supported in
$V/4$, then we set
a cut-off function $\beta_{u_o}$ on $U_{u_o}(\delta,\epsilon_o, V)$
as
$$
\beta_{u_o}(\exp_{u_{\mkj_\rho}}\zeta)
=\beta_1(\|\zeta\|)\beta_2(|\rho|)\beta_3(\mkj).
$$

\subsection{Obstruction bundles}\label{section_19.3}
For any $u\in \M_{\mkj}(X,A)\subset
\om_{g,m}(X,A)$ we let
$\mathrm{coker}_u$ be the cokernel of operator
$D_{u,\mkj}$. Choose a cut-off function $\beta$
on $\mkj$ such that it is support away from
nodal points and $\aut(\mkj)$-invariant. For a proper
choice of $\beta$, namely, if the support of $\beta-1$
is near nodal points, then
the space
$$
O_u:=\beta\cdot \mathrm{coker}_u
$$
is complement to the image of $D_{u,\mkj}$.
For any $u$ and its neighborhood $U_{u}(\delta,\epsilon
,V)$ we define the local obstruction bundle
$$
\mc O_u = U_u(\delta,\epsilon, V)\times O_u.
$$

\subsection{Local stabilization}\label{section_19.4}
We now can follow the argument in \S\ref{section_14}.
Be equipped
with $\mc O_{u_o}$ and cut-off functions
$\beta_{u_o}$ we
can construct the local stabilization at
$U_{u_o}(\delta,\epsilon, V)$
if $\delta,\epsilon$ and $V$
are small.

Be precise, for small $\delta,\epsilon$ and $V$,
we can embed  $O_{u_o}$ into $\mc E_u$ properly
for any
$u\in U_{u_o}(\delta,\epsilon, V)$. Then we define
the stabilized equation over $\mc O_{u_o}$
by
$$
\hat S_{u_o}(u,\xi)=
\bar\partial_{J,\mkj}u + \beta_{u_o}(u)\xi=:
(\bar\partial + s_{u_o})(u,\xi).
$$
This finishes the construction of local stabilization
for $\om_{g,m}(X,A)$.

Let $W_{u_o}$ be the moduli space $\hat S_{u_o}\inv(0)$. If
$d\hat S_{u_o}|_u$ is surjective for any $u$.
Then $W_{u_o}$ is a topological orbifold. The proof
is same as that in part II. In fact, we can parallelly copy
the argument in \S\ref{section_13} and show that
$W_{u_o}$ admits a smooth structure.

\section{Virtual  structures for $\om_{g,m}(X,A)$ and the
Gromov-Witten invariants}
\label{section_20}

\subsection{Virtual orbifold structures on $\om_{g,m}(X,A)$}
\label{section_20.1}
As explained in \S\ref{section_14},
the existence of local
stabilization and the compactness of $\om_{g,m}(X,A)$
imply that there is a virtual (topological) orbifold
for $\om_{g,m}(X,A)$. We formulate notations.

Suppose that there are $n$ points
$$
\Lambda=\{u_1,\ldots,u_n\}\subset \om_{g,m}(X,A)
$$
with neighborhoods $U_{u_i}(2\delta_i,2\epsilon_i,2V_i)$
such that
$$
\bigcup_{i=1}^n U_{u_i}(\delta_i,\epsilon_i,V_i)
\supset \om_{g,m}(X,A).
$$
Following the construction given in \S\ref{section_14}, we have
a sequence of orbifolds, (which may not be smooth,)
$$
\{W_I\}_{I\subset \{1,\ldots, n\}}.
$$
 Hence we have
a virtual orbifold $\mc W$ given by
$$
W_k=\{W_I\}_{|I|=k}.
$$
The goal is to show that
\begin{theorem}
$\mc W$ admits a smooth structure. Hence it can be a
smooth virtual orbifold.
\end{theorem}
{\bf Proof. } Let
$$
p: W_I\to \overline{\chi}_{g,m}(X,A)
$$
be the projection. For each $W_I$ set
$$
W_I(\strata)= p\inv(\chi_\strata)\cap W_I.
$$
Then
$$
\mc W_\strata= \{W_I(\strata)\}
$$
forms a smooth virtual orbifold for each $\strata$.
This is due to the construction in \S\ref{section_14} for
a Fredholm system.
Since we are working within a stratum, the smooth structure
exists automatically.

Next we show that $\mc W$ admits a smooth structure
at $\mc W_\strata$.
For each $W_I$ the smooth structure
at $W_I(\strata)$ is induced by gluing maps.
Let $\mc L_\strata$ be the gluing bundle over
$\chi_\strata$. It induces a bundle over each
$W_I(\strata)$, we denote the bundle by $\mc L_{I,\strata}$.
Then we have gluing maps
$$
Gl_{I,\strata}: \mc L_{I,\strata,\epsilon_o}
\to W_I.
$$
Note that
$$
\{\mc L_{I,\strata}\}
$$
itself is a smooth virtual orbifold. If $Gl_{I,\strata}$
is compatible with the overlapping maps, then
the smooth structures induced on $W_I$ are compatible
with the virtual structure on $W_I$. To be precise, this is what we
mean:
suppose we have $I\subset J$ and $x\in W_{I,J}, y\in W_{J,I}$ with
$x=\pi_{J,I}(y)$. We denote them by
$$
x=(u,o_1), y=(u,o_1,o_2).
$$
For any gluing parameter $\rho$ we want
\begin{equation}\label{eqn_20.1}
Gl_{J,\strata}(y,\rho)= Gl_{I,\strata}(x,o_2).
\end{equation}
To make \eqref{eqn_20.1} available, we should require
that the pre-gluing maps and right inverses
$Q_x, Q_y$ used for gluing map are same. There is no problem
for
the consistency of pre-gluing maps.
For  right inverses,
this can be easily
achieved as well: let $\mc Q_{I,\strata}$ be right inverses
used for $W_I(\strata)$, we can use partition of unity
to reproduce a new group of right inverses $\mc Q_{I,\strata}'$
such that for any $x$ and $y$ as above
$$
Q_x\in \mc Q'_{I,\strata}, Q_y\in \mc Q'_{J,\strata}
$$
are equal. This allows us to give a smooth structure
of $\mc W$ at $\mc W_\strata$.

As before, since the smooth structures on $\mc W$
induced by gluing maps from different strata
may be different, we should apply the technique given
in \S\ref{section_13}:
Let $\mc S_0$ be the set of smallest strata.  for any
$\strata\in \mc S_0$ let
$$
\mc W_\strata=\{W_I| p(W_I)\cap \chi_\strata\not=\emptyset
\}.
$$
It is a smooth virtual orbifold. We may assume that
$$
\mc W_\strata\cap\mc  W_{\strata'}=\emptyset.
$$
Hence,
$$
\mc W_{\mc S_0}:=\bigcup_{\strata\in \mc S_0} W_\strata
$$
still form a smooth virtual orbifold.

Next we consider the set $\mc S_1$ of smallest strata next to
those in $\mc S_0$. Then for $\strata\in \mc S_1$,
$\mc W_\strata$ is still a smooth virtual orbifold. However,
on $\mc W_\strata\cap \mc W_{\mc S_0}$ they may have two different
smooth structures due to the discrepancy of gluing
maps on different strata. We can then apply the
argument in \S\ref{section_13.2} to perturb the gluing maps
on $\mc W_\strata$ such that
its smooth structure is compatible with that induced
from $\mc W_{\mc S_0}$. By this way, we have a modified smooth
structure on
$$
\mc W_{\mc S_1}=\bigcup_{\strata\in \mc S_1}\mc  W_\strata
$$
such that
$$
\mc W_{\mc S_0\cup\mc S_1}=\mc W_{\mc S_0}\cup\mc  W_{\mc S_1}
$$
forms a smooth virtual orbifold. We continue the
process, then we have a smooth structure on $\mc W$.
q.e.d.

\subsection{The Gromov-Witten invariants}

For the moduli space $\om_{g,m}(X,A)$ we have constructed
an associated virtual orbifold $\mc W$.  As explained in
\S\ref{section_14}, we have a transition data on $\mc W$
$$
\Theta=\{
\Theta_{J,I}= \bigwedge_{j\in J-I}\Theta_j\}_{I\subset
J}
$$
and a $\Theta$-form $\theta=(\Theta_I)$.

Suppose that the virtual dimension of $\om_{g,m}(X,A)$
then for any degree $d$ form $\alpha$ on $\mc W$ we define
the Gromov-Witten invariants to be
$$
\mu_\theta(\alpha).
$$
In general, $\alpha$ is induced from forms on $X$ (by
evaluation maps) or from forms on $\bar M_{g,m}$. Moreover,
the invariant is independent of the construction
of $\mc W$.

\section{Symplectic virtual localization}\label{section_21}

We now derive the symplectic virtual localization formula
for Gromov-Witten invariants.

Let $G=S^1$ act on $(X,\omega)$ symplectomorphically.
It then induces an action on $\overline
\chi_{g,m}(X,A)$ and on $\om_{g,m}(X,A)$. First  we can modify the construction of virtual orbifold
$\mc W$ such that it is an $S^1$-virtual orbifold.
The forms $\Theta$, $\theta$ and
$\alpha$ are then replaced by equivariant
forms $\Theta^G$, $\theta_G$ and $\alpha_G$ .

Then applying the virtual
localization formula for $G$-virtual orbifolds, we have
\begin{theorem}
Suppose that $(X,\omega)$ admits an $S^1$ symplectomorphic
action. Then the virtual localization formula for
Gromov-Witten invariant $\mu_\theta(\alpha)$ is given by
$$
\mu_{\theta_G}(\alpha_G)=
\int_{\mc W^G}\frac{ i_{\mc W^G}^*(\alpha_G\wedge
\theta)}{ e_{G}(\mc W^G)}.
$$
Here $\mc W^G$ is the virtual orbifold for $\om_{g,m}^G(X,A)$
and $e_{G}(\mc W^G)$ is the equivariant Euler form
of the normal bundle of $\mc W^G$ in $\mc W$.
\end{theorem}

\section{An application of the virtual localization formula}
\label{section_22}

\subsection{Models $W_k$ and their Gromov-Witten invariants }
\label{sect_22.1}

Let
$$
V_k=\{(u_1,u_2,u_3,u_4)|u_1^2+u_2^2+u_3^2+ u_4^{2k}=0\}
$$
for $k=1,2,\ldots$. $V_k$ contains a singularity at 0.  By blowing-up
at $0$, we have $W_k$ with an exceptional line $A=\mathbb{P}^1$.
$W_k$ can be given by two coordinate patches $(w,z_1,z_2)$ and
$(x, y_1,y_2)$ and by a transition map between the coordinate
patches. Here $A$ is given by $\{z_1=z_2=0\}=\{y_1=y_2=0\}$. The transition
map is given by
$$
\left\{
\begin{array}{l}
z_1=x^2y_1+xy_2^k\\
z_2=y_2\\
w=1/x.
\end{array}
\right.
$$
The normal bundle of $A$ in $W_k$ is known as
$$
\left\{
\begin{array}{lll}
\mc O(-1)\oplus \mc O(-1), & \mbox{when} & k=1\\
\mc O\oplus \mc O(-2), && k\geq 2.
\end{array}
\right.
$$
$W_k$ are Calabi-Yau threefolds.
When $k=1$, this is well-known conifold.
Since $A$ is extremal ray, the moduli spaces
$$
\om_{g,0}(W_k, d[A])=\om_{g,0}(A,d[A]).
$$
Hence, we are allowed to define (local) Gromov-Witten invariants.
When $k=1$, the invariants on $W_1$ is computed by
Faber-Pandharipande in \cite{FP} by localization
techniques. When $k>2$, the invariants are computed by
 Bryan-Katz-Leung
\cite{BKL} by using deformation arguments. The results are given in the
following theorem.
\begin{theorem}[Faber-Pandharipande, Bryan-Katz-Leung]\label{thm_22.1}
Let $C_k(g,d)$ be the Gromov-Witten invariants for moduli spaces
$\om_{g,0}(W_k,d[A])$. Then
\begin{equation}\label{eqn_22.1}
C_1(g,d)= \frac{|B_{2g}|d^{2g-3}}{2g\cdot(2g-2)!}
\end{equation}
and
\begin{equation}\label{eqn_22.2}
C_k(g,d)=kC_1(g,d).
\end{equation}
\end{theorem}
In this paper, we use the localization formula to verify \eqref{eqn_22.2}.
Such a model is closely related to the framework of
Li-Ruan's study on Gromov-Witten theory with respect to flops.
Such a problem was first proposed and solved in
\cite{LR}, and then later reconsidered in \cite{LY}.
In \cite{LY}, a computation of \eqref{eqn_22.2} without
using deformation is also asked.  On the other hand,
in orbifolds, there is a similar problem in this framework.
It is known that the deformation technique can not be applied for orbifold case.
Partial results have been  considered in
\cite{CLZZ}. The localization technique would be a key
to understand the orbifold Gromov-Witten invariants. These
are the motivations for recompute \eqref{eqn_22.2} by using
localization.

\subsection{Localization set-up}
There is a $T^2$-action on $W_k$ given by
$$
(t_1,t_2)\cdot(w,z_1,z_2)=
(t_1^\lambda w, t_2^uz_2, t_1^{-\lambda}t_2^{ku}z_1).
$$
The weights of the action is said to be $(\lambda, u, -\lambda+ku)$.

The moduli space is $\om_{g,0}(A, d[A])$. The fix loci of the action
in this moduli space are associated with graphs(\cite{K}, \cite{GP}).
For each graph $\Gamma$, we denote the fixed loci by $M_\Gamma$.
The Gromov-Witten invariant is given by
\begin{equation}\label{eqn_22.3}
\int_{V_k}1=
\sum_{\Gamma}\int_{M_\Gamma}\frac{\theta}{e_{T}(N_{M_\Gamma})}=
\sum_{\Gamma}\int_{M_\Gamma}\frac{1}{e^{vir}_{T}(N_{M_\Gamma})},
\end{equation}
where $V_k$ is a virtual neighborhood of the moduli
space $\om_{g,0}(A,d[A])$ and $\theta$
is a $\Theta$-form constructed from cokernels. Unlike the well-known case $k=1$, neither $V_k$
is the moduli space, nor $V_k=V_{k'}$ when $k\not=k'$.
On $M_\Gamma$ we have
is  a $K$-bundle $\mathbb{H}^0-\mathbb{H}^1$. The fiber of $\mathbb{H}^i,
i=0,1$ over $f:\Sigma\to W_k$ is given by
$
H^i(\Sigma, f^\ast TW_k).
$ Then
$$
\theta|_{M_\Gamma}= e_T(\mathbb{H}^1), \;\;\; N_{M_\Gamma}=\mathbb{H}^0.
$$
\subsection{Proof of Theorem \ref{thm_22.1}}
We follow the computation in \cite{FP}. We denote the left hand side
of \eqref{eqn_22.3} by $J$, each term on right hand side by $I_{\Gamma}(\lambda,u)$
and the sum by
 $I(\lambda, u)$. Clearly,
$J=I(\lambda,u)$ implies that $I$ is independent of choice of $u$.
We will compute
$$
\lim_{u\to 0}I(1, u).
$$
By the same reason as the computation in \cite{FP}, we know that
$$
\lim_{u\to 0}I_{\Gamma}(1,u)=0
$$
unless $\Gamma$ is the graph that consists of a single edge. Now
let $\Gamma$ be such a graph that consists of one edge.
Its two ends are marked by $p_1$ and
$p_2$, the fixed point on $A$ of the action. ($p_1$ is the point with
$w=0$ and $p_2$ is the other one.) Suppose the corresponding genus
are $g_1,g_2$ with $g_1+g_2=g$. Such a graph is denoted  by $\Gamma_{g_1,g_2}$. Then by
a direct computation,
we find that
$$
\lim_{u\to 0}I_{\Gamma_{g_1,g_2}}(1,u)
= kd^{2g-3} b_{g_1}b_{g_2}.
$$
Therefore,
$$
\lim_{u\to 0}I(1,u)= kd^{2g-3}\sum_{g_1+g_2=g} b_{g_1}b_{g_2}
=k C(g,d).
$$
The last equation is proved in \cite{FP}. This proves the theorem.

\end{document}